\documentclass[a4paper,12pt,twoside]{article}
\usepackage{graphicx}
\usepackage{a4wide}
\usepackage[top=2.5cm,bottom=2.5cm,right=1cm,left=1cm]{geometry}
\usepackage{amssymb,amsmath,amsthm}
\usepackage{mathrsfs}
\usepackage{siunitx}
\usepackage{subfigure,color,colordvi,graphicx,xcolor}
\usepackage[only,llbracket,rrbracket,interleave]{stmaryrd}
\usepackage[superscript,sort,compress]{cite}

\makeatletter
\newcommand*\wt[2][0.2ex]{%
        \begingroup
        \mathchoice{\wt@helper{#1}{#2}{\displaystyle}{\textfont}}
                   {\wt@helper{#1}{#2}{\textstyle}{\textfont}}
                   {\wt@helper{#1}{#2}{\scriptstyle}{\scriptfont}}
                   {\wt@helper{#1}{#2}{\scriptscriptstyle}{\scriptscriptfont}}%
        \endgroup
        #2%
}
\newcommand*\wt@helper[4]{%
        \def\currentfont{\the#41}%
        \def\currentskewchar{\char\the\skewchar\currentfont}%
        \setbox\tw@\hbox{\currentfont#2\currentskewchar}%
        \dimen@ii\wd\tw@
        \setbox\tw@\hbox{\currentfont#2{}\currentskewchar}%
        \advance\dimen@ii-\wd\tw@
        \rlap{\raisebox{-#1}{$\m@th#3\kern\dimen@ii\widetilde{\phantom{#2}}$}}%
}
\makeatother

\newcommand{\bm}[1]{\text{\boldmath $#1$\unboldmath}}

\newcommand{\vect}[1]{\mathbf{#1}}
\newcommand{\mat}[1]{\mathbf{#1}}

\newcommand{\grad}{\bm{\nabla}}

\newcommand{\eltwo}{\ensuremath{\mathcal{L}_2}}

\newcommand{\Vone}{\mathcal{V}^1}
\newcommand{\Vzero}{\mathcal{V}^0}
\newcommand{\VzeroHat}{\mathcal{\hat{V}}^0}

\newcommand{\VoneElem}{\Vone(\Omega_e)}
\newcommand{\VzeroElem}{\Vzero(\Omega_e)}
\newcommand{\VzeroFaces}{\VzeroHat(\Gamma\cup\Gamma_N)}

\newcommand{\nen}  {\ensuremath{\texttt{n}_{\texttt{en}}}}
\newcommand{\nfn}  {\ensuremath{\texttt{n}_{\texttt{fn}}}}
\newcommand{\nsd}  {\ensuremath{\texttt{n}_{\texttt{sd}}}}
\newcommand{\numel}{\ensuremath{\texttt{n}_{\texttt{e}}}}

\newcommand{\numfael}{\ensuremath{\texttt{n}_{\texttt{f}}^{\texttt{e}}}}
\newcommand{\niadapt}  {\ensuremath{\texttt{n}_{\texttt{i}}}}

\newcommand{\nmo}{M}
\newcommand{\basis}{\widetilde{N}}

\newcommand{\hu}{\hat{u}}
\newcommand{\bu}{\bm{u}}

\newcommand{\bhu}{\widehat{\bu}}
\newcommand{\bq}{\bm{q}}

\newcommand{\bn}{\bm{n}}
\newcommand{\bx}{\bm{x}}
\newcommand{\bL}{\bm{L}}

\newcommand{\qe}{\vect{q}_e}
\newcommand{\ue}{\text{u}_e}

\newcommand{\uHj}{\hat{\text{u}}_j}

\newcommand{\Le}{\mat{L}_e}
\newcommand{\bue}{\vect{u}_e}
\newcommand{\pe}{\text{p}_e}

\newcommand{\buHj}{\hat{\vect{u}}_j}

\newcommand{\Insd}{\mat{I}_{\nsd}}

\newcommand{\jump}[1]{\llbracket #1\rrbracket}

\newcommand{\Pu}{\mathbb{P}_0}

\newcommand{\projvi}{\vect{p}_{e,i}}

\newcommand{\projcj}{p_{e,j}}
\newcommand{\projvj}{\vect{p}_{e,j}}

\newcommand{\projviT}{\tilde{\vect{p}}_{e,i}}
\newcommand{\projmiT}{\widetilde{\mat{P}}_{e,i}}
\newcommand{\projcjT}{\tilde{p}_{e,j}}
\newcommand{\projvjT}{\tilde{\vect{p}}_{e,j}}

\newcommand{\bw}{\bm{w}}

\newcommand{\Aset}{\mathcal{A}_e}
\newcommand{\Bset}{\mathcal{B}_e}
\newcommand{\Dset}{\mathcal{D}_e}
\newcommand{\Nset}{\mathcal{N}_e}

\newenvironment{keywords}{\begin{quote}\emph{\textbf{Keywords:}}}{\end{quote}}

\theoremstyle{definition}

\newtheorem{remark}{Remark}

\graphicspath{{figsPDF/}{figsPNG/}{figs_arXiv/}}

\begin{document}
\title{A second-order face-centred finite volume method on general meshes with automatic mesh adaptation}

\author{
\renewcommand{\thefootnote}{\arabic{footnote}}
			  M. Giacomini\footnotemark[1] \ and
			  R. Sevilla\footnotemark[2]
}

\date{}
\maketitle

\renewcommand{\thefootnote}{\arabic{footnote}}

\footnotetext[1]{Laboratori de C\`alcul Num\`eric (LaC\`aN), ETS de Ingenieros de Caminos, Canales y Puertos, Universitat Polit\`ecnica de Catalunya, Barcelona, Spain}
\footnotetext[2]{Zienkiewicz Centre for Computational Engineering, College of Engineering, Swansea University, Wales, UK
\vspace{5pt}\\
Corresponding author: Ruben Sevilla. \textit{E-mail:} \texttt{r.sevilla@swansea.ac.uk}
}

\begin{abstract}
A second-order face-centred finite volume strategy on general meshes is proposed. The method uses a mixed formulation in which a constant approximation of the unknown is computed on the faces of the mesh. Such information is then used to solve a set of problems, independent cell-by-cell, to retrieve the local values of the solution and its gradient. The main novelty of this approach is the definition of a new basis function, utilised for the linear approximation of the primal variable in each cell, suitable for computations on general meshes, including meshes with different element types. The resulting approach provides second-order accuracy for the solution and first-order for its gradient, without the need of reconstruction procedures, is robust in the incompressible limit and insensitive to cell distortion and stretching. The second-order accuracy of the solution is exploited to devise an automatic mesh adaptivity strategy. An efficient error indicator is obtained from the computation of one extra local problem, independent cell-by-cell, and is used to drive mesh adaptivity. Numerical examples illustrating the approximation properties of the method and of the mesh adaptivity procedure are presented. The potential of the proposed method with automatic mesh adaptation is demonstrated in the context of microfluidics.
\end{abstract}


\section{Introduction}
\label{sc:Intro}

Finite volume (FV) methods are one of the most popular computational methods for solving systems of conservation laws~\cite{MR1925043,EYMARD2000713,MR2417929,barth2018finite}. These methods are usually classified into two families, namely cell-centred FVs and vertex-centred FVs depending on the definition of the unknowns at the centroids or at the vertices of the cells respectively. The main attractive properties of FV methods are their numerical efficiency, local conservation and robustness which make them appealing solutions to treat flow problems of industrial interest~\cite{cueto2007finite,Othmer-14,Parolini-CACPSF-13}. However, one of the drawbacks of low-order cell-centred and vertex-centred FVs is the need for a reconstruction of the gradient. In this context, the quality of the reconstruction is directly linked to the quality of the mesh, leading to an important loss of accuracy, and even second-order convergence, in the presence of highly distorted or stretched cells~\cite{diskin2010comparison,diskin2011comparison}.

A new class of FV methods, named as face-centred finite volume (FCFV) method, was recently introduced in~\cite{FCFV2018}. The method is based on the hybridisable discontinuous Galerkin (HDG) method by Cockburn and co-workers \cite{cockburn2004characterization,Jay-CG:05,Jay-CG:05b,Jay-CG:05-GAMM,Jay-CGL:09} and defines the global unknowns on the cell faces (edges in two dimensions). The main attractive properties of this scheme are the ability to produce a first-order accurate approximation of the solution and its gradient without the need of a reconstruction. Therefore, the method is insensitive to mesh distortion and cell stretching. In addition, the FCFV method inherits the convergence properties of HDG and it passes the LBB condition using equal-order approximations for velocity and pressure in the context of incompressible flows~\cite{Jay-CGNPS-11,giacomini2018superconvergent} and it is robust when solving linear elasticity problems in the incompressible limit~\cite{Fu-FCS-15,HDGelasticityVoigt,FCFVelas}. The main drawback of this method is that, even for a sufficiently regular mesh, it provides a first-order approximation of the solution, compared to the second-order provided by cell-centred and vertex-centred FV methods.

In~\cite{FCFV2} the authors proposed a second-order FCFV method with a computational cost almost identical to the cost of the original first-order FCFV method. The main idea is to use a piecewise linear approximation of the solution in the cells, but maintain a piecewise constant approximation for its gradient in the cells and for the solution on the cell faces. In addition, the method introduces a projection operator in the definition of the numerical fluxes, following the work of~\cite{oikawa2015hybridized,oikawa2016analysis,Shi-QS-16}. However, the second-order FCFV method proposed in~\cite{FCFV2} is only applicable on simplicial meshes. Furthermore, the mesh adaptivity process proposed in~\cite{FCFV2} is expensive as it requires the solution of two global problems to compute a local error indicator to drive the mesh adaptivity process.

This paper proposes a new second-order FCFV method applicable to general and hybrid meshes. The key idea is to introduce a new approximation space for the primal variable that leads to second-order convergence on general meshes. Numerical examples involving triangular and quadrilateral cells in two dimensions and tetrahedral, hexahedral, prismatic and pyramidal cells in three dimensions are presented to demonstrate the optimal convergence properties of the method in the context of second-order elliptic problems. In addition, this paper proposes a new and efficient error indicator to drive a mesh adaptivity process. Contrary to the error indicator proposed in~\cite{FCFV2}, the new strategy does not require the solution of two global problems and only involves local quantities. More precisely, the proposed error indicator only requires cell-by-cell calculations and therefore can be easily computed in parallel.

The remainder of this paper is organised as follows. The rationale of the second-order FCFV method is recalled in section~\ref{sc:FCFV2} for a scalar second-order elliptic problem. After introducing the new basis functions required for general and hybrid meshes, the novel second-order FCFV formulation for Poisson and Stokes equations is derived in section~\ref{sc:FCFV2general}. Section~\ref{sc:meshAdaptivity} proposes an efficient strategy to perform mesh adaptivity via the computation of an inexpensive local error indicator. Extensive numerical tests are discussed in section~\ref{sc:studies} to validate the proposed method, in two and three dimensions, and to verify its optimal approximation properties for general and hybrid meshes. Section~\ref{sc:examples} presents the application of the proposed mesh adaptivity strategy to a two dimensional thermal problem using both triangular and quadrilateral meshes and to an incompressible Stokes flow around a complex three dimensional geometry of a microswimmer. Eventually, section~\ref{sc:Conclusion} summarises the main results and novelties of the paper.

\section{Fundamentals of the second-order FCFV method} 
\label{sc:FCFV2}

This section briefly recalls the second-order FCFV method introduced in~\cite{FCFV2}. To simplify the presentation, the method is described using the Poisson equation as a model problem.

\subsection{Problem statement and mixed formulation}
\label{sc:PoissonStatement}

An open bounded domain $\Omega\in\mathbb{R}^{\nsd}$ is considered, where $\nsd$ denotes the number of spatial dimensions. The boundary of the domain is split into the non-overlapping Dirichlet boundary, $\Gamma_D$, where the solution is known, and the Neumann boundary, $\Gamma_N$, where the normal flux is known. 

The computational domain is assumed to be partitioned in $\numel$ non-overlapping cells $\Omega_e$, for $e=1,\ldots,\numel$. The boundary of each cell is expressed as the union of a set of faces (edges in two dimensions),  $\Gamma_{e,j}$, for $j=1,\ldots,\numfael$, where $\numfael$ denotes the number of faces (edges in two dimensions) of the cell $\Omega_e$.

The FCFV method considers the strong form of the Poisson equation written in mixed form, via the introduction of the variable $\bq$, and in a cell-by-cell fashion, namely
\begin{equation} \label{eq:PoissonBrokenFirstOrder}
\left\{\begin{aligned}
\bq+\grad u &= \bm{0} &&\text{in $\Omega_e$, and for $e=1,\dotsc ,\numel$,}\\	
\grad\cdot\bq &= s          &&\text{in $\Omega_e$, and for $e=1,\dotsc ,\numel$,}\\
u &= u_D     &&\text{on $\partial\Omega_e\cap\Gamma_D$,}\\
\bn\cdot\bq &= -t         &&\text{on $\partial\Omega_e\cap\Gamma_N$,}\\
\jump{u\bn} &=\bm{0}  &&\text{on $\Gamma$,}\\
\jump{\bn\cdot \bq} &= 0  &&\text{on $\Gamma$.}\\
\end{aligned} \right.
\end{equation}
where $s$ is a source term, $\bn$ is the outward unit normal to the boundary, $u_D$ is the known value of the solution on the Dirichlet boundary, $t$ is the known value of the flux on the Neumann boundary and $\Gamma$, defined by
\begin{equation}\label{eq:Gamma}
\Gamma := \Big[ \bigcup_{e=1}^{\numel} \partial\Omega_e \Big]\setminus\partial\Omega
\end{equation}
is the so-called \textit{internal mesh skeleton}.

It is worth noting that the last two equations in~\eqref{eq:PoissonBrokenFirstOrder} impose the continuity of the solution and the normal flux, respectively, across the internal faces of the mesh, the \textit{jump operator} being defined as
\begin{equation}
\jump{\odot} = \odot_l + \odot_r ,
\end{equation}
that is, the sum of the quantity inside the left and right cell, $\Omega_l$ and $\Omega_r$ respectively, sharing a face of the mesh skeleton~\cite{AdM-MFH:08}.

\subsection{Strong form of the local and global problems}
\label{sc:PoissonStrong}

Following the standard rationale of HDG~\cite{Jay-CGL:09,cockburn2009hybridizable,Nguyen-NPC:09,Nguyen-NPC:10,RS-SH:16} and FCFV~\cite{FCFV2018,FCFVelas,FCFV2} methods, the strong mixed form~\eqref{eq:PoissonBrokenFirstOrder} is split into the so-called local and global problems. The local problems are defined independently in each cell as
\begin{equation} \label{eq:Dlocal-strong}
\left\{\begin{aligned}
\bq_e + \grad u_e &=\bm{0}  &&\text{in $\Omega_e$, }\\
\grad\cdot\bq_e &= s          &&\text{in $\Omega_e$,}\\		
u_e &= u_D     &&\text{on $\partial\Omega_e\cap\Gamma_D$,}\\
u_e &=\hu  &&\text{on $\partial\Omega_e\setminus\Gamma_D$,}
\end{aligned} \right.
\end{equation}
for $e=1,\dotsc ,\numel$. It is worth noting that each local problem contains only Dirichlet boundary conditions and it introduces a new independent variable, $\hu$, called the \textit{hybrid} variable, that corresponds to the solution at the cell faces (edges in two dimensions).

The global problem is defined on the mesh skeleton and the Neumann boundary as
\begin{equation} \label{eq:Dtransmission}
\left\{\begin{aligned}
\jump{u\bn} &=\bm{0}  &&\text{on $\Gamma$,}\\
\jump{\bn\cdot \bq} &= 0  &&\text{on $\Gamma$,}\\
\bn\cdot\bq &= -t         &&\text{on $\Gamma_N$.}\\
\end{aligned} \right.
\end{equation}
The first equation in~\eqref{eq:Dtransmission} can be henceforth omitted because the continuity of the solution is automatically satisfied due to the imposition of the Dirichlet boundary condition in the local problems and to the uniqueness of the hybrid variable on $\Gamma$.

\subsection{Weak form of the local and global problems}
\label{sc:PoissonWeak}

The recently proposed second-order FCFV~\cite{FCFV2} introduces a linear approximation of the primal variable in each cell, $u_e^h$, and a piecewise constant approximation of the mixed and hybrid variables, $\bq_e^h$ and $\hu^h$ respectively.

The weak formulation of the local problem in each cell is: find $(u_e^h,\bq_e^h) \in  \VoneElem \times [\VzeroElem]^{\nsd}$ such that
\begin{subequations}\label{eq:weakPoisson1}
	\begin{align} 
	- \int_{\Omega_e} \bq_e^h d\Omega = \int_{\partial\Omega_e \cap \Gamma_D} u_D \bn_e  d \Gamma + \int_{\partial\Omega_e\setminus\Gamma_D} \hu^h \bn_e  d \Gamma , \label{eq:weakPoisson1Q}
	\\
	- \int_{\Omega_e} \grad v \cdot \bq_e^h d\Omega + \int_{\partial\Omega_e}  v (\bn_e\cdot\widehat{\bq}_e^h ) d \Gamma = \int_{\Omega_e} v s d\Omega \label{eq:weakPoisson1U}
	\end{align}	
\end{subequations}
for all test functions $v \in \VoneElem$ and for $e=1,\dotsc ,\numel$. Here, $\VoneElem$ is the space of at most linear functions in $\Omega_e$ and $\VzeroElem$ is the space of constant functions in $\Omega_e$. It is worth mentioning that the weak formulation~\eqref{eq:weakPoisson1Q} is obtained by selecting an arbitrary test function in $[\VzeroElem]^{\nsd}$.

The numerical flux, $\widehat{\bq}_e^h$, introduced in~\eqref{eq:weakPoisson1U} is defined as 
\begin{equation} \label{eq:EBENumFlux}
\bn_e\cdot\widehat{\bq}_e^h := \begin{cases}
\bn_e\cdot\bq_e^h + \tau_e (\Pu u_e^h - u_D    ) & \text{on $\partial\Omega_e\cap\Gamma_D$,} \\
\bn_e\cdot\bq_e^h + \tau_e (\Pu u_e^h - \hu^h) & \text{elsewhere,}  
\end{cases}
\end{equation}
where $\tau_e>0$ is the so-called \textit{stabilisation parameter}~\cite{Jay-CGL:09,cockburn2009hybridizable,Nguyen-NPC:09,Nguyen-NPC:10,RS-SH:16} and, similar to~\cite{oikawa2015hybridized,oikawa2016analysis}, the projection operator $\Pu$ over the space of constant functions is introduced.

\begin{remark}
As discussed in~\cite{FCFV2}, the two ingredients required to obtain a second-order FCFV method are the use of piecewise linear functions to approximate the primal variable and the introduction of the projection operator in the definition of the numerical flux.
\end{remark}

The weak formulation of the local problems is obtained after introducing the definition of the numerical flux into~\eqref{eq:weakPoisson1U}, performing an integration by parts of the first term and exploiting that $\grad \cdot \bq_e^h = 0$, being $\bq_e^h$ a piecewise constant function in each cell $\Omega_e$. Hence, it reads: find $(u_e^h,\bq_e^h) \in  \VoneElem \times [\VzeroElem]^{\nsd}$ such that
\begin{subequations}\label{eq:weakPoisson}
	\begin{align} 
	- \int_{\Omega_e} \bq_e^h d\Omega & = \int_{\partial\Omega_e \cap \Gamma_D} u_D \bn_e  d \Gamma + \int_{\partial\Omega_e\setminus\Gamma_D} \hu^h \bn_e  d \Gamma, \label{eq:weakPoissonQ}
	\\
	\int_{\partial\Omega_e}  v \tau_e \Pu u_e^h d \Gamma &  = \int_{\Omega_e} v s d\Omega + \int_{\partial\Omega_e \cap \Gamma_D}  v \tau_e u_D d \Gamma + \int_{\partial\Omega_e\setminus\Gamma_D}  v \tau_e \hu^h d \Gamma \label{eq:weakPoissonU}
	\end{align}	
\end{subequations}
for all $v \in \VoneElem$ and for $e=1,\dotsc ,\numel$.

The weak formulation of the global problem is derived by following a similar procedure. It reads: find $\hu^h \in \VzeroFaces$ such that
\begin{equation} \label{eq:HDG-Poisson-Dglobal1}
\sum_{e=1}^{\numel} 
\int_{\partial\Omega_e\setminus\Gamma_D} \bn_e\cdot \widehat{\bq}_e^h d\Gamma
= -\sum_{e=1}^{\numel}  \int_{\partial\Omega_e\cap\Gamma_N} t  d\Gamma,
\end{equation}
where an arbitrary test function has been selected from $\VzeroFaces$. By using the expression of the numerical flux introduced in~\eqref{eq:EBENumFlux}, the discrete weak formulation becomes: find $\hu^h \in \VzeroFaces$ such that
\begin{equation}\label{eq:HDG-Poisson-Dglobal}
\sum_{e=1}^{\numel} 
\int_{\partial\Omega_e\setminus\Gamma_D} \left(  \bn_e\cdot\bq_e^h + \tau_e (\Pu u_e^h - \hu^h) \right) d\Gamma 
= -\sum_{e=1}^{\numel}  \int_{\partial\Omega_e\cap\Gamma_N} t  d\Gamma.
\end{equation}

\subsection{FCFV discretisation}
\label{sc:PoissonFV}

Introduce the notation $\Dset$ for the set of faces of cell $\Omega_e$ on the Dirichlet boundary $\Gamma_D$ and $\Bset$ for the faces of $\Omega_e$ in $\Gamma \cup \Gamma_N$, that is the ones not on the Dirichlet boundary. Henceforth, the stabilisation parameter, Dirichlet and Neumann data are considered to assume constant values $\tau_j$, $u_{D,j}$ and $t_j$ respectively on each face/edge $\Gamma_{e,j}$ of the cell $\Omega_e$.

The discrete local problem obtained from the weak form~\eqref{eq:weakPoisson} provides an explicit expression of the primal and mixed variables in each cell as a function of the hybrid variable on the cell faces/edges, namely
\begin{subequations}\label{eq:HDG-Poisson-DlocalK0Exp}
	\begin{align}
	\qe & =  - |\Omega_e|^{-1}\vect{z}_e  - |\Omega_e|^{-1} \sum_{j\in\Bset} |\Gamma_{e,j}| \bn_j \uHj , \label{eq:HDG-Poisson-DlocalK0QExp}
	\\
	\bue & =  \mat{m}_e^{-1} \mat{b}_e + \mat{m}_e^{-1} \sum_{j\in\Bset} \tau_j \vect{r}_j \uHj, \label{eq:HDG-Poisson-DlocalK0UExp}
	\end{align}
\end{subequations}
where $\qe$ denotes the value of the mixed variable at the centroid of the cell and $\bue$ denotes the nodal values of the primal variable. The right hand side vectors depending on the problem data are defined as
\begin{equation} \label{eq:poissonPrecomp}
\mat{b}_e := \mat{f}_e  + \sum_{j\in\Dset} \tau_j \vect{d}_j, \qquad
\vect{z}_e := \sum_{j\in\Dset} |\Gamma_{e,j}| \bn_j u_{D,j} 
\end{equation}
and the remaining matrices and vectors in the discrete equations are given by
\begin{align} 
(m_e)_{IJ} & := \sum_{j=1}^{\numfael} (\projcj)_J \tau_j \frac{1}{\nfn^{e,j}} |\Gamma_{e,j}| \chi_{\mathcal{F}_{e,j}}(I) , &  (f_e)_I & := \frac{1}{\nen} s_e |\Omega_e|,      \label{eq:matsPoisson1}  \\
(d_j)_I & := \frac{1}{\nfn^{e,j}} u_{D,j} |\Gamma_{e,j}|, &  (r_j)_I & := \frac{1}{\nfn^{e,j}} |\Gamma_{e,j}| \delta_{Ij}.    \label{eq:matsPoisson2}
\end{align}
In the above expressions, $\nen$ and $\nfn$ are the numbers of nodes in each cell and face respectively and $\delta_{Ij}$ is the Kronecker delta. The vector $\projvj$, defined as
\begin{equation} \label{eq:vectPoissonProj}
(\projcj)_I := \frac{1}{\nfn^{e,j}} \chi_{\mathcal{F}_{e,j}}(I),
\end{equation}
is used to compute the projection of the primal variable from the space of linear to the one of constant polynomial functions. In addition, $\mathcal{F}_{e,k}$ denotes the set of nodes of the cell $\Omega_e$ that belong to the face $\Gamma_{e,k}$ and the indicator function of a set $\square$ is defined as
\begin{equation} \label{eq:indicatorFunction}
\chi_{\square}(I) = 
\bigg\{
\begin{array}{ll}
1 & \text{ if } \ I\in\square \\
0 & \text{ otherwise}.
\end{array}
\end{equation}

The discrete global problem obtained by plugging the explicit expressions~\eqref{eq:HDG-Poisson-DlocalK0Exp} in the weak form~\eqref{eq:HDG-Poisson-Dglobal} results in a global system of equations where the unknown vector corresponds to the hybrid variable at the cell faces/edges. It can be written as
\begin{equation} \label{eq:globalSystemPoisson}
\mat{\widehat{K}}\vect{\hu}=\vect{\hat{f}} ,
\end{equation}
where the global matrix $\mat{\widehat{K}}$ and vector $\vect{\hat{f}}$ are the result of assembling the contribution from each cell, given by
\begin{subequations}\label{HDG-Poisson-globalSystem}
	\begin{align}
	&
	{\widehat{K}}^e_{i,j} :=  |\Gamma_{e,i}| \Big(  \tau_i \tau_j \projvj \cdot \left( \mat{m}_e^{-1} \vect{r}_j \right) - |\Omega_e|^{-1} |\Gamma_{e,j}| \bn_i \cdot \bn_j - \tau_i \delta_{ij} \Big) 
	, \\
	&
	{\widehat{f}}^e_i :=  |\Gamma_{e,i}| \Big( |\Omega_e|^{-1}  \bn_i \cdot \vect{z}_e - \tau_i \projvi \cdot \left( \mat{m}_e^{-1} \mat{b}_e\right) - t_i \, \chi_{\Nset}(i)  \Big),
	\end{align}
\end{subequations}
for $i,j \in \Bset$ and with $\delta_{ij}$ denoting the Kronecker delta.

\section{New basis functions for the second-order FCFV on general meshes}
\label{sc:FCFV2general}

The second-order FCFV method proposed in~\cite{FCFV2}, summarised in the previous section, employs nodal shape functions to define the approximation of the primal variable, namely
\begin{equation} \label{eq:lagrangeApprox}
u_e^h(\bx) = \sum_{J=1}^{\nen} N^e_J(\bx) u^e_J,
\end{equation}
where $\{ N^e_J \}_{J=1}^{\nen}$ is the set of linear Lagrange polynomials in the cell $\Omega_e$ and $u^e_J$, for $J=1,\ldots,\nen$, are the corresponding nodal values of the unknown function.

This section shows that this approach is only applicable with simplicial (triangular and tetrahedral) cells and proposes a new approximation space for general polygons and polyhedrons in two and three dimensions respectively. The proposed formulation is first presented for the Poisson equation introduced in section~\ref{sc:FCFV2}, and is extended to Stokes equations in section~\ref{sc:FCFV2generalStokesFV}.

\subsection{The second-order FCFV with nodal basis functions}
\label{sc:FCFVquaProblem}

Consider the Poisson model problem~\eqref{eq:PoissonBrokenFirstOrder} in two dimensions. The matrix $\mat{m}_e$ used in the local problem to write the solution in the cell as a function of the solution on the faces can thus be computed analytically. For a triangular cell, and assuming a constant value of the stabilisation parameter $\tau_e$ for all the faces, the matrix is given by
\begin{equation} \label{eq:meTri}
\mat{m}_e = \frac{\tau_e}{4}
\begin{bmatrix}
|\Gamma_{e,1}| + |\Gamma_{e,3}| & |\Gamma_{e,1}| & |\Gamma_{e,3}| \\
|\Gamma_{e,1}| & |\Gamma_{e,2}| + |\Gamma_{e,1}| & |\Gamma_{e,2}| \\
|\Gamma_{e,3}| & |\Gamma_{e,2}| & |\Gamma_{e,3}| + |\Gamma_{e,2}| 
\end{bmatrix} .
\end{equation}
This matrix is invertible, with determinant equal to $(\tau_e^3/16) |\Gamma_{e,1}||\Gamma_{e,2}||\Gamma_{e,3}|$.

In a similar fashion, assume a constant value of the stabilisation parameter $\tau_e$ for all the faces of a quadrilateral cell. Using a nodal approximation with bilinear Lagrange polynomials, the matrix is given by
\begin{equation} \label{eq:meQua}
\mat{m}_e = \frac{\tau_e}{4}
\begin{bmatrix}
|\Gamma_{e,1}| + |\Gamma_{e,4}| & |\Gamma_{e,1}| & 0 & |\Gamma_{e,4}| \\
|\Gamma_{e,1}| & |\Gamma_{e,2}| + |\Gamma_{e,1}| & |\Gamma_{e,2}| & 0 \\
0 & |\Gamma_{e,2}| & |\Gamma_{e,3}| + |\Gamma_{e,2}| & |\Gamma_{e,3}| \\
|\Gamma_{e,4}| & 0 & |\Gamma_{e,3}| & |\Gamma_{e,4}| + |\Gamma_{e,3}|
\end{bmatrix} .
\end{equation}
Contrary to the matrix obtained in~\eqref{eq:meTri} for a triangular cell, the matrix $\mat{m}_e$ for a quadrilateral cell is singular and the local problem in~\eqref{eq:HDG-Poisson-DlocalK0Exp} cannot be solved.
Hence, in order to devise a second-order FCFV method suitable to handle quadrilateral cells in 2D and hexahedral, prismatic and pyramidal cells in 3D, an alternative description of the primal variable needs to be considered.

\subsection{New linear basis functions for the second-order FCFV method}
\label{sc:FCFVmodalBasis}

This work proposes the use of a linear approximation of the primal variable in each cell, irrespective of its shape. More precisely, the approximation of the primal variable is defined as
\begin{equation} \label{eq:modalApprox}
u_e^h(\bx) = \sum_{J=1}^{\nmo} \basis^e_J(\bx) c^e_J,
\end{equation}
where the number of terms of the expansion is selected as $\nmo = \nsd+1$, $\{ \basis^e_J \}_{J=1}^{\nmo}$ is a set of basis functions that span the space of polynomials of, at most, degree one and $c^e_J$, for $J=1,\ldots,M$ are coefficients appropriately defined to describe the unknown function.

The following basis functions are proposed in this work
\begin{equation} \label{eq:modalBasis}
\basis^e_1(\bx) = 1 \quad \text{and} \quad \basis^e_{k}(\bx) = x_{k-1} - \bar{x}_{k-1}^e \quad \text{for $k=2,\ldots\nmo$}.
\end{equation}
where $\bar{\bx}^e = (\bar{x}_1^e, \ldots, \bar{x}_{\nsd}^e)$ denotes the coordinates of the centroid of the cell $\Omega_e$.

It is worth noting that other choices for the basis functions are possible. The choice made here ensures that the coefficients in the approximation~\eqref{eq:modalApprox} have a physical interpretation. More precisely, $c^e_1$ is the value of the primal variable at the centroid of the cell, whereas each coefficient $c^e_k$, for $k=2,\ldots\nmo$, corresponds to the value of the derivative of the primal variable in the $x_k$ direction, at the centroid of the cell.

\begin{remark}
The proposed FCFV scheme that results from considering the new set of linear basis functions reduces to the original second-order FCFV in~\cite{FCFV2} for triangular and tetrahedral cells as in both cases the basis functions span the space of polynomials of, at most, degree one. 
\end{remark}

\subsection{Second-order FCFV discretisation of the Poisson equation}
\label{sc:FCFV2generalPoissonFV}

Considering the approximation of the primal variable proposed in the previous section and a constant approximation of the mixed and hybrid variables, the discrete local problem obtained from the weak form~\eqref{eq:weakPoisson} is
\begin{subequations}\label{eq:HDG-Poisson-DlocalModal}
	\begin{align}
	\qe & =  - |\Omega_e|^{-1}\vect{z}_e  - |\Omega_e|^{-1} \sum_{j\in\Bset} |\Gamma_{e,j}| \bn_j \uHj , \label{eq:HDG-Poisson-DlocalQModal}
	\\
	\bue & =  \mat{\widetilde{m}}_e^{-1} \mat{\tilde{b}}_e + \mat{\widetilde{m}}_e^{-1} \sum_{j\in\Bset} \tau_j \vect{\tilde{r}}_j \uHj. \label{eq:HDG-Poisson-DlocalUModal}
	\end{align}
\end{subequations}

It is worth noting that equation~\eqref{eq:HDG-Poisson-DlocalQModal}, expressing the mixed variable in terms of the hybrid variable, is identical to the discrete equation~\eqref{eq:HDG-Poisson-DlocalK0QExp} of the original second-order FCFV, whereas equation~\eqref{eq:HDG-Poisson-DlocalUModal} requires the following definitions
\begin{equation} \label{eq:poissonPrecompModal}
\mat{\tilde{b}}_e := \mat{\widetilde{f}}_e  + \sum_{j\in\Dset} \tau_j \vect{\tilde{d}}_j, \qquad
\end{equation}
and 
\begin{align} 
(\widetilde{m}_e)_{IJ} & := \sum_{j=1}^{\numfael} (\projcjT)_J  \tau_j \int_{\Gamma_{e,j}}  \basis_I d\Gamma , &  (\tilde{f}_e)_I & := \int_{\Omega_e}  \basis_I s d\Omega,      \label{eq:matsPoisson1Modal}  \\
(\tilde{d}_j)_I & := u_{D,j} \int_{\Gamma_{e,j}}  \basis_I d\Gamma , &  (\tilde{r}_j)_I & := \int_{\Gamma_{e,j}}  \basis_I d\Gamma.    \label{eq:matsPoisson2Modal}
\end{align}

The vector $\projvjT$, given by
\begin{equation} \label{eq:vectPoissonProjModal}
(\projcjT)_1 = 1 \quad \text{and} \quad (\projcjT)_k = \bar{x}_{k-1}^{e,j} - \bar{x}_{k-1}^e \quad \text{for $k=2,\ldots\nmo$} 
\end{equation}
is introduced to compute the projection of the primal linear variable on the space of the constant functions over a face $\Gamma_{e,j}$, where $\bar{\bx}^{e,j} = (\bar{x}_1^{e,j}, \ldots, \bar{x}_{\nsd}^{e,j})$ denotes the centroid of the face $\Gamma_{e,j}$.

Analogously, the discrete global problem obtained from the weak form~\eqref{eq:HDG-Poisson-Dglobal} using the expressions of primal and mixed variable in~\eqref{eq:HDG-Poisson-DlocalModal}, results in a global system of equations where the unknown vector corresponds to the hybrid variable at the cell faces/edges. It can be written as
\begin{equation} \label{eq:globalSystemPoissonModal}
\mat{\widehat{K}}\vect{\hu}=\vect{\hat{f}} ,
\end{equation}
where the global matrix $\mat{\widehat{K}}$ and vector $\vect{\hat{f}}$ are the result of assembling the contribution from each cell, given by
\begin{subequations}\label{HDG-Poisson-globalSystemModal}
	\begin{align}
	&
	{\widehat{K}}^e_{i,j} :=  |\Gamma_{e,i}| \Big(  \tau_i \tau_j \projvjT \cdot \left( \mat{\widetilde{m}}_e^{-1} \vect{\tilde{r}}_j \right) - |\Omega_e|^{-1} |\Gamma_{e,j}| \bn_i \cdot \bn_j - \tau_i \delta_{ij} \Big) 
	, \\
	&
	{\widehat{f}}^e_i :=  |\Gamma_{e,i}| \Big( |\Omega_e|^{-1}  \bn_i \cdot \vect{z}_e - \tau_i \projviT \cdot \left( \mat{\widetilde{m}}_e^{-1} \mat{\tilde{b}}_e\right) - t_i \, \chi_{\Nset}(i)  \Big),
	\end{align}
\end{subequations}
for $i,j \in \Bset$.

\subsection{Second-order FCFV discretisation of the Stokes equations}
\label{sc:FCFV2generalStokesFV}

In this section, the procedure described above for the derivation of the second-order FCFV formulation of the Poisson problem is extended to Stokes equations.
By introducing the mixed variable $\bL$, the strong form of the mixed problem is written cell-by-cell as
\begin{equation} \label{eq:StokesBroken}
\left\{\begin{aligned}
\bL + \sqrt{\nu} \grad\bu & = \bm{0} &&\text{in $\Omega_e$ and for $e=1,\ldots,\numel$,}\\
\grad\cdot\bigl(\sqrt{\nu} \bL + p\Insd \bigr) &= \bm{s}         &&\text{in $\Omega_e$ and for $e=1,\ldots,\numel$,}\\
\grad\cdot\bu &= 0  &&\text{in $\Omega_e$ and for $e=1,\ldots,\numel$,}\\
\bu &= \bu_D     &&\text{on $\partial\Omega_e \cap \Gamma_D$,}\\
\bn \cdot \bigl(\sqrt{\nu} \bL + p \Insd\bigr)  &= - \bm{t}  &&\text{on $\partial\Omega_e \cap \Gamma_N$,}\\  		
\jump{\bu \otimes \bn} &= \bm{0}  &&\text{on $\Gamma$,}\\
\jump{\bn \cdot \bigl(\sqrt{\nu} \bL + p \Insd\bigr)} &= \bm{0}  &&\text{on $\Gamma$,}\\
\end{aligned} \right.
\end{equation}
with the solvability constraint for the uniqueness of pressure given by
\begin{equation}\label{eq:pressureConstraint}
\frac{1}{|\partial\Omega_e|} \int_{\partial\Omega_e}  p_e d\Gamma  = \rho_e ,
\end{equation}
where $\rho_e$ is the mean value of the pressure on the boundary of the cell $\Omega_e$~\cite{donea2003finite,Nguyen-CNP:10}.
As for the scalar problem, the last two equations in~\eqref{eq:StokesBroken} impose the continuity of the velocity and the normal flux, respectively, across the internal faces of the mesh.

Following the rationale in~\cite{FCFV2018,FCFV2}, an additional unknown, the hybrid velocity $\bhu$, is introduced on the cell faces and the condition $\bu = \bhu$ is enforced on the boundary $\partial\Omega_e \setminus \Gamma_D$ for $e=1,\ldots,\numel$.
The weak formulation of the Stokes problem in each cell is: given $\bu_D$ on $\partial\Omega_e \cap \Gamma_D$ and $\bhu^h$ on $\partial\Omega_e \setminus \Gamma_D$, find $(\bu_e^h,p_e^h,\bL_e^h) \in  [\VoneElem]^{\nsd} \times \VzeroElem \times [\VzeroElem]^{\nsd\times\nsd}$ such that
\begin{subequations}\label{HDG-Stokes-DlocalSymm}
	\begin{align}
	&
	-\int_{\Omega_e} \bL^h_e d\Omega = \int_{\partial\Omega_e\cap\Gamma_D}  \sqrt{\nu} \bn_e \otimes \bu_D d\Gamma + \int_{\partial\Omega_e\setminus\Gamma_D}  \sqrt{\nu} \bn_e \otimes \bhu^h d\Gamma, \label{eq:HDG-Stokes-Dlocal-LSymm} 
	\\
	&
	\int_{\partial \Omega_e} \tau_e \bw \cdot \Pu \bu^h_e  d\Gamma
	=  \int_{\Omega_e}\bw \cdot \bm{s} d\Omega   + \int_{\partial\Omega_e\cap\Gamma_D}  \tau_e \bw \cdot \bu_D d\Gamma + \int_{\partial\Omega_e\setminus\Gamma_D}  \tau_e \bw \cdot \bhu^h d\Gamma,
 	 \label{eq:HDG-Stokes-Dlocal-MomSymm}
	\\
	&
	\int_{\partial\Omega_e\setminus\Gamma_D} \bhu \cdot \bn_e d\Gamma + \int_{\partial\Omega_e\cap\Gamma_D} \bu_D \cdot \bn_e  d\Gamma = 0 , 
	\label{eq:HDG-Stokes-Dlocal-USymm} 
	\\
	&
	\frac{1}{|\partial\Omega_e|} \int_{\partial\Omega_e}  p^h_e d\Gamma  = \rho_e^h , \label{eq:HDG-Stokes-Dlocal-PSymm}
	\end{align}
\end{subequations}
for all test functions $\bw \in [\VoneElem]^{\nsd}$, where the definition of the numerical normal flux featuring the projection operator $\Pu$ is utilised
\begin{equation} \label{eq:traceStokes}
\bn_e \! \cdot \!\bigl(\widehat{\sqrt{\nu} \bL^h_e \!+\! p^h_e \Insd} \bigr) \!:= \!
\begin{cases}
\bn_e \! \cdot \! \bigl( \sqrt{\nu} \bL^h_e \!+\! p^h_e \Insd \bigr) \!+\! \tau_e (\Pu \bu^h_e \!-\! \bu_D) & \text{on $\partial\Omega_e\cap\Gamma_D$,} \\
\bn_e \!\cdot \!\bigl( \sqrt{\nu} \bL^h_e \!+\! p^h_e \Insd \bigr) \!+\! \tau_e (\Pu \bu^h_e \!-\! \bhu^h) & \text{elsewhere.}  
\end{cases}
\end{equation}
In equation~\eqref{HDG-Stokes-DlocalSymm}, $[\VoneElem]^{\nsd}$ is the space of $\nsd$ dimensional vectors whose components are at most linear functions in $\Omega_e$, whereas $[\VzeroElem]^{\nsd\times\nsd}$ and $\VzeroElem$ are the spaces of constant $\nsd \times \nsd$ tensorial and scalar functions in the cell $\Omega_e$, respectively.
It is worth mentioning that the weak formulations in equations~\eqref{eq:HDG-Stokes-Dlocal-LSymm} and \eqref{eq:HDG-Stokes-Dlocal-USymm} are obtained by selecting an arbitrary constant test function in the spaces $[\VzeroElem]^{\nsd \times \nsd}$ and $\VzeroElem$, respectively. 

The weak formulation of the global problem is derived in a similar way by imposing Neumann and transmission conditions on $\Gamma_N$ and $\Gamma$, respectively, and a compatibility condition for the incompressiblity constraint in each cell $\Omega_e$, for $e=1,\ldots,\numel$. It reads: find $(\bhu^h,\rho_e^h) \in [\VzeroFaces]^{\nsd} \times \mathbb{R}$ such that
\begin{subequations}\label{HDG-Stokes-Dglobal}
	\begin{align}
	&
	\sum_{e=1}^{\numel} 
	\int_{\partial\Omega_e\setminus\Gamma_D} \left(  \sqrt{\nu} \bn_e\cdot \bL^h_e + p^h_e \bn_e + \tau_e (\Pu \bu^h_e - \bhu^h) \right) d\Gamma 
	= -\sum_{e=1}^{\numel}  \int_{\partial\Omega_e\cap\Gamma_N} \bm{t}  d\Gamma,
	\label{eq:HDG-Stokes-Dglobal1}
	\\
	&
	\int_{\partial\Omega_e\setminus\Gamma_D} \bhu \cdot \bn_e d\Gamma + \int_{\partial\Omega_e\cap\Gamma_D} \bu_D \cdot \bn_e  d\Gamma = 0 \quad \text{ for } e=1,\dotsc,\numel ,
	\label{eq:weakCompatibilityStokes}
	\end{align}
\end{subequations}
where an arbitrary test function in the space $[\VzeroFaces]^{\nsd}$ has been selected in equation~\eqref{eq:HDG-Stokes-Dglobal1}.

\begin{remark}
Equation~\eqref{eq:HDG-Stokes-Dlocal-USymm} of the local problem coincides with the compatibility condition~\eqref{eq:weakCompatibilityStokes} enforcing the divergence-free nature of the velocity field cell-by-cell.
Following the strategy discussed for the first and second-order FCFV method~\cite{FCFV2018,FCFV2}, this equation is omitted from the local problems and imposed solely in the global one since it involves only the global variable $\bhu^h$.
\end{remark}

The discrete FCFV local problems are obtained starting from the linear discretisation described in section~\ref{sc:FCFVmodalBasis} for the velocity $\bu_e^h$, a constant cell-by-cell approximation of the pressure $p_e^h$, the mixed variable $\bL_e^h$ and the mean pressure $\rho_e^h$ and a constant face-by-face approximation of the hybrid velocity $\bhu^h$. It follows that
\begin{subequations}\label{HDG-Stokes-localModal}
	\begin{align}
	&
	 \Le =  -| \Omega_e |^{-1}\sqrt{\nu}\mat{Z}_e    -| \Omega_e |^{-1} \sqrt{\nu} \sum_{j \in \Bset} | \Gamma_{e,j} | \bn_j \otimes \buHj  , \label{eq:HDG-Stokes-Dlocal-LModal} 
	\\
	& \bue =  \widetilde{\mat{M}}_e^{-1} \widetilde{\mat{B}}_e + \widetilde{\mat{M}}_e^{-1} \sum_{j\in\Bset} \tau_j \widetilde{\mat{R}}_j \buHj, \label{eq:HDG-Stokes-Dlocal-MomModal} 
	\\
	&
	\pe = \rho_e , \label{eq:HDG-Stokes-Dlocal-PModal}
	\end{align}
\end{subequations}
where 
\begin{equation}
\mat{Z}_e := \sum_{j\in\Dset} |\Gamma_{e,j}| \bn_j \otimes \bu_{D,j} 
\quad \text{and} \quad
\widetilde{\mat{B}}_e := \widetilde{\mat{F}}_e  + \sum_{j\in\Dset} \tau_j \widetilde{\mat{D}}_j .
\end{equation}
It is worth noticing that equations~\eqref{eq:HDG-Stokes-Dlocal-LModal} and~\eqref{eq:HDG-Stokes-Dlocal-PModal} are identical to the original second-order FCFV~\cite{FCFV2} and only equation~\eqref{eq:HDG-Stokes-Dlocal-MomModal} is affected by the change of basis discussed above. More precisely,
\begin{align} 
	(\widetilde{\mat{M}}_e)_{IJ} & := \Insd \sum_{j=1}^{\numfael} (\projcjT)_J \tau_j \int_{\Gamma_{e,j}}  \basis_I d\Gamma, &  (\widetilde{\mat{F}}_e)_I & := \int_{\Omega_e}  \basis_I \bm{s} d\Omega,      \label{eq:matsStokes1Modal}  \\
	(\widetilde{\mat{D}}_j)_I & := \bu_{D,j} \int_{\Gamma_{e,j}}  \basis_I  d\Gamma, &  (\widetilde{\mat{R}}_j)_{IJ} & := \Insd \int_{\Gamma_{e,j}}  \basis_I d\Gamma.    \label{eq:matsStokes2Modal}
\end{align}

The discrete FCFV global system is obtained by plugging the expressions~\eqref{HDG-Stokes-localModal} of $\bue$, $\pe$ and $\Le$ into equation~\eqref{HDG-Stokes-Dglobal}, leading to 
\begin{equation} \label{eq:globalSystemStokes}
\begin{bmatrix}
\mat{\widehat{K}}_{\hu \hu} & \mat{\widehat{K}}_{\hu \rho} \\
\mat{\widehat{K}}_{\hu \rho}^T & \mat{0}_{\numel}
\end{bmatrix}
\begin{Bmatrix}
\vect{\hat{u}}  \\
\bm{\rho}
\end{Bmatrix}
=
\begin{Bmatrix}
\vect{\hat{f}}_{\hu} \\
\vect{\hat{f}}_{\rho}
\end{Bmatrix} .
\end{equation}
It is straightforward to observe that the matrix of problem~\eqref{eq:globalSystemStokes} features a saddle-point structure with a symmetric block $\mat{\widehat{K}}_{\hu \hu}$ in the top-left position, as classical in the approximation of incompressible Stokes equations.
Both the left and right hand sides of the linear system above are assembled by computing for $i,j \in \Bset$ the contributions of each cell as
\begin{subequations}\label{HDG-Stokes-globalSystemModal}
	\begin{align}
	&
	(\mat{\widehat{K}}_{\hu \hu})^e_{i,j} :=  |\Gamma_{e,i}| \left[ \tau_i \tau_j \projmiT \left( \widetilde{\mat{M}}_e^{-1} \widetilde{\mat{R}}_j \right) - \nu | \Omega_e |^{-1} |\Gamma_{e,j}| (\bn_i \cdot \bn_j)\mat{I}_{\nsd}  - \tau_i \delta_{ij} \mat{I}_{\nsd} \right] 
	, \\
	&
	(\mat{\widehat{K}}_{\hu \rho})^e_i := |\Gamma_{e,i}| \bn_i , \\
	&
	(\vect{\hat{f}}_{\hu})^e_i :=  |\Gamma_{e,i}| \left( \nu | \Omega_e |^{-1} \bn_i \cdot \mat{Z}_e  - \tau_i \projmiT \left( \widetilde{\mat{M}}_e^{-1} \widetilde{\mat{B}}_e \right)  - \bm{t}_i \, \chi_{\Nset}(i) \right) , \\
	&
	(\hat{\text{f}}_{\rho})^e := -\sum_{j \in \Dset} | \Gamma_{e,j} | \bu_{D,j} \cdot \bn_j .
	\end{align}
\end{subequations}
The operator utilised to project the linear velocity in each cell on the space of constant functions on the face $\Gamma_{e,i}$ is defined via the matrix $\projmiT$
\begin{equation} \label{eq:matStokesProj-2D}
\projmiT := 
\begin{bmatrix}
\projviT^T &  \bm{0}_{1 \times 3}  \\
\bm{0}_{1 \times 3} &  \projviT^T
\end{bmatrix}
\end{equation}
in 2D and
\begin{equation} \label{eq:matStokesProj-3D}
\projmiT := 
\begin{bmatrix}
\projviT^T &  \bm{0}_{1 \times 4} &  \bm{0}_{1 \times 4}  \\
\bm{0}_{1 \times 4} &  \projviT^T &  \bm{0}_{1 \times 4}  \\
\bm{0}_{1 \times 4}  &  \bm{0}_{1 \times 4}&  \projviT^T
\end{bmatrix}
\end{equation}
in 3D, $\projviT$ being the vector introduced in~\eqref{eq:vectPoissonProjModal} for the projection in the scalar case.

\subsection{Computational aspects}
\label{sc:computational}

The implementation of a FCFV method features three steps. First, a preprocess routine to compute all the elemental quantities required by the method. Second, the computation, assembly and solution of a global problem on the mesh faces with the hybrid variable as unknown. Finally, the solution of $\numel$ local problems to retrieve the primal and mixed variable in each cell.

The current implementation assumes that the stabilisation parameter is constant in the whole domain. Therefore, for the Poisson problem, the entries of the matrices $\mat{\tilde{m}}_e$ in equation~\eqref{eq:matsPoisson1Modal} and the vectors $\vect{\tilde{d}}_j$ and $\vect{\tilde{r}}_j$ in equation~\eqref{eq:matsPoisson2Modal} can be precomputed. Analogously, for the Stokes problem, the entries of the matrices $\widetilde{\mat{M}}_e$ and $\widetilde{\mat{R}}_j$ in equation~\eqref{eq:matsStokes1Modal} and~\eqref{eq:matsStokes2Modal} respectively and the vector $\widetilde{\mat{D}}_j$ in equation~\eqref{eq:matsStokes2Modal} can be precomputed. These terms only require the calculation of the integral of the new basis functions over a generic face of a cell. More precisely, for a generic face in two dimensions and for a triangular face in three dimensions, the integrals are given by
\begin{equation} \label{eq:basisIntegralFace}
\int_{\Gamma_{e,j}} \basis_1 d\Gamma = | \Gamma_{e,j} |, \quad \text{and} \quad 
\int_{\Gamma_{e,j}} \basis_k d\Gamma = | \Gamma_{e,j} | \left( \bar{x}_{k-1}^{e,j} - \bar{x}_{k-1}^e \right), \quad \text{for $k=2,\ldots\nmo$}.
\end{equation}
For a generic quadrilateral face in three dimensions, analytical integration is only feasible for the first basis function, $\basis_1$, whereas for the remaining basis functions, $\basis_k$ for $k=2,\ldots\nmo$, a numerical quadrature is employed.

Similarly, the entries of the vectors $\mat{\tilde{f}}_e$ and $\widetilde{\mat{F}}_e$ in equations~\eqref{eq:matsPoisson1Modal} and ~\eqref{eq:matsStokes1Modal}, respectively for the Poisson and Stokes problems, can be precomputed via an integral of the new basis functions over a generic cell. For a triangular or tetrahedral cell, such integral is given by
\begin{equation} \label{eq:basisIntegralCell}
\int_{\Omega_e} \basis_1 d\Omega = | \Omega_e |, \quad \text{and} \quad 
\int_{\Omega_e} \basis_k d\Omega = 0, \quad \text{for $k=2,\ldots\nmo$},
\end{equation}
whereas for other cell types, the integral is computed numerically.
Hence, the computational cost of the preprocess routine to determine the terms in~\eqref{eq:matsPoisson1Modal}-\eqref{eq:matsPoisson2Modal} and~\eqref{eq:matsStokes1Modal}-\eqref{eq:matsStokes2Modal} is limited.

The proposed second-order FCFV method requires the solution of a global system of equations with exactly the same number of unknowns and non-zero entries of the global matrix of the first-order FCFV proposed in~\cite{FCFV2018,FCFVelas}. This is because the second-order FCFV uses the same constant approximation space for the hybrid variable, which is the only variable featuring in the global systems of equations~\eqref{HDG-Poisson-globalSystemModal} and \eqref{HDG-Stokes-globalSystemModal}, respectively for the Poisson and Stokes problems. 

The extra cost of the second-order FCFV compared to the first-order FCFV is due to the extra operations required to assemble the global system. Tables~\ref{tab:operationsPoisson} and \ref{tab:operationsStokes} detail the number of operations required to compute one elemental matrix and one elemental right hand side vector for the Poisson and Stokes problems respectively, for the proposed second-order FCFV method and the first-order FCFV method introduced in~\cite{FCFV2018}. 
More precisely, table~\ref{tab:operationsPoisson} presents the cost of computing one elemental contribution of equation~\eqref{HDG-Poisson-globalSystemModal} and one of equation~22 in~\cite{FCFV2018}, whereas table~\ref{tab:operationsStokes} compares the cost of building one elemental block of equation~\eqref{HDG-Stokes-globalSystemModal} with one of equation~39 in~\cite{FCFV2018}.
\begin{table}[hbt]
	\centering
	\begin{tabular}[hbt]{| c || c | c | c | c | c | c |}
		\hline
		Method & Triangle & Quadrilateral & Tetrahedron & Hexahedron & Prism & Pyramid \\
		\hline & & & & & &
		\\ [-1em] 
		\hline
		First-order & 126 & 212 & 252 & 534  & 380 & 380 \\
		\hline
		Second-order & 354 & 592 & 932 & 1,962 & 1,400 & 1,400 \\
		\hline
	\end{tabular}
	\caption{Number of operations required to compute the elemental matrix and right hand side of the Poisson problem by the first and second-order FCFV methods for different cell types.}
	\label{tab:operationsPoisson}
\end{table}

\begin{table}[hbt]
	\centering
	\begin{tabular}[hbt]{| c || c | c | c | c | c | c |}
		\hline
		Method & Triangle & Quadrilateral & Tetrahedron & Hexahedron & Prism & Pyramid \\
		\hline & & & & & &
		\\ [-1em] 
		\hline
		First-order & 168 & 272 & 364 & 714  & 525 & 525 \\
		\hline
		Second-order &  2,340 & 3,900 & 21,424 & 44,988 & 32,135 & 32,135 \\
		\hline
	\end{tabular}
	\caption{Number of operations required to compute the elemental matrix and right hand side of the Stokes problem by the first and second-order FCFV methods for different cell types.}
	\label{tab:operationsStokes}
\end{table}

\begin{figure}[!bt]
	\centering
	\subfigure[Triangular cells]{\includegraphics[width=0.32\textwidth]{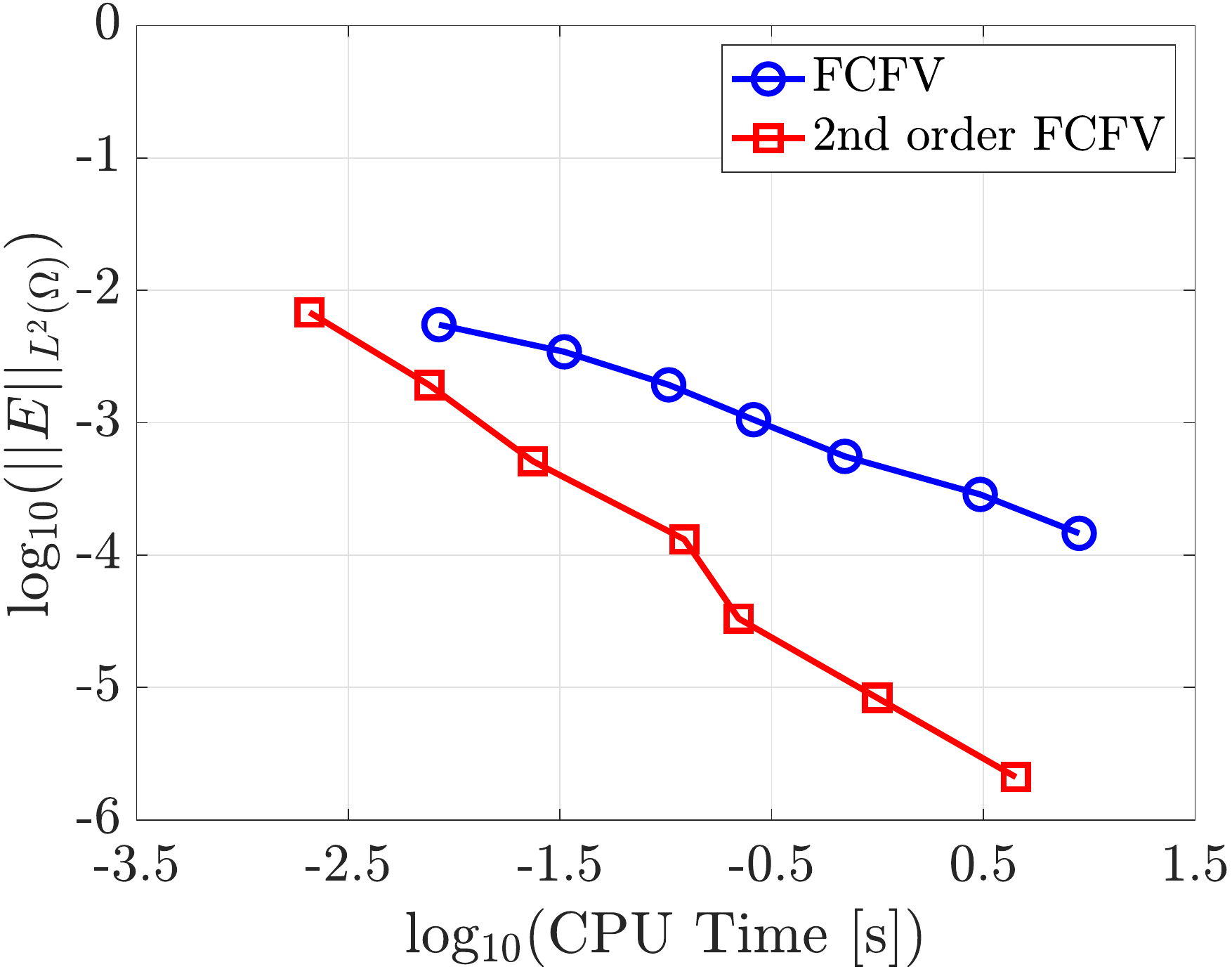}}
	\subfigure[Quadrilateral cells]{\includegraphics[width=0.32\textwidth]{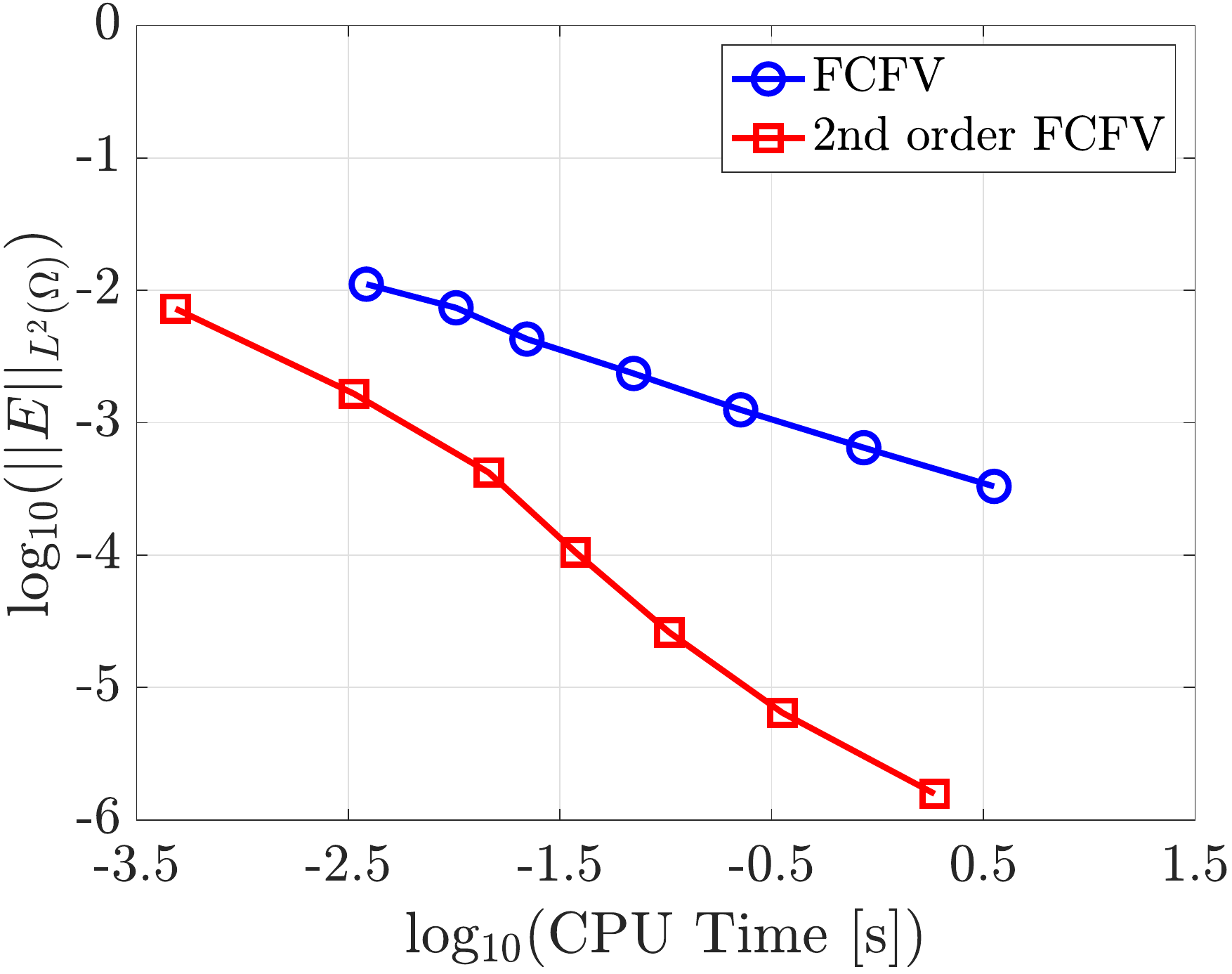}}
	\subfigure[Tetrahedral cells]{\includegraphics[width=0.32\textwidth]{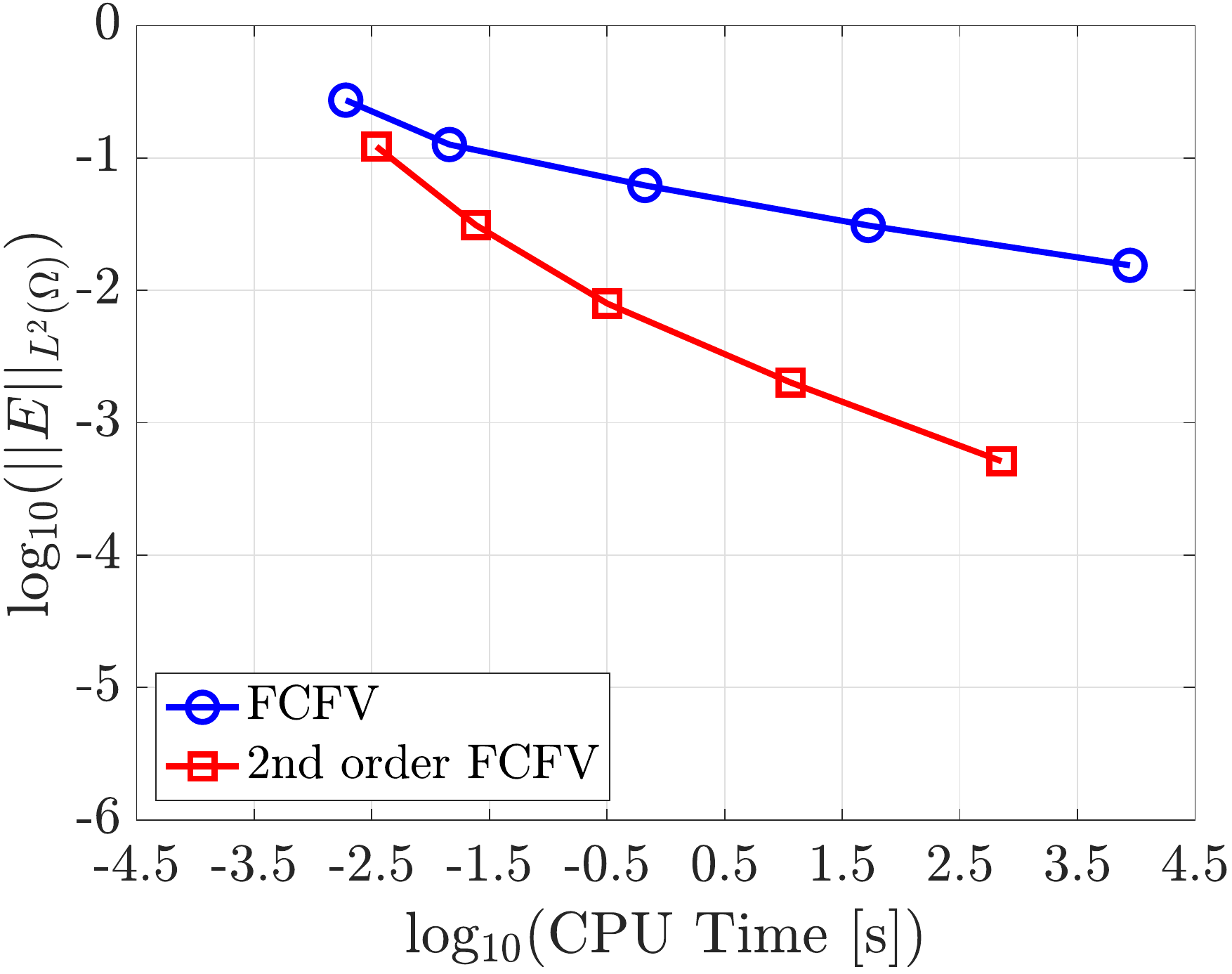}}
	\subfigure[Hexahedral cells]{\includegraphics[width=0.32\textwidth]{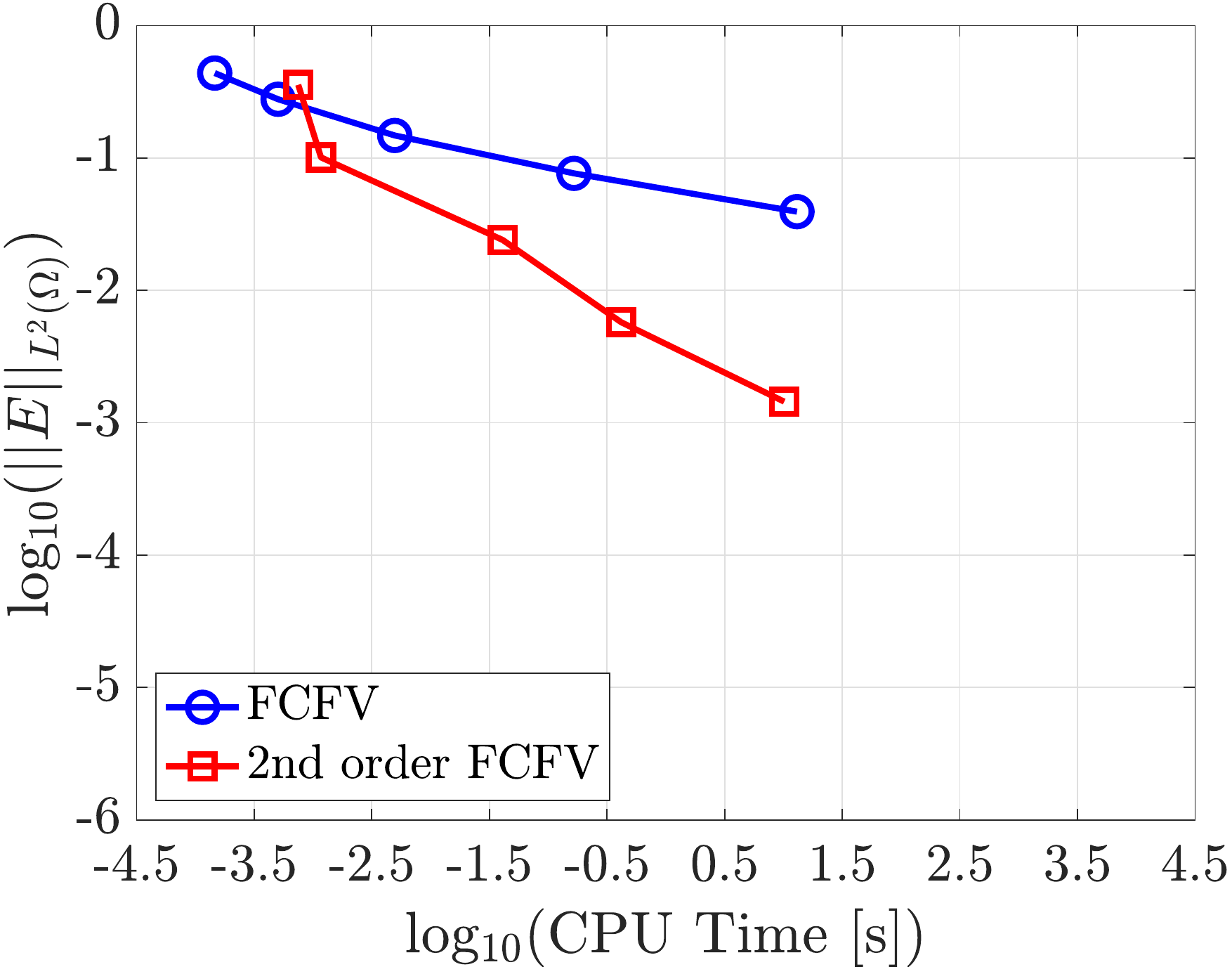}} 
	\subfigure[Prismatic cells]{\includegraphics[width=0.32\textwidth]{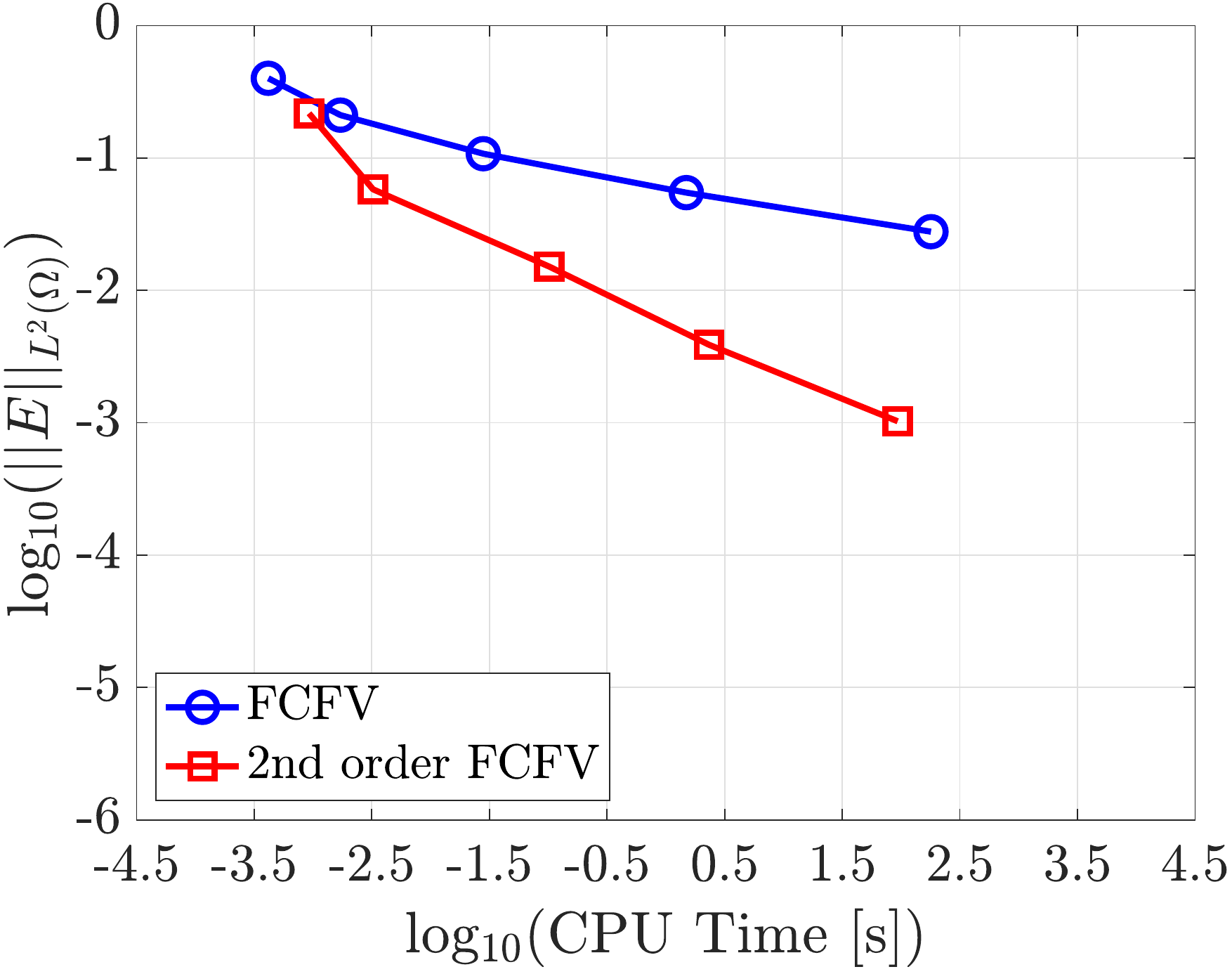}}
	\subfigure[Pyramidal cells]{\includegraphics[width=0.32\textwidth]{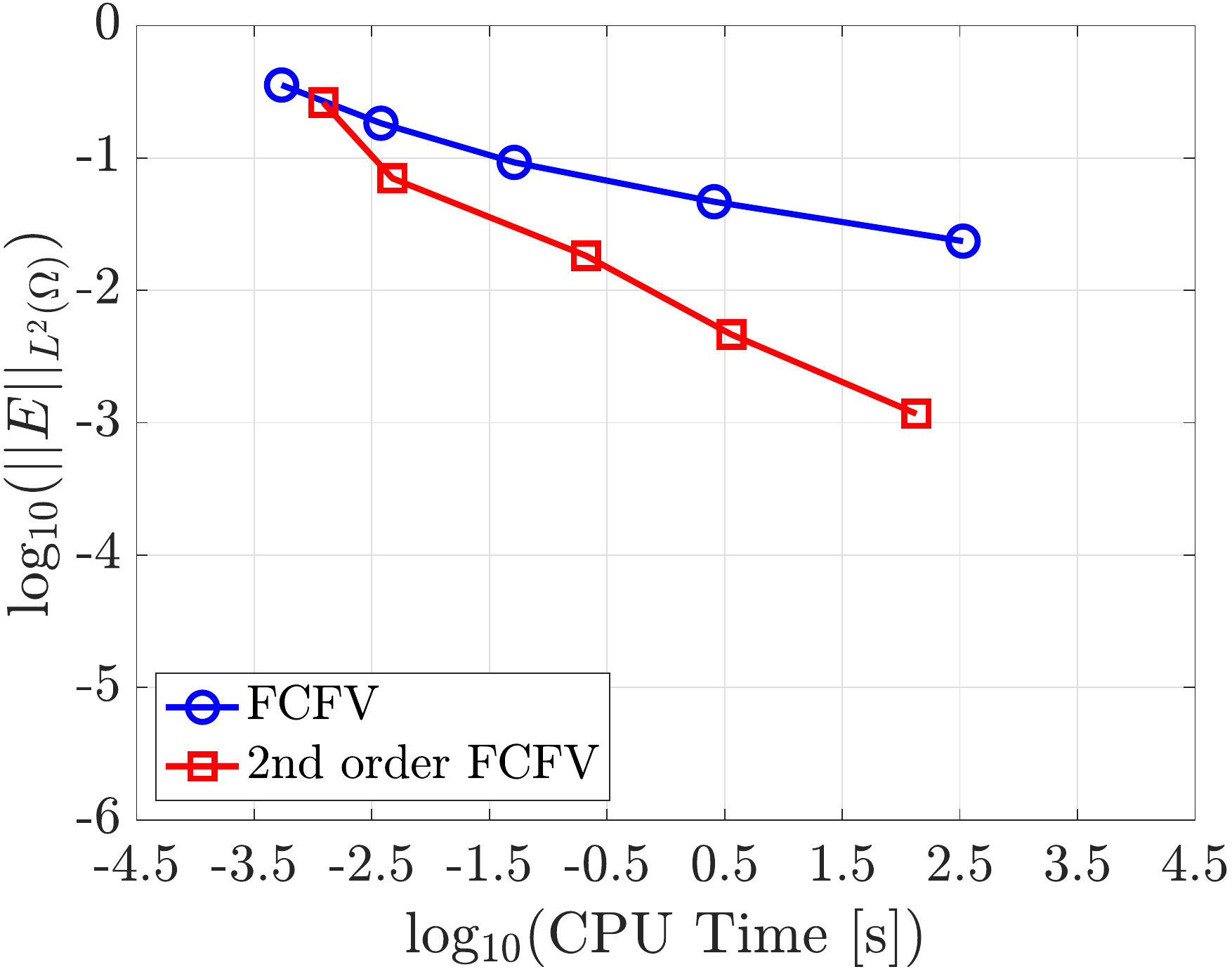}}	
	\caption{Relative error of the velocity measured in the $\eltwo(\Omega)$ norm as a function of the CPU time for the solution of the Stokes problem using (a) triangular, (b) quadrilateral, (c) tetrahedral, (d) hexahedral, (e) prismatic and (f) pyramidal cells.}
	\label{fig:CPUtimeStokesK0vsK0K1}
\end{figure}

Although the second-order method requires more operations for a given spatial discretisation, the extra accuracy provided results in a more efficient method, when the computational cost required to achieve a given accuracy is considered. To quantify this gain, Figure~\ref{fig:CPUtimeStokesK0vsK0K1} shows the relative error of the velocity field measured in the $\eltwo(\Omega)$ norm as a function of the CPU time for the first and second-order methods using different cell types. The test case considered involves the solution of the Stokes equations, in two and three dimensions, for a problem described in the next section, where the analytical solution is known.
The results show that to achieve a 1\% error, the second-order method is almost one order of magnitude faster in two dimensions whereas in three dimensions the second-order method is more than four order of magnitude faster.

Finally, the solution of equation~\eqref{eq:HDG-Poisson-DlocalModal} and~\eqref{HDG-Stokes-localModal} for the local Poisson and Stokes problems respectively, relies on precomputed quantities and can be easily performed in parallel being such computation independent cell-by-cell.

\section{An automatic mesh adaptivity strategy for the second-order FCFV}
\label{sc:meshAdaptivity}

In~\cite{FCFV2}, the authors devised an error indicator using the higher convergence rate of the second-order FCFV method with respect to the original first-order FCFV approach~\cite{FCFV2018,FCFVelas}.
The main drawback of the proposed error indicator is the high computational cost that induces the solution of an extra global problem to estimate the error of a numerical solution.

This work proposes a new error indicator that does not require the solution of an extra global problem. Instead, only a local computation cell-by-cell is required to devise an accurate and efficient indicator to drive mesh adaptivity. 
The error indicator proposed here is significantly cheaper because it exploits the solution of the global system already computed for the second-order FCFV method to solve an extra local problem cell-by-cell. Hence, its cost is negligible when compared to the cost of assembling and solving an extra global problem as discussed in~\cite{FCFV2}.

Hereafter, the proposed strategy is described for the Poisson problem but it is also applicable to the Stokes equations.
After the global problem given by equation~\eqref{eq:globalSystemPoissonModal} is solved, two approximations of the primal variable are computed using the same hybrid variable $\hat{\text{u}}$. 
On the one hand, a first approximation of the primal variable, $\ue$, is computed by solving the local, cell-by-cell, problem of equation~\eqref{eq:HDG-Poisson-DlocalModal}. 
On the other hand, a second approximation $\ue^\star$ is obtained by solving an extra local, cell-by-cell, problem corresponding to the first-order FCFV, first presented in~\cite{FCFV2018}, namely
\begin{equation}
\ue^\star = \alpha_e^{-1} \beta_e + \alpha_e^{-1} \sum_{j\in\Bset} |\Gamma_{e,j}| \tau_j \uHj,
\label{eq:Poisson-locU}
\end{equation}
where
\begin{equation} \label{eq:poissonPrecompK0}
\alpha_e := \sum_{j\in\Aset} |\Gamma_{e,j}| \tau_j 
, \quad
\beta_e := |\Omega_e| s_e  + \sum_{j\in\Dset} |\Gamma_{e,j}| \tau_j u_{D,j} ,
\end{equation}
and $\Aset$ denotes the set of all faces of cell $\Omega_e$.

The local error indicator for the cell $\Omega_e$ is thus defined as
\begin{equation} \label{eq:errorMeasure}
E_e := \left[ \frac{1}{|\Omega_e|} \int_{\Omega_e} \left( \ue - \ue^\star \right)^2 d\Omega \right]^{1/2} ,
\end{equation}
and it is employed to devise an automatic mesh adaptivity process. 
First, the following \emph{a priori} local error estimate for elliptic problems is recalled~\cite{diez1999unified,giorgiani2014hybridizable,HDGNEFEMstokes,sevilla2019hdg}
\begin{equation} \label{eq:aPrioriError}
\varepsilon_e := \| u^{\text{ex}} - u^h \|_{\eltwo(\Omega_e)} \leq C h_e^{1 + \nsd/2},
\end{equation}
where $u^{\text{ex}}$ and $u^h$ are the exact solution and its constant approximation in the cell $\Omega_e$ respectively, $h_e$ is the characteristic cell size and $C$ is an unknown constant. Then, using a classical Richardson extrapolation, the desired cell size is calculated as
\begin{equation} \label{eq:newElemSize}
h_e^\star = h_e \left(\frac{\varepsilon}{E_e}\right)^{1/\left(1 + \nsd/2 \right)}.
\end{equation}
where $\varepsilon$ is the user defined target error in each cell.
It is worth noting that the corresponding formula in~\cite{FCFV2}, namely equation (47), presents a typo in the exponent of $\varepsilon/E_e$.

\begin{remark}
The accuracy of the error indicator~\eqref{eq:errorMeasure} is guaranteed because $\ue$ is a second-order approximation of the solution, whereas $\ue^\star$ converges with first-order only. 
Moreover, this indicator provides information about the error between the solution $\ue^\star$ and the exact solution but, as it will be shown in the numerical examples, it also provides information about the error between solution $\ue$ and the exact solution.
\end{remark}

To illustrate the efficiency of the proposed strategy, table~\ref{tab:operationsErrorPoisson} reports the number of operations required by the error indicator in~\cite{FCFV2} and by the new error indicator proposed in this work for the Poisson problem, and using different cell types. Similarly, table~\ref{tab:operationsErrorStokes} show the corresponding number of operations for the Stokes problem. 
The tables display the number of operations required in each case to compute the error indicator for a mesh of $1,000$ cells of different types.
Concerning the indicator in~\cite{FCFV2}, its cost is given by the operations required to assemble and solve the extra global problem and to compute the corresponding $\ue$ for each cell using the new hybrid variable. 
The strategy proposed here only involves the extra computation of equation~\eqref{eq:Poisson-locU} for each cell.
\begin{table}[hbt]
	\centering
	\begin{tabular}[hbt]{| c || c | c | c | c | c | c |}
		\hline
		Method & Triangle & Quadrilateral & Tetrahedron & Hexahedron & Prism & Pyramid \\
		\hline & & & & & &
		\\ [-1em] 
		\hline
		Error indicator in~\cite{FCFV2} & $5.6 \times 10^8$ & $1.3 \times 10^9$ & $1.3 \times 10^9$ & $4.5 \times 10^9$  & $2.6 \times 10^9$ & $2.6 \times 10^9$ \\
		\hline
		Error indicator~\eqref{eq:errorMeasure} & $1.0 \times 10^4$ & $1.3 \times 10^4$ & $1.3 \times 10^4$ & $1.9 \times 10^4$ & $1.6 \times 10^4$ & $1.6 \times 10^4$ \\
		\hline
	\end{tabular}
	\caption{Number of operations required to compute the error indicator of the Poisson problem for a mesh with 1,000 cells of different types.}
	\label{tab:operationsErrorPoisson}
\end{table}
\begin{table}[hbt]
	\centering
	\begin{tabular}[hbt]{| c || c | c | c | c | c | c |}
		\hline
		Method & Triangle & Quadrilateral & Tetrahedron & Hexahedron & Prism & Pyramid \\
		\hline & & & & & &
		\\ [-1em] 
		\hline
		Error indicator in~\cite{FCFV2} & $1.1 \times 10^{10}$ & $2.1 \times 10^{10}$ & $5.7 \times 10^{10}$ & $1.7 \times 10^{11}$  & $1.0 \times 10^{11}$ & $1.0 \times 10^{11}$ \\
		\hline
		Error indicator~\eqref{eq:errorMeasure} & $1.7 \times 10^4$ & $2.2 \times 10^4$ & $3.1 \times 10^4$ & $4.5 \times 10^4$ & $3.8 \times 10^4$ & $3.8 \times 10^4$  \\
		\hline
	\end{tabular}
	\caption{Number of operations required to compute the error indicator of the Stokes problem for a mesh with 1,000 cells of different types.}
	\label{tab:operationsErrorStokes}
\end{table}
The results show that the new error indicator is several orders of magnitude less expensive. In addition, the number of operations of the proposed error indicator scales linearly with the number of cells, whereas the number of operations required by the error indicator proposed in~\cite{FCFV2} scales with the cube of the number of cells, due to the required extra solution of a global linear system.

\section{Numerical studies}
\label{sc:studies}

This section presents an extensive set of experiments to numerically validate the optimal convergence properties of the proposed method and to show its robustness against the choice of the stabilisation parameter and the mesh properties such as cell distortion and stretching. The examples presented involve meshes of different cell types, namely triangular and quadrilateral cells in two dimensions, and tetrahedral, hexahedral, prismatic and pyramidal cells in three dimensions. In addition, results with hybrid meshes are presented for the first time in the context of the FCFV method.

The model problem considered for the Poisson equation involves the numerical solution of~\eqref{eq:PoissonBrokenFirstOrder} in $\Omega=[0,1]^{\nsd}$. In two dimensions, the source term and boundary data are selected such that the analytical solution is known and given by
\begin{equation}
u^{\text{ex}}(x_1,x_2) = \exp \big( \alpha \sin(a x_1 + c x_2) + \beta\cos(b x_1 + d x_2) \big),
\end{equation}
with $\alpha=0.1$, $\beta=0.3$, $a=5.1$, $b=4.3$, $c=-6.2$ and $d=3.4$. Neumann boundary conditions are imposed on the bottom part of the boundary, on $\Gamma_N = \{(x_1,x_2) \in \mathbb{R}^2 \; | \; x_2=0\}$, and Dirichlet boundary conditions are set on the rest of the boundary. 
For the three dimensional Poisson problem, the source term and boundary data are selected such that the analytical solution is 
\begin{equation}
u^{\text{ex}}(x_1,x_2,x_3) = \exp \big( \alpha \sin(a x_1 + c x_2 + e x_3) + \beta\cos(b x_1 + d x_2 + f x_3) \big),
\end{equation}
with $\alpha=0.1$, $\beta=0.3$, $a=5.1$, $b=4.3$, $c=-6.2$, $d=3.4$, $e=1.8$ and $f=1.7$. Neumann boundary conditions are imposed on $\Gamma_N = \{(x_1,x_2,x_3) \in \mathbb{R}^3 \; | \; x_3=0\}$, whereas on the remaining boundary surfaces, Dirichlet conditions are enforced.

The domain $\Omega=[0,1]^{\nsd}$ is also utilised for the Stokes equations~\eqref{eq:StokesBroken} with viscosity $\nu=1$. For the two dimensional case, the source term and boundary conditions are devised in order for the analytical velocity and pressure fields to be
\begin{equation}
\left\{
\begin{aligned}
u^{\text{ex}}_1(x_1,x_2) & = x_1^2 (1-x_1)^2 (2x_2-6x_2^2+4x_2^3) , \\
u^{\text{ex}}_2(x_1,x_2) & = - x_2^2 (1-x_2)^2 (2x_1-6x_1^2+4x_1^3) , \\
p^{\text{ex}}(x_1,x_2) & = x_1 (1-x_1) .
\end{aligned}
\right.
\end{equation}
On $\Gamma_N = \{(x_1,x_2) \in \mathbb{R}^2 \; | \; x_2=0\}$, a Neumann condition representing a pseudo-traction is imposed, whereas the analytical velocity enforcing Dirichlet conditions is set on the rest of the boundary. 
Similarly, in three dimensions, Neumann boundary conditions are imposed on $\Gamma_N = \{(x_1,x_2,x_3) \in \mathbb{R}^3 \; | \; x_3=0\}$ and Dirichlet conditions on the rest of the boundary, to match the analytical expressions of velocity and pressure given by
\begin{equation}
\left\{
\begin{aligned}
u^{\text{ex}}_1(x_1,x_2,x_3) & = \tfrac{1}{2} + (x_3-x_2) \sin\left(x_1 - \tfrac{1}{2} \right) ,\\
u^{\text{ex}}_2(x_1,x_2,x_3) & = 1 - x_2 \left(x_3- \tfrac{1}{2} x_2 \right)\cos\left(x_1-\tfrac{1}{2} \right) - x_2 \left(x_1-\tfrac{1}{2}x_2 \right)\cos\left(x_3-\tfrac{1}{2} \right) , \\
u^{\text{ex}}_3(x_1,x_2,x_3) & = \tfrac{1}{2} + (x_1-x_2) \sin\left(x_3 -\tfrac{1}{2} \right) , \\
p^{\text{ex}}(x_1,x_2,x_3) & = x_1 (1-x_1) + x_2 (1-x_2) + x_3 (1-x_3).
\end{aligned}
\right.
\end{equation}

\begin{figure}[!tb]
	\centering
	\subfigure[Quadrilateral mesh]{\includegraphics[width=0.28\textwidth]{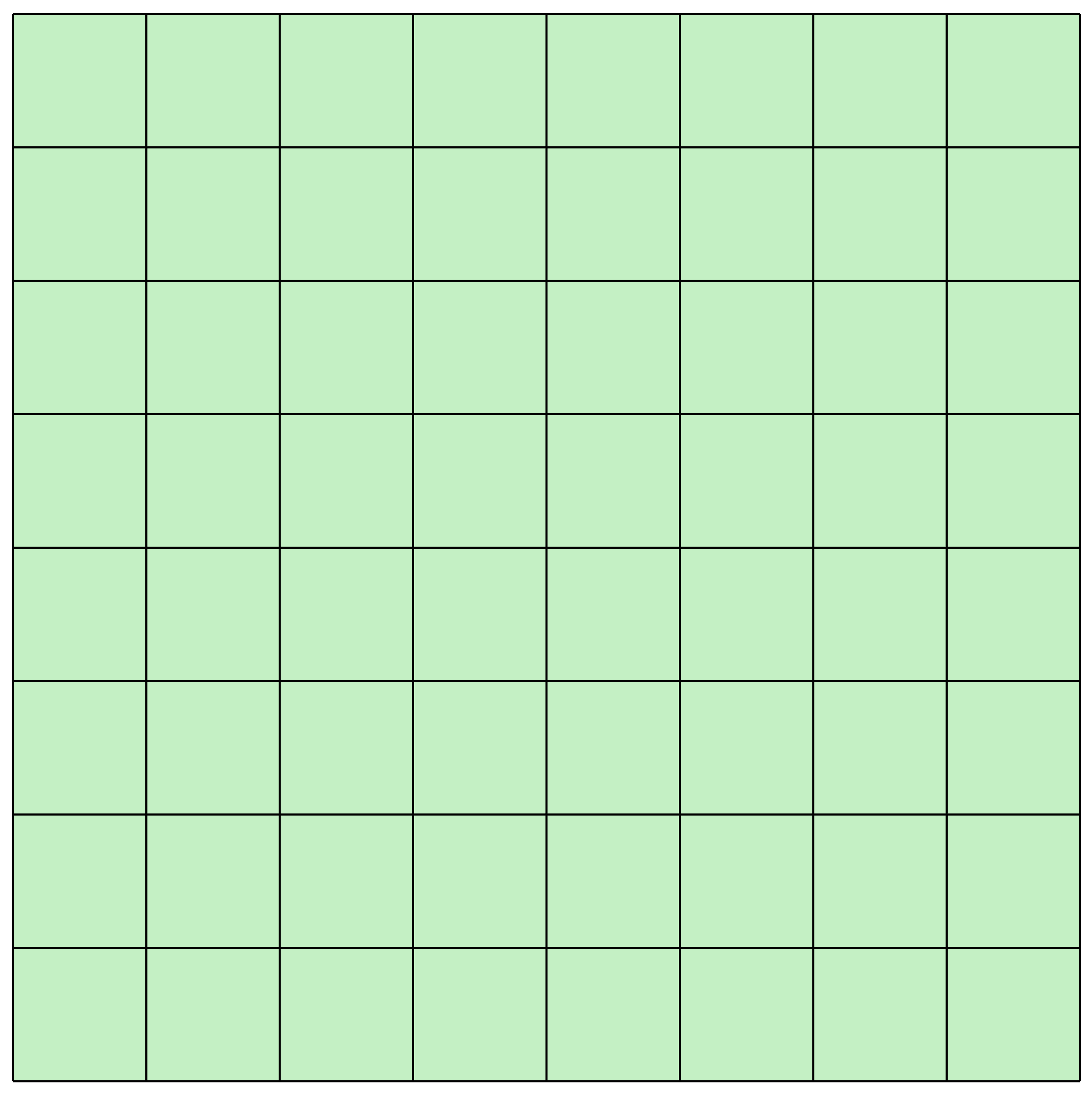}}
	\subfigure[Triangular mesh]{\includegraphics[width=0.28\textwidth]{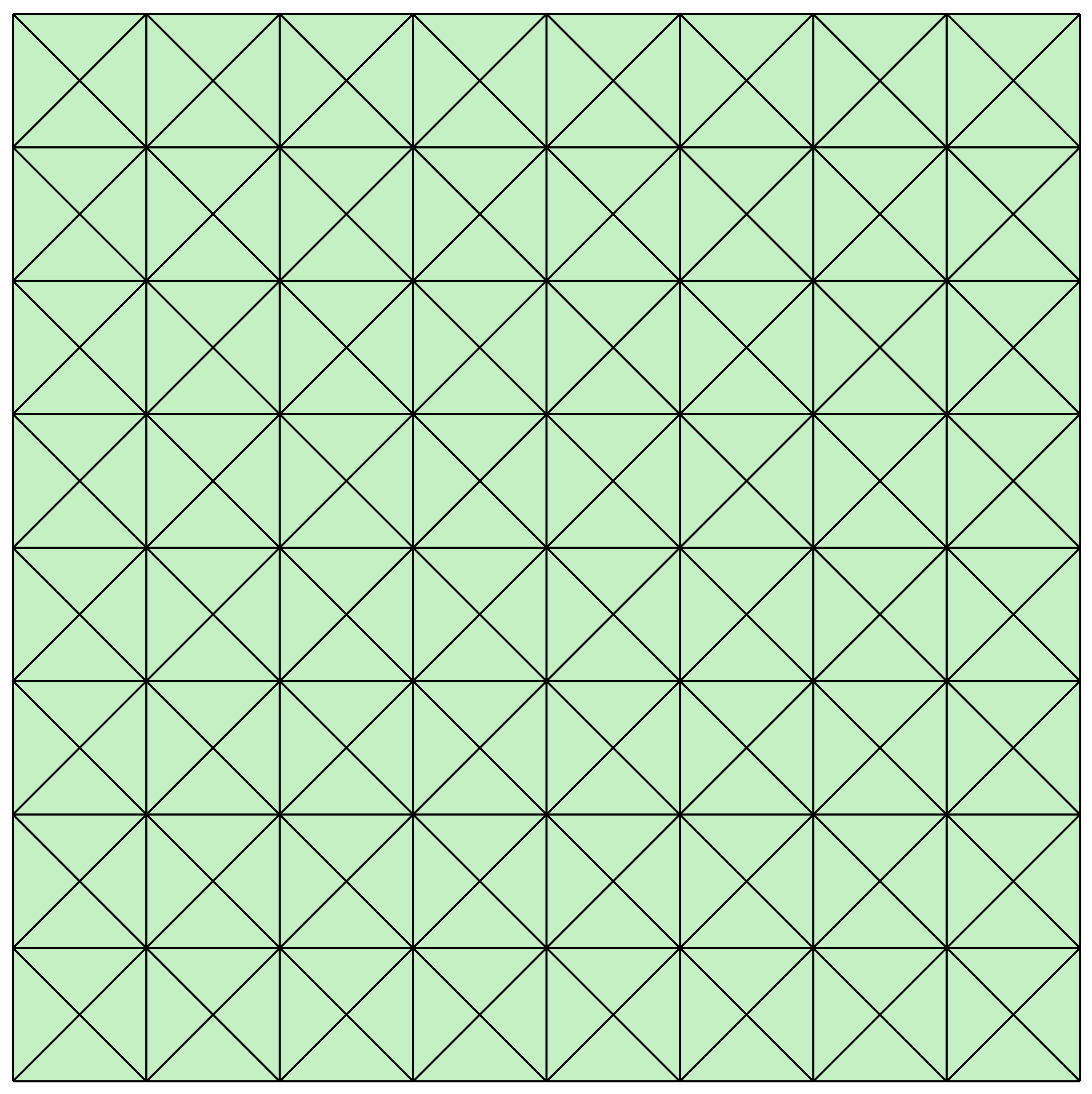}}
	\subfigure[Hybrid mesh]{\includegraphics[width=0.28\textwidth]{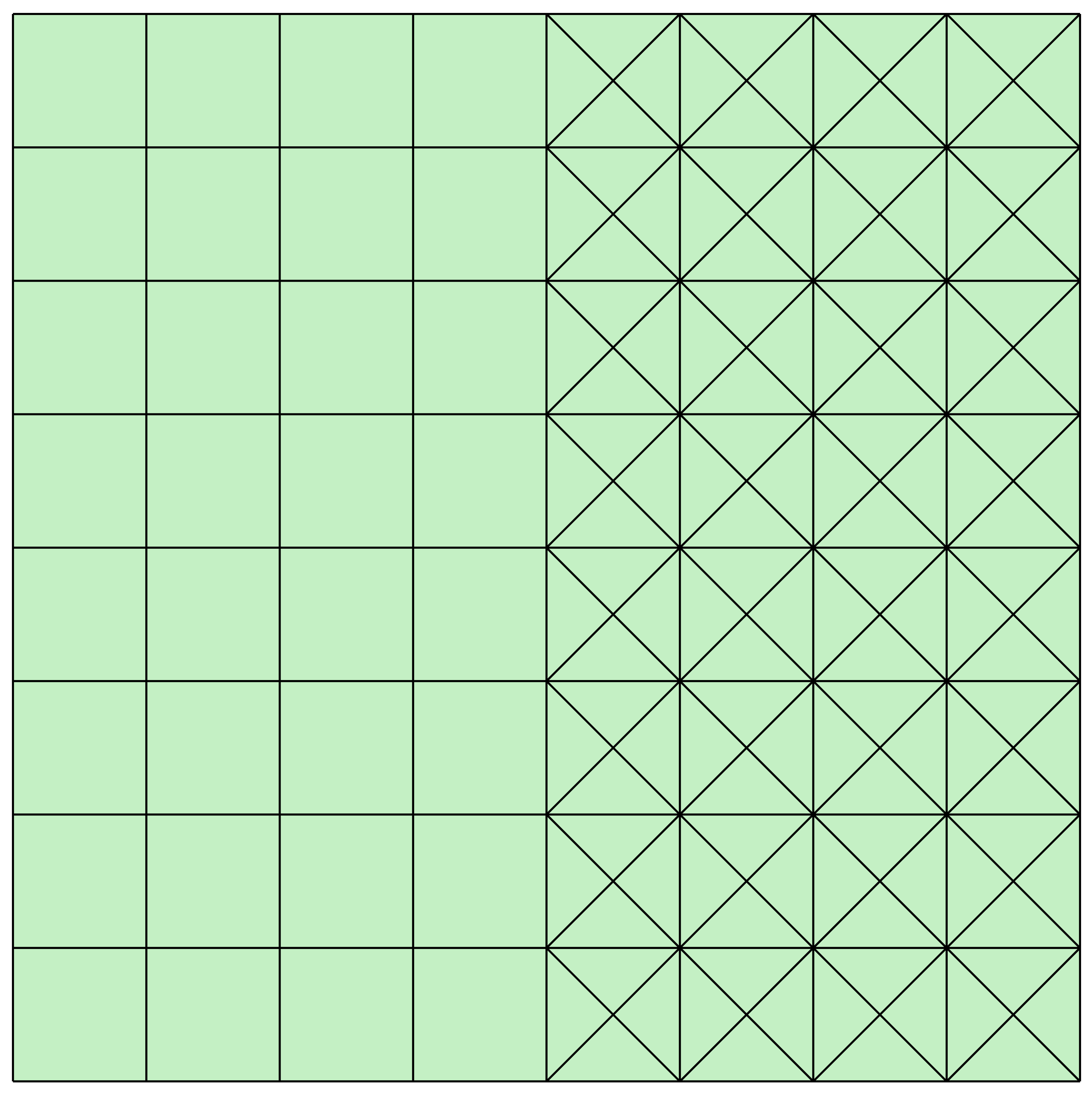}}	
	\caption{Meshes corresponding to the third level of refinement for the domain $\Omega=[0,1]^2$ using (a) quadrilateral, (b) triangular and (c) hybrid cells.}
	\label{fig:2Dmeshes}
\end{figure}

To perform the mesh convergence study in two dimensions, a set of eight uniform triangular and quadrilateral meshes are generated. The meshes corresponding to the third level of refinement are displayed in Figures~\ref{fig:2Dmeshes} (a) and (b). 
To demonstrate the flexibility of the proposed method, hybrid meshes made of quadrilateral and triangular cells are also considered. The third hybrid mesh is displayed in Figure~\ref{fig:2Dmeshes} (c).

\begin{figure}[!tb]
	\centering
	\subfigure[Hexahedral mesh] {\includegraphics[width=0.19\textwidth]{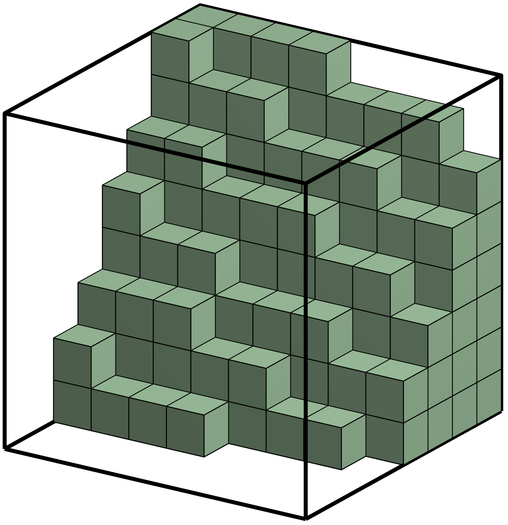}}
	\subfigure[Tetrahedral mesh]{\includegraphics[width=0.19\textwidth]{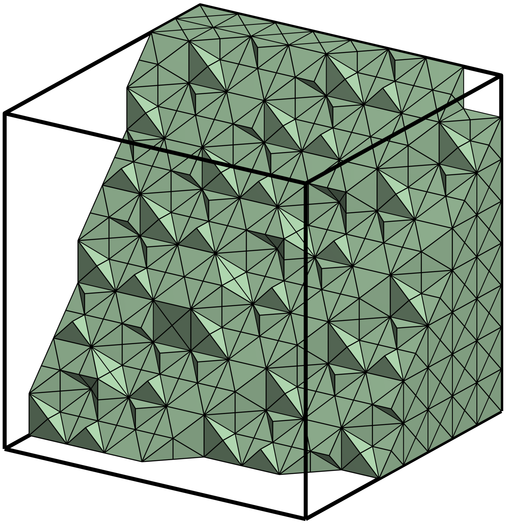}}
	\subfigure[Prismatic mesh]  {\includegraphics[width=0.19\textwidth]{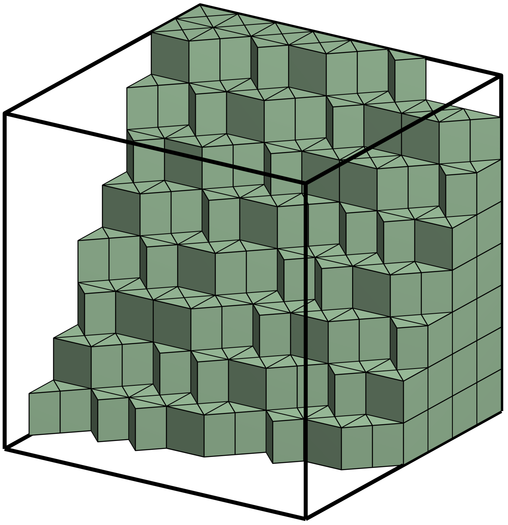}}	
	\subfigure[Pyramidal mesh]  {\includegraphics[width=0.19\textwidth]{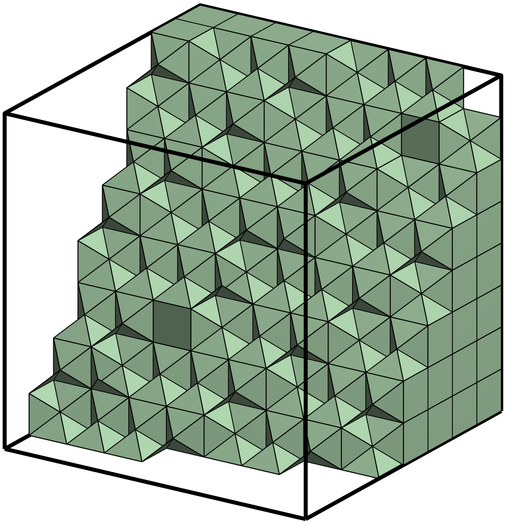}}	\subfigure[Hybrid mesh]     {\includegraphics[width=0.19\textwidth]{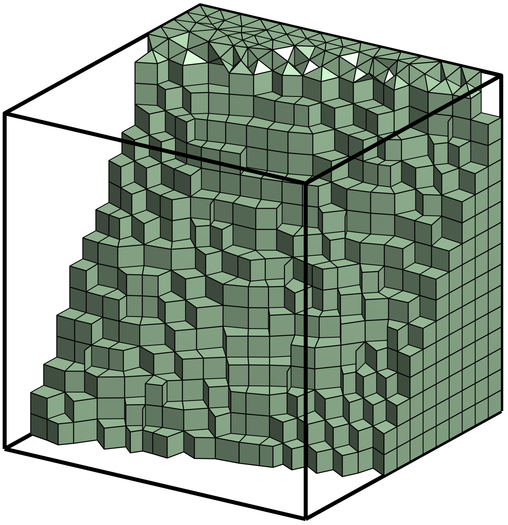}}	
	\caption{Internal view of the meshes corresponding to the fourth level of refinement for the domain $\Omega=[0,1]^3$ featuring (a) hexahedral, (b) tetrahedral, (c) prismatic, (d) pyramidal and (e) hybrid cells.}
	\label{fig:3Dmeshes}
\end{figure}

In three dimensions, a set of six uniform tetrahedral, hexahedral, prismatic and pyramidal meshes are generated. In addition, hybrid meshes containing a mixture of different cell types are also considered. Figure~\ref{fig:3Dmeshes} shows the fourth level of mesh refinement for the different types of meshes considered and table~\ref{tab:hybridMeshes3D} shows the statistics for the hybrid meshes.
\begin{table}[hbt]
\centering
\begin{tabular}[hbt]{| c || c | c | c | c | c |}
	\hline
	Mesh & Number of cells & Hexahedral cells & Tetrahedral cells & Pyramidal cells & $h$ \\
	\hline & & & &
	\\ [-1em] 
	\hline
	1 & 8 & 1 & 2 & 5 & 1.414 \\
	\hline
	2 & 140 & 42 & 28 & 70 & 0.413 \\
	\hline
	3 & 868 & 434 & 124 & 310 & 0.209 \\
	\hline
	4 & 6,226 & 4,245 & 566 & 1,415 & 0.104 \\
	\hline
	5 & 20,340 & 15,594 & 1,356 & 3,390 & 0.071 \\
	\hline
	6 & 64,638 & 53,865 & 3,078 & 7,695 & 0.047 \\
	\hline
\end{tabular}
\caption{Details of the six hybrid meshes for the convergence study in three dimensions.}
\label{tab:hybridMeshes3D}
\end{table}

\subsection{Optimal convergence of the FCFV method}
\label{sc:meshConvergence}

The first experiment involves a mesh convergence study for the Poisson and Stokes problems both in two and three dimensions and using meshes with different cell types. In all cases the stabilisation parameter is selected to be $\tau = 10^4$ in two dimensions and $\tau=10^2$ for three dimensional problems. A detailed study of the influence of the stabilisation parameter on the accuracy of the proposed FCFV method is provided in section~\ref{sc:InfluenceTau}.

\begin{figure}[!bt]
	\centering
	\subfigure[$u$]{\includegraphics[width=0.32\textwidth]{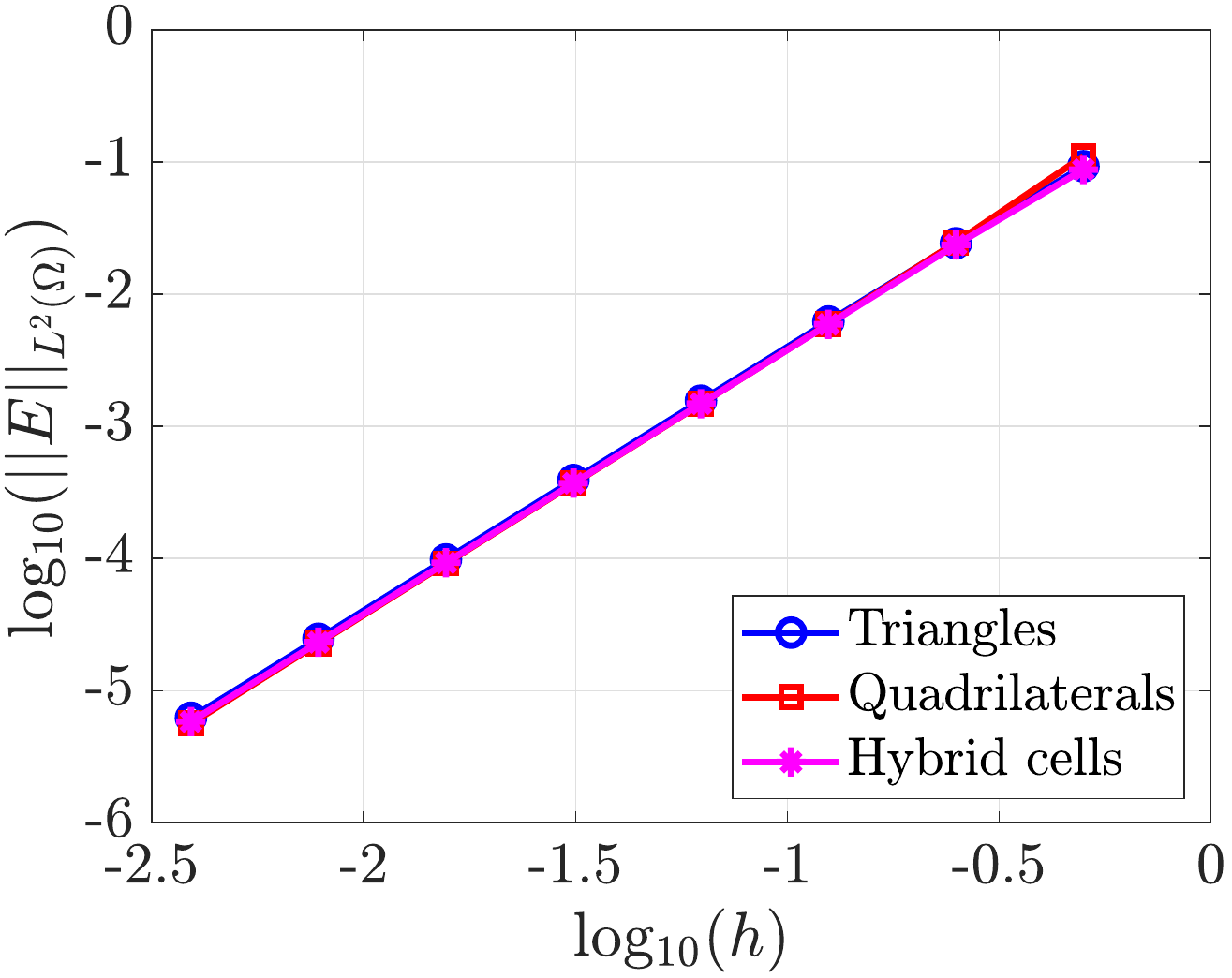}}
	\subfigure[$\bq$]{\includegraphics[width=0.32\textwidth]{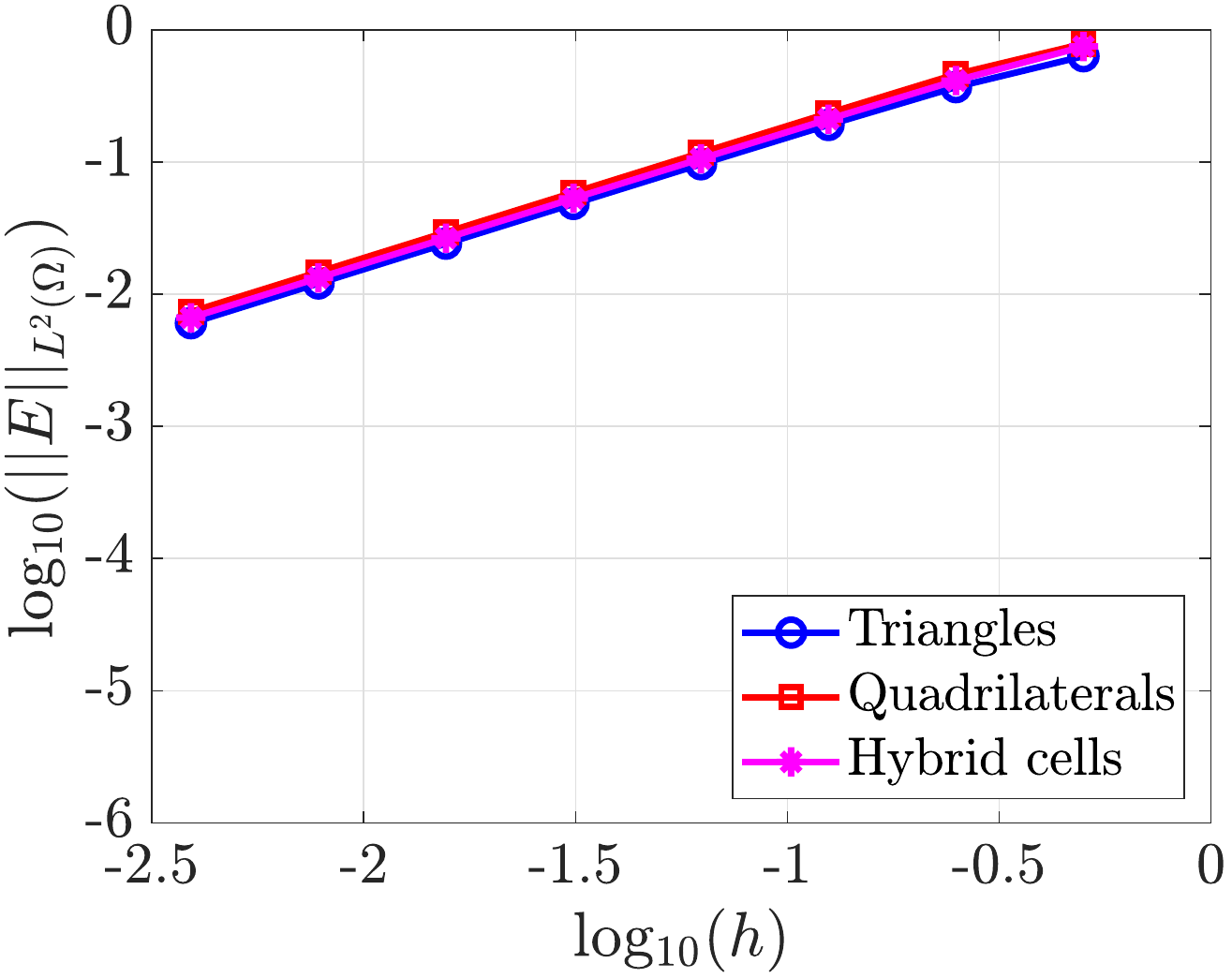}}
	\caption{Mesh convergence of the error of the solution $u$ and its gradient $\bq$ in the $\eltwo(\Omega)$ norm as a function of the cell size $h$ for two dimensional Poisson problem on regular meshes using different cell types.}
	\label{fig:Poisson_Conv_2D}
\end{figure}

For the two dimensional Poisson problem, Figure~\ref{fig:Poisson_Conv_2D} shows the relative error, measured in the $\eltwo(\Omega)$ norm, of the primal and mixed variables as a function of the characteristic cell size. 
The results show optimal, quadratic, convergence of the primal variable for triangular, quadrilateral and hybrid meshes with almost identical accuracy in the three cases. For the mixed variable an optimal, linear, convergence is also observed, again, with almost identical accuracy in the three cases.

\begin{figure}[!bt]
	\centering
	\subfigure[$u$]{\includegraphics[width=0.32\textwidth]{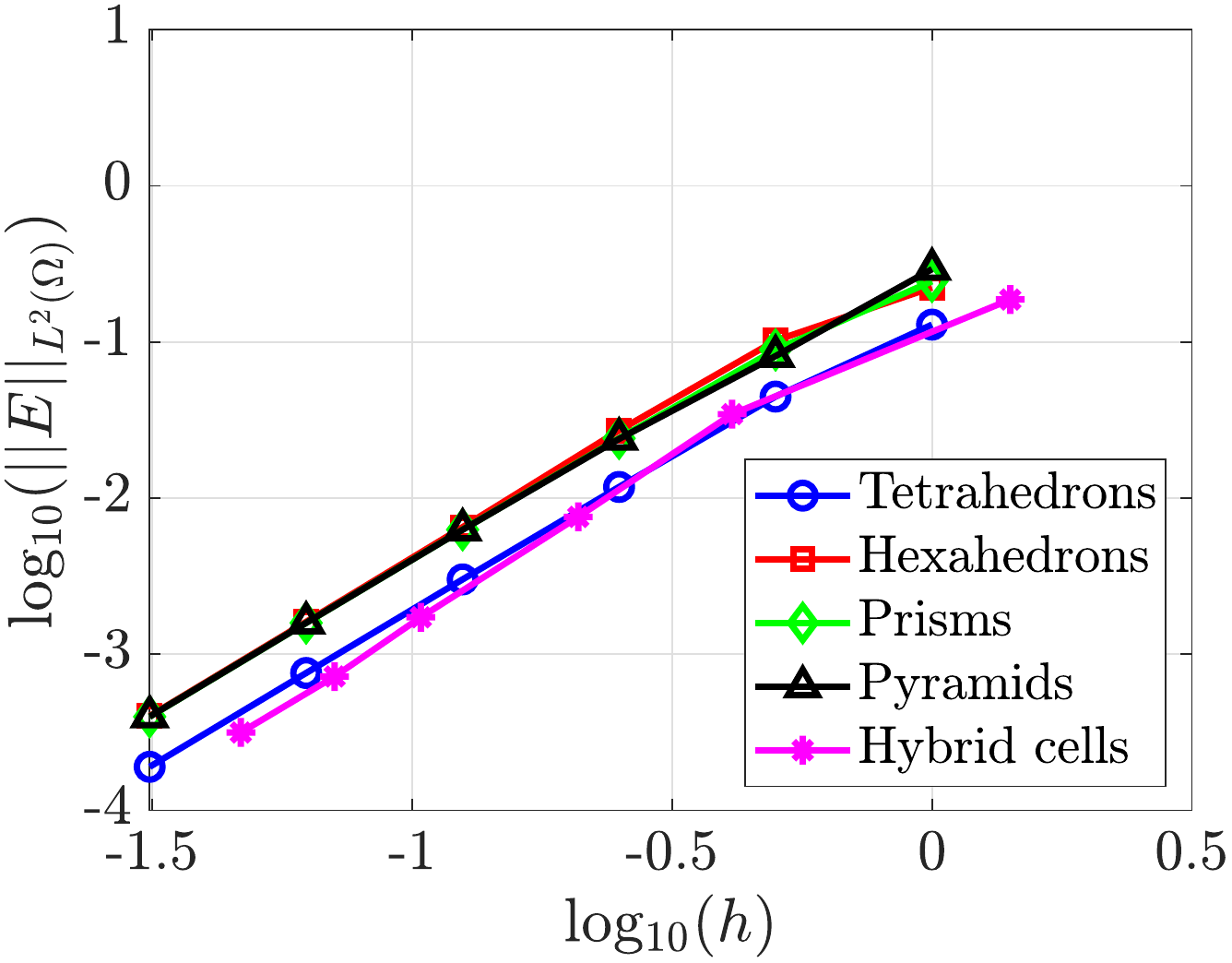}}
	\subfigure[$\bq$]{\includegraphics[width=0.32\textwidth]{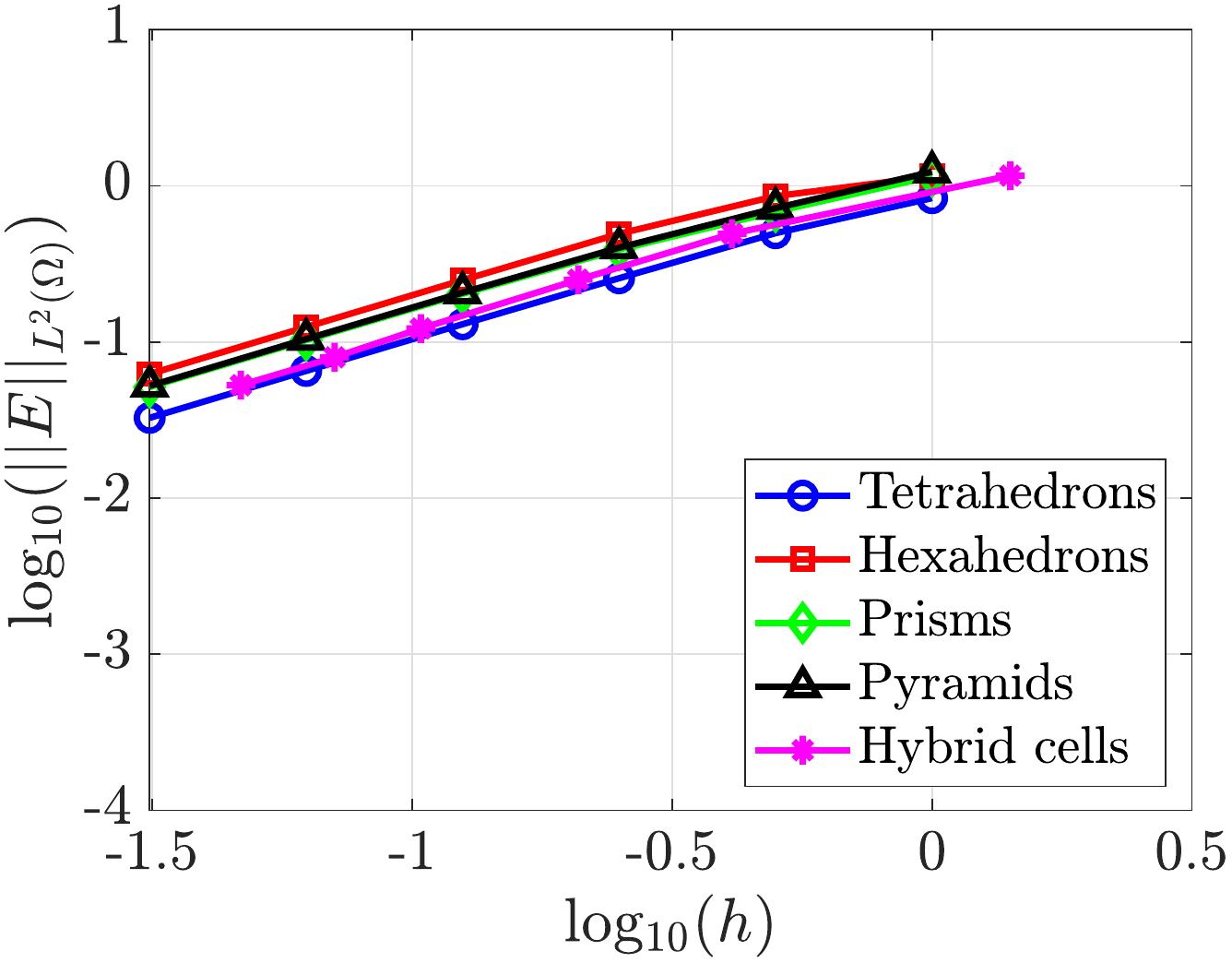}}
	\caption{Mesh convergence of the error of the solution $u$ and its gradient $\bq$ in the $\eltwo(\Omega)$ norm as a function of the cell size $h$ for three dimensional Poisson problem on regular meshes using different cell types.}
	\label{fig:Poisson_Conv_3D}
\end{figure}

The same mesh convergence study is performed for the three dimensional Poisson problem.  Figure~\ref{fig:Poisson_Conv_3D} shows the relative error of the primal and mixed variables, measured in the $\eltwo(\Omega)$ norm, as a function of the characteristic cell size. 
The results display a similar qualitative behaviour when compared to the two dimensional case. An optimal, quadratic, convergence is observed for the primal variable and an optimal, linear, rate of convergence is observed for the mixed variable, for all cell types and hybrid meshes. Tetrahedral and hybrid meshes provide slightly more accurate results when compared to hexahedral, prismatic and pyramidal meshes.

\begin{figure}[!bt]
	\centering
	\subfigure[$\bu$]{\includegraphics[width=0.32\textwidth]{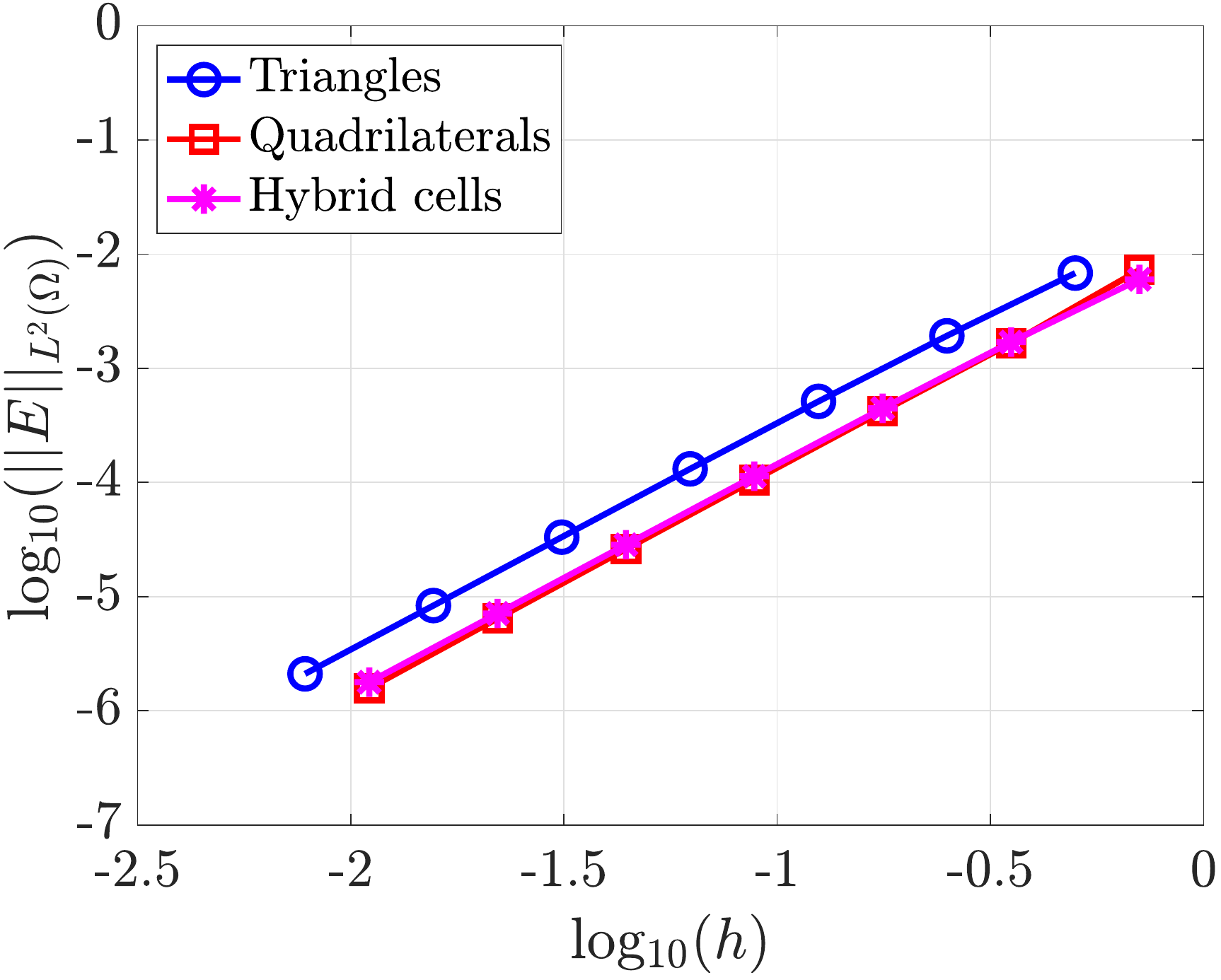}}
	\subfigure[$p$]{\includegraphics[width=0.32\textwidth]{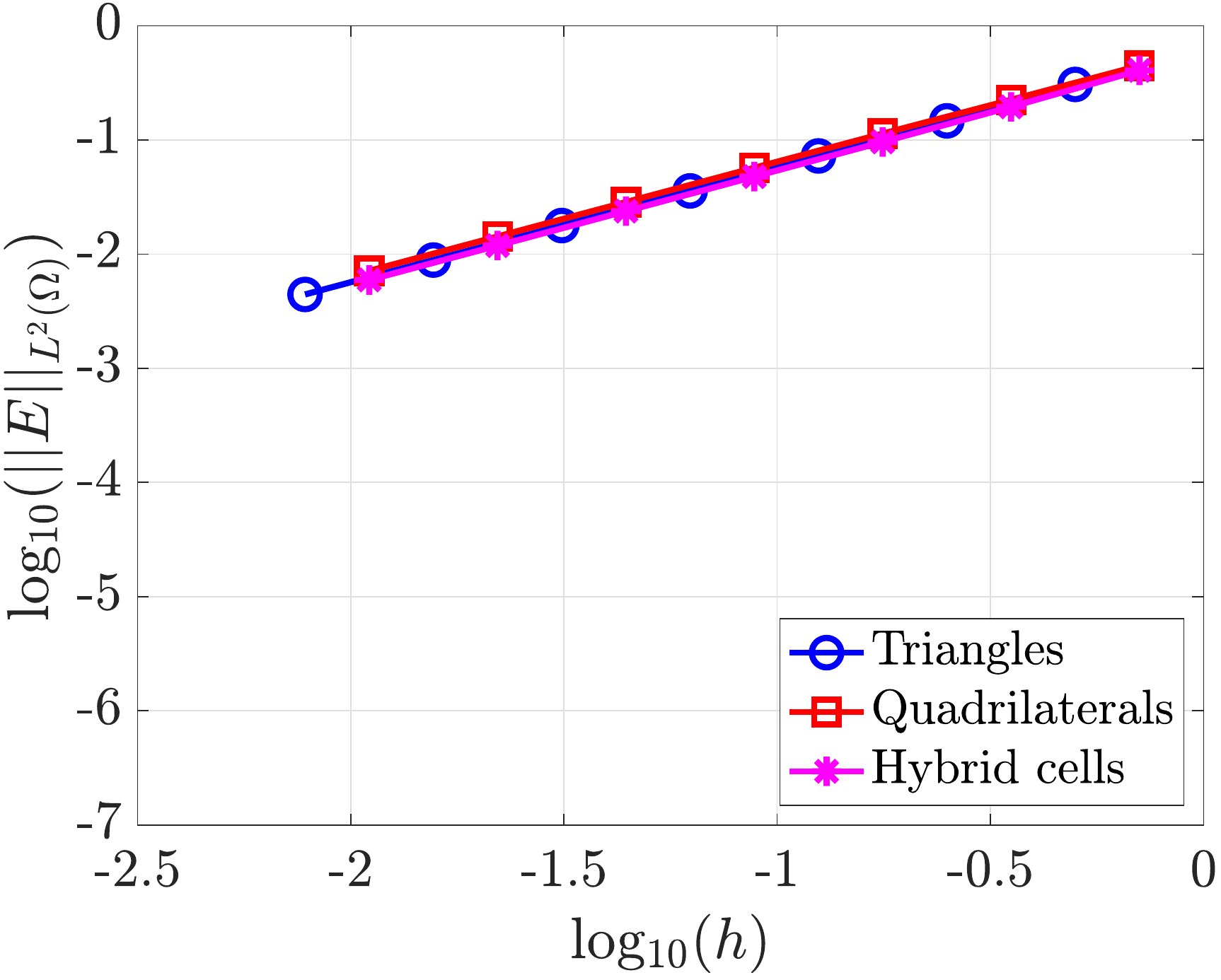}}
	\subfigure[$\bL$]{\includegraphics[width=0.32\textwidth]{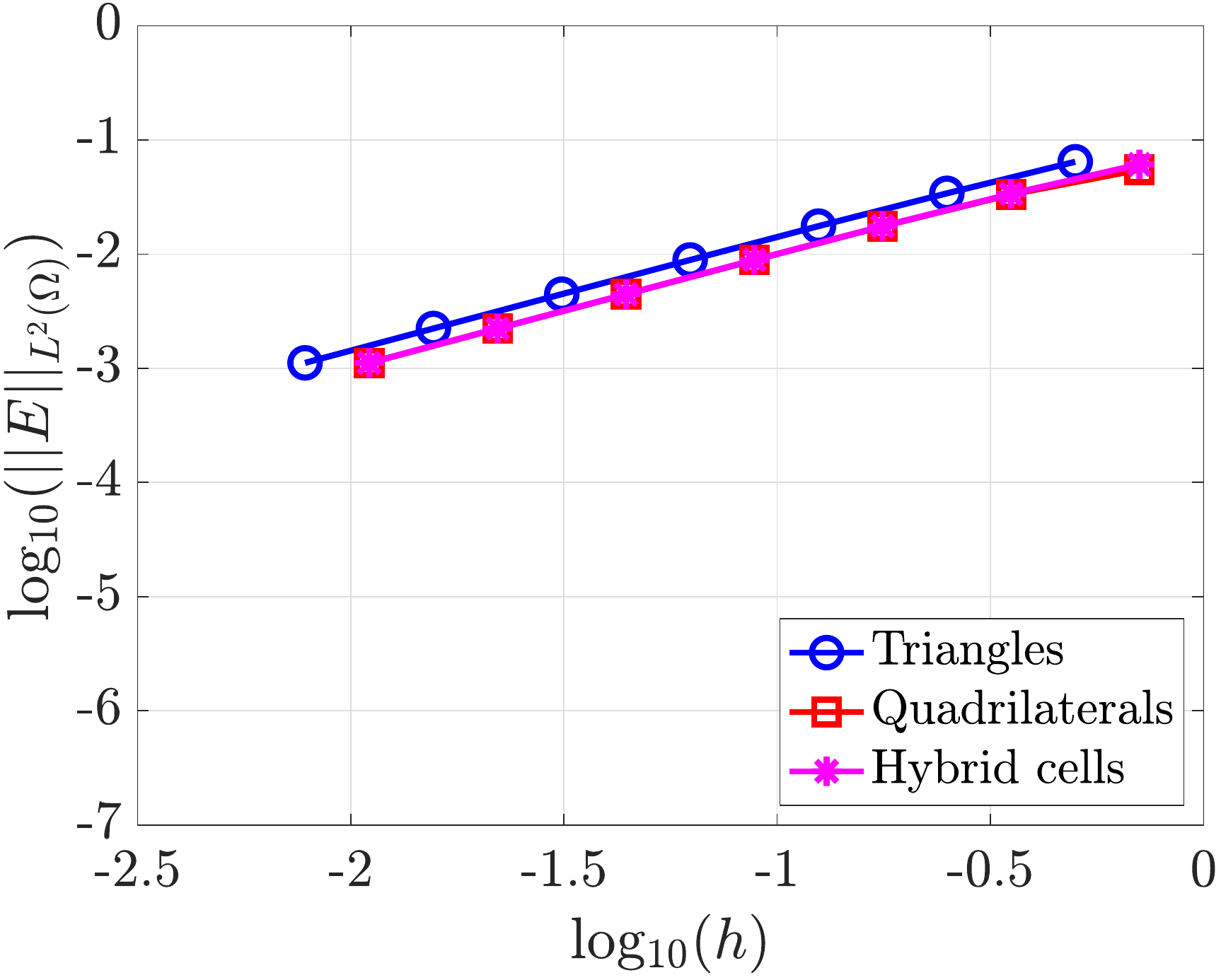}}
	\caption{Mesh convergence of the error of velocity $\bu$, pressure $p$ and gradient of velocity $\bL$ in the $\eltwo(\Omega)$ norm as a function of the cell size $h$ for two dimensional Stokes problem on regular meshes using different cell types.}
	\label{fig:Stokes_Conv_2D}
\end{figure}

\begin{figure}[!bt]
	\centering
	\subfigure[$\bu$]{\includegraphics[width=0.32\textwidth]{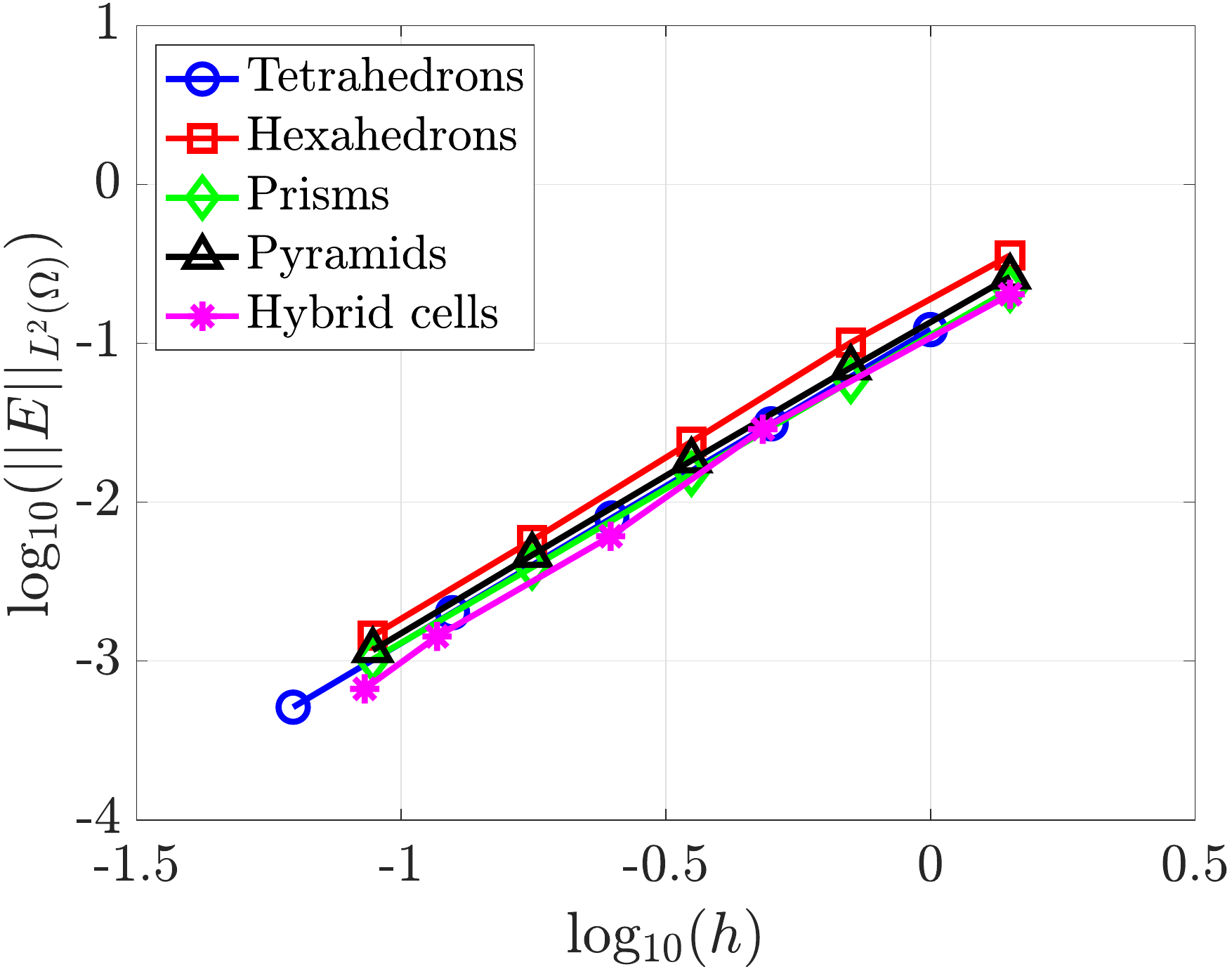}}
	\subfigure[$p$]{\includegraphics[width=0.32\textwidth]{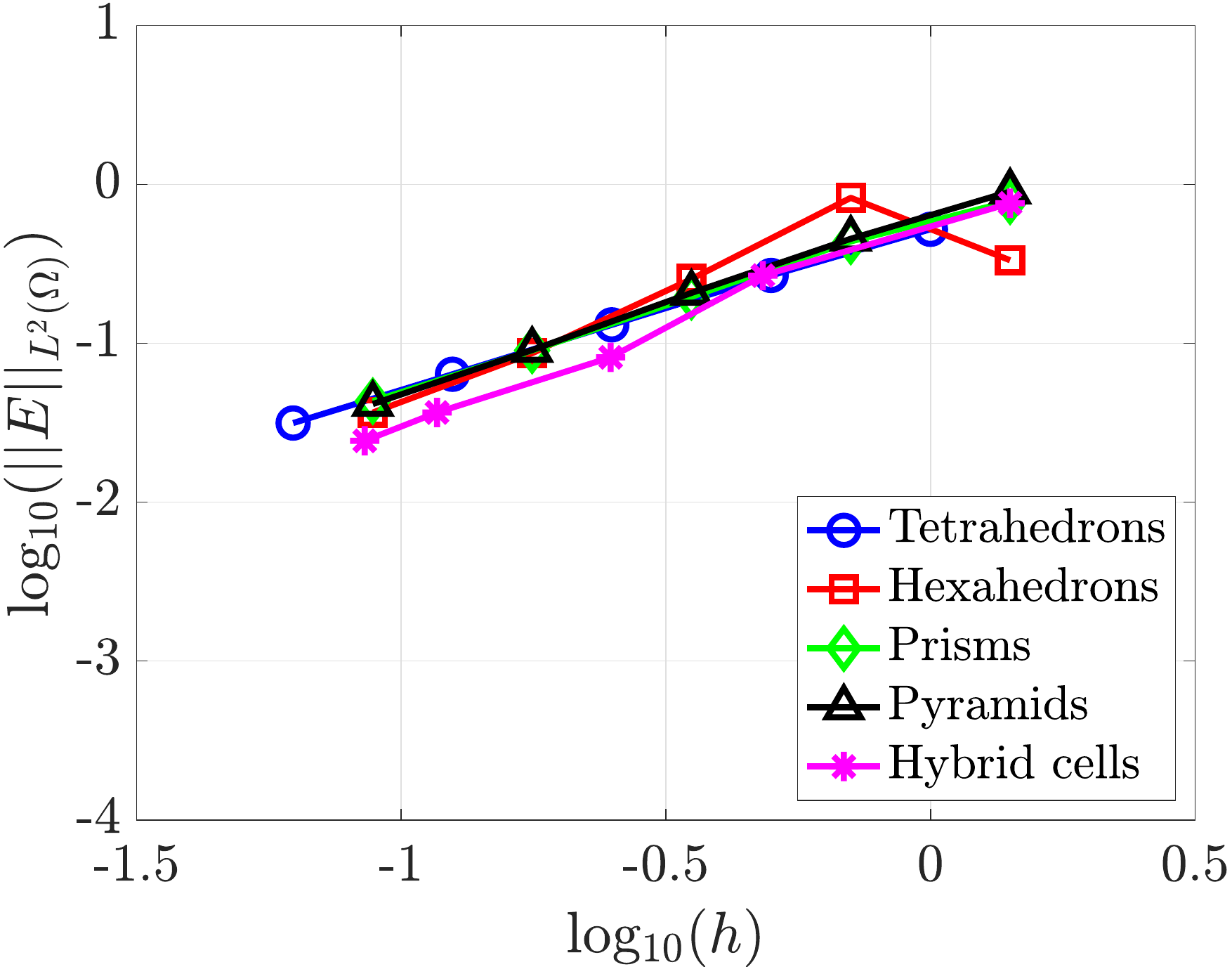}}
	\subfigure[$\bL$]{\includegraphics[width=0.32\textwidth]{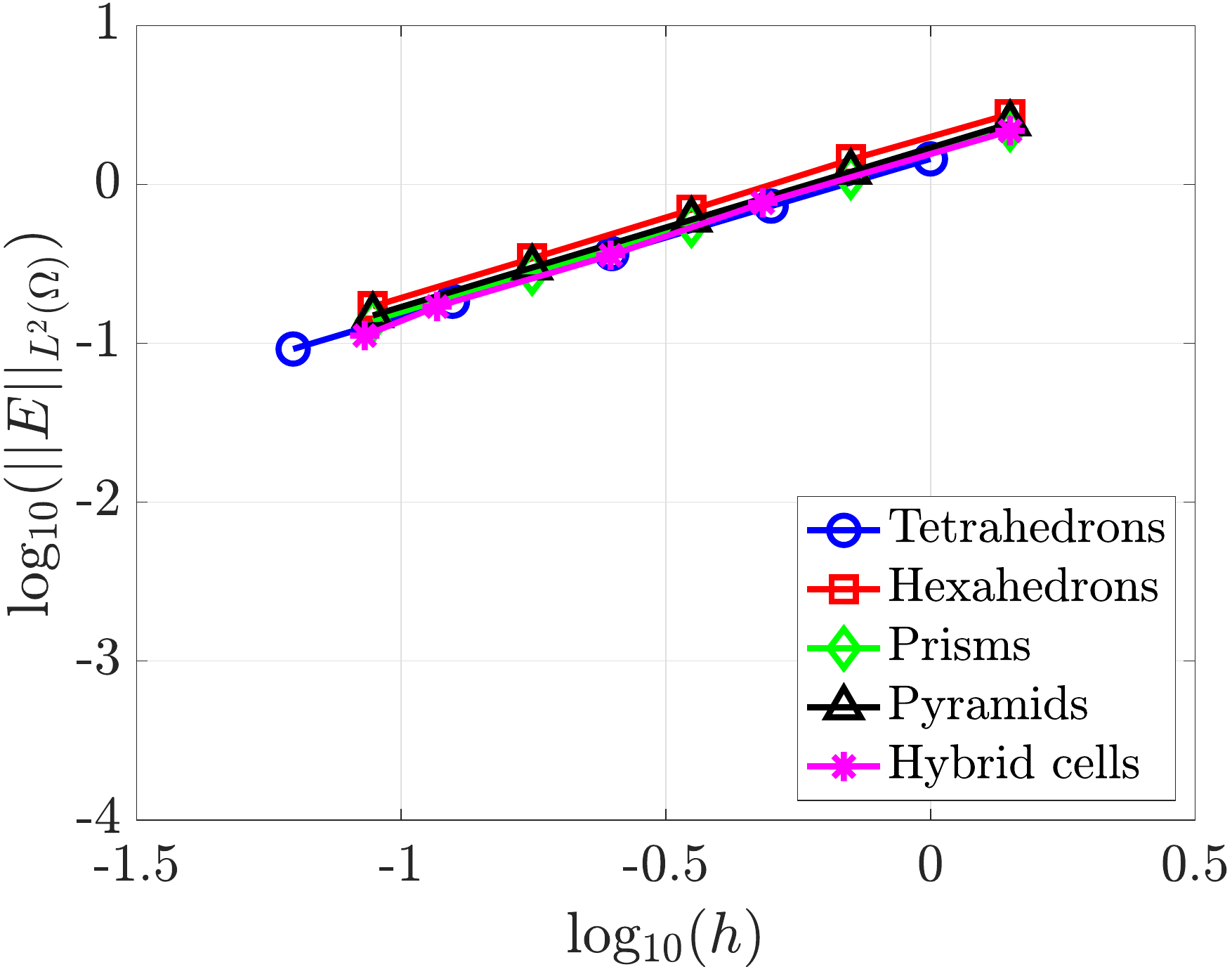}}
	\caption{Mesh convergence of the error of velocity $\bu$, pressure $p$ and gradient of velocity $\bL$ in the $\eltwo(\Omega)$ norm as a function of the cell size $h$ for three dimensional Stokes problem on regular meshes using different cell types.}
	\label{fig:Stokes_Conv_3D}
\end{figure}

Next, the mesh convergence study is performed for the Stokes problem in two and three dimensions. Figure~\ref{fig:Stokes_Conv_2D}  displays the relative error of velocity, pressure and gradient of velocity, measured in the $\eltwo(\Omega)$ norm, as a function of the characteristic cell size. 
The results show again an optimal quadratic convergence of the error of the velocity for all the different types of meshes. For the pressure and the gradient of the velocity, optimal, linear, convergence of the error is also observed for all types of meshes.
The same conclusions are observed from the results of the Stokes problem in three dimensions, displayed in Figure~\ref{fig:Stokes_Conv_3D}.

\subsection{Influence of the stabilisation parameter}
\label{sc:InfluenceTau}

The influence of the stabilisation parameter $\tau$ is studied numerically. The results in this section only consider the Stokes problem as further numerical examples, not reported here for brevity, have shown that identical conclusions are obtained for the Poisson problem.

\begin{figure}[!bt]
	\centering
	\subfigure[Triangular cells]{\includegraphics[width=0.32\textwidth]{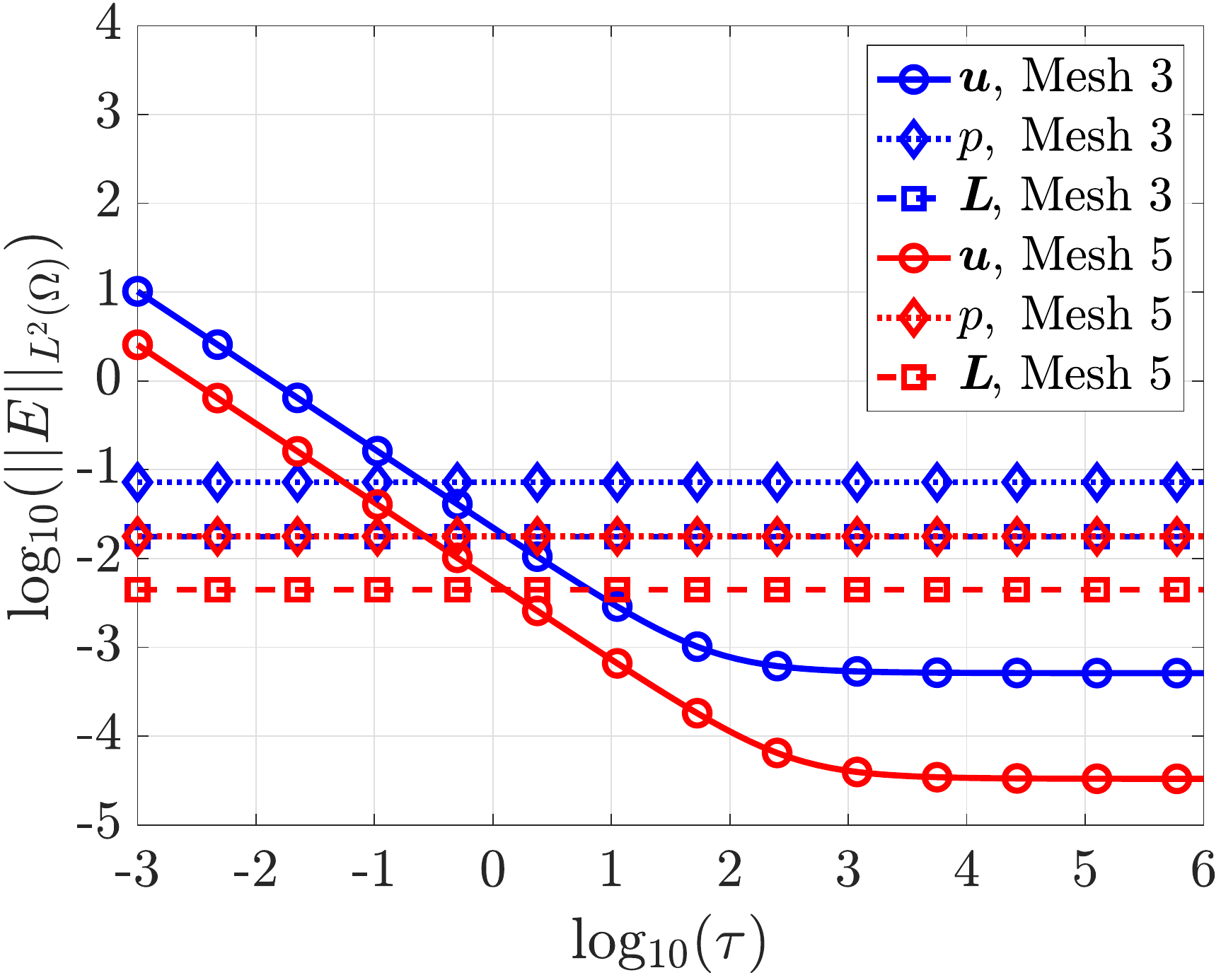}}
	\subfigure[Quadrilateral cells]{\includegraphics[width=0.32\textwidth]{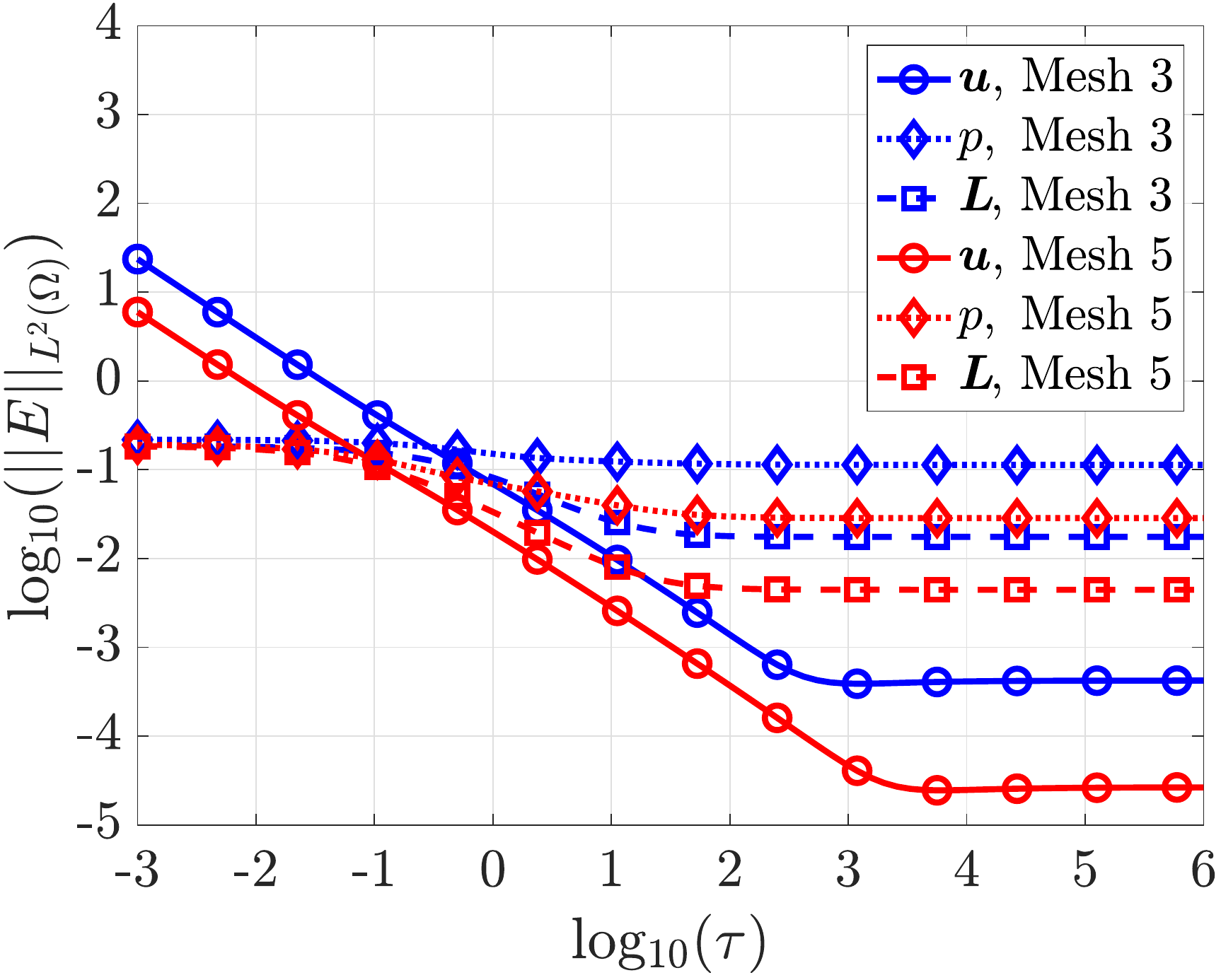}}
	\subfigure[Tetrahedral cells]{\includegraphics[width=0.32\textwidth]{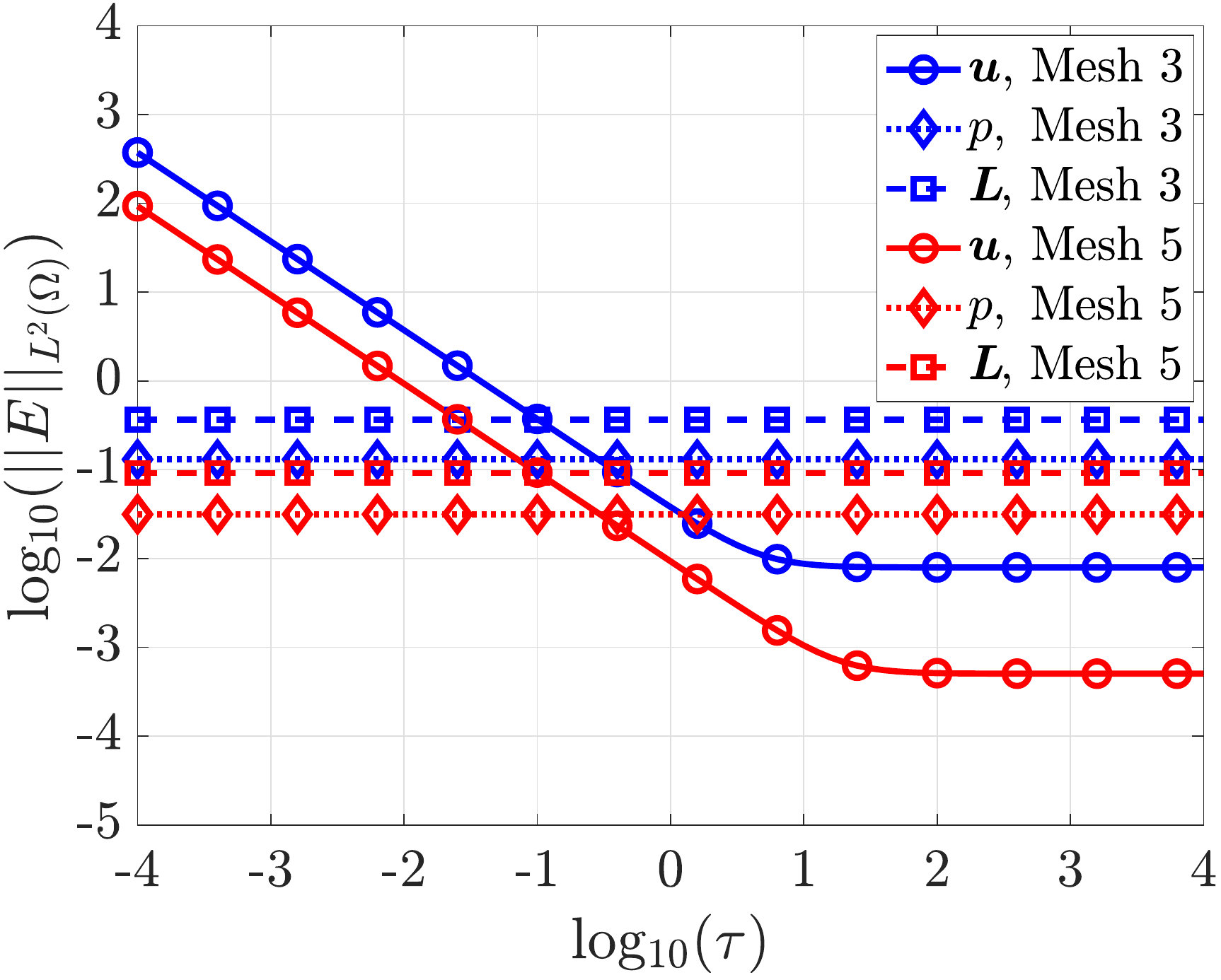}}
	\subfigure[Hexahedral cells]{\includegraphics[width=0.32\textwidth]{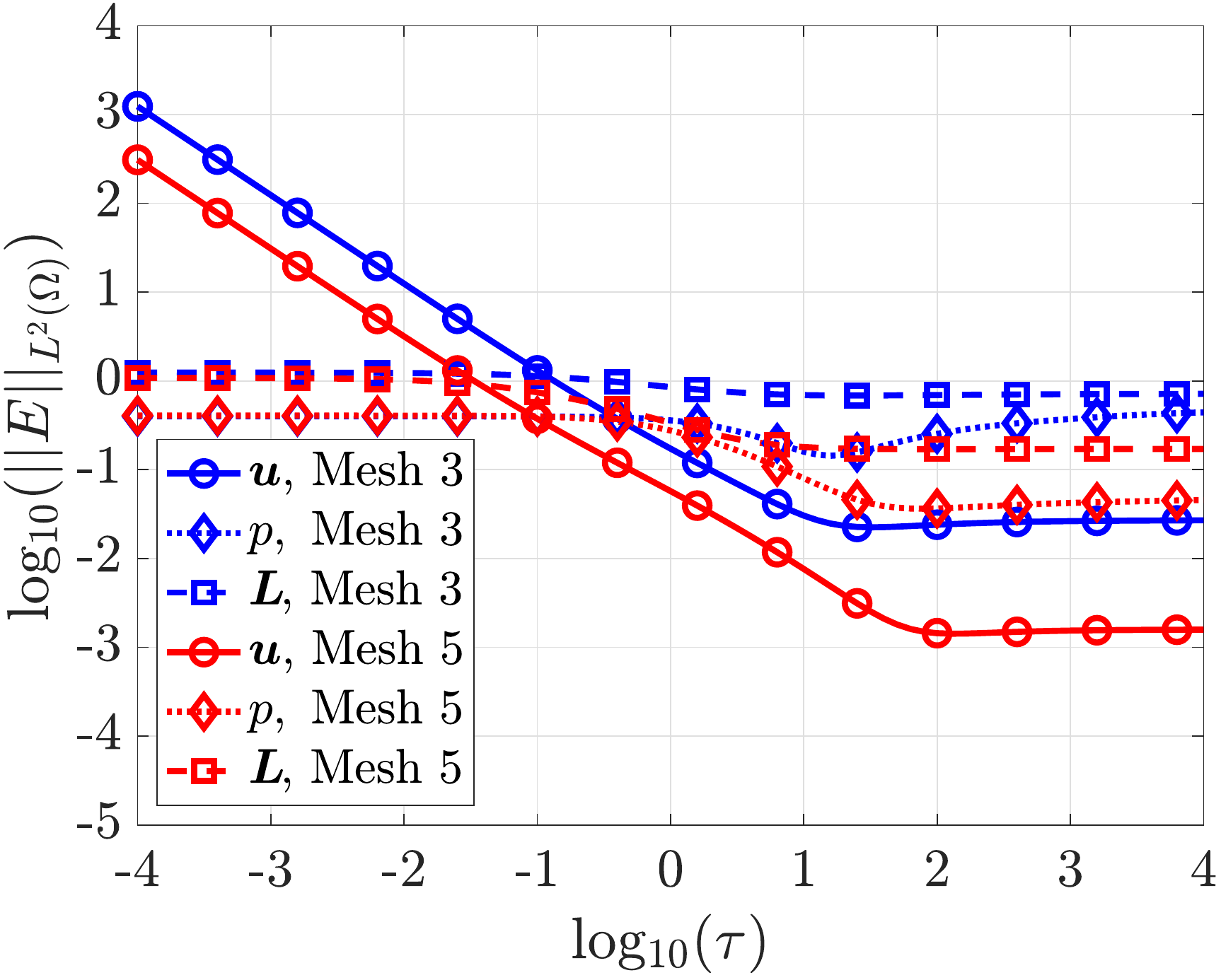}} 
	\subfigure[Prismatic cells]{\includegraphics[width=0.32\textwidth]{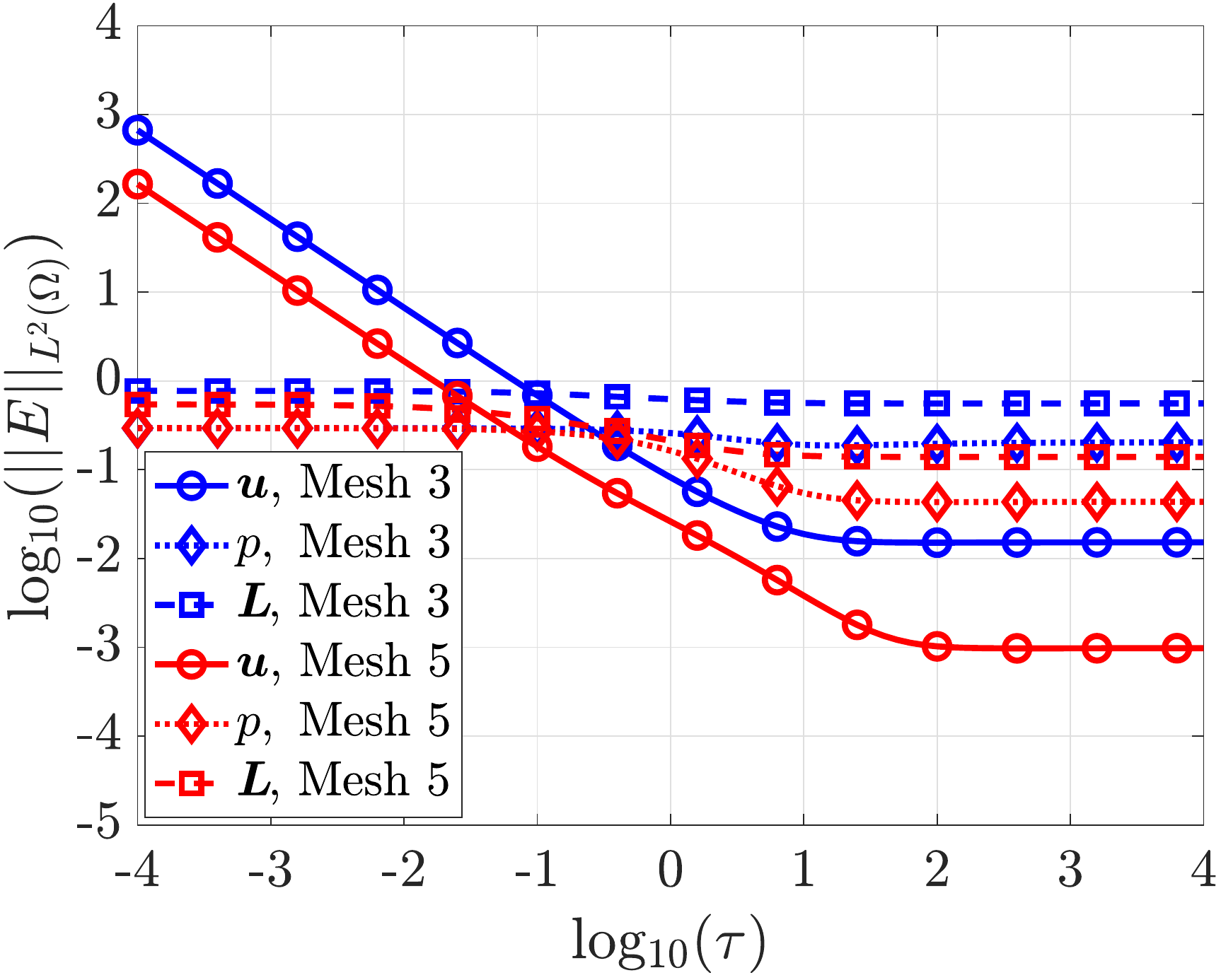}}
	\subfigure[Pyramidal cells]{\includegraphics[width=0.32\textwidth]{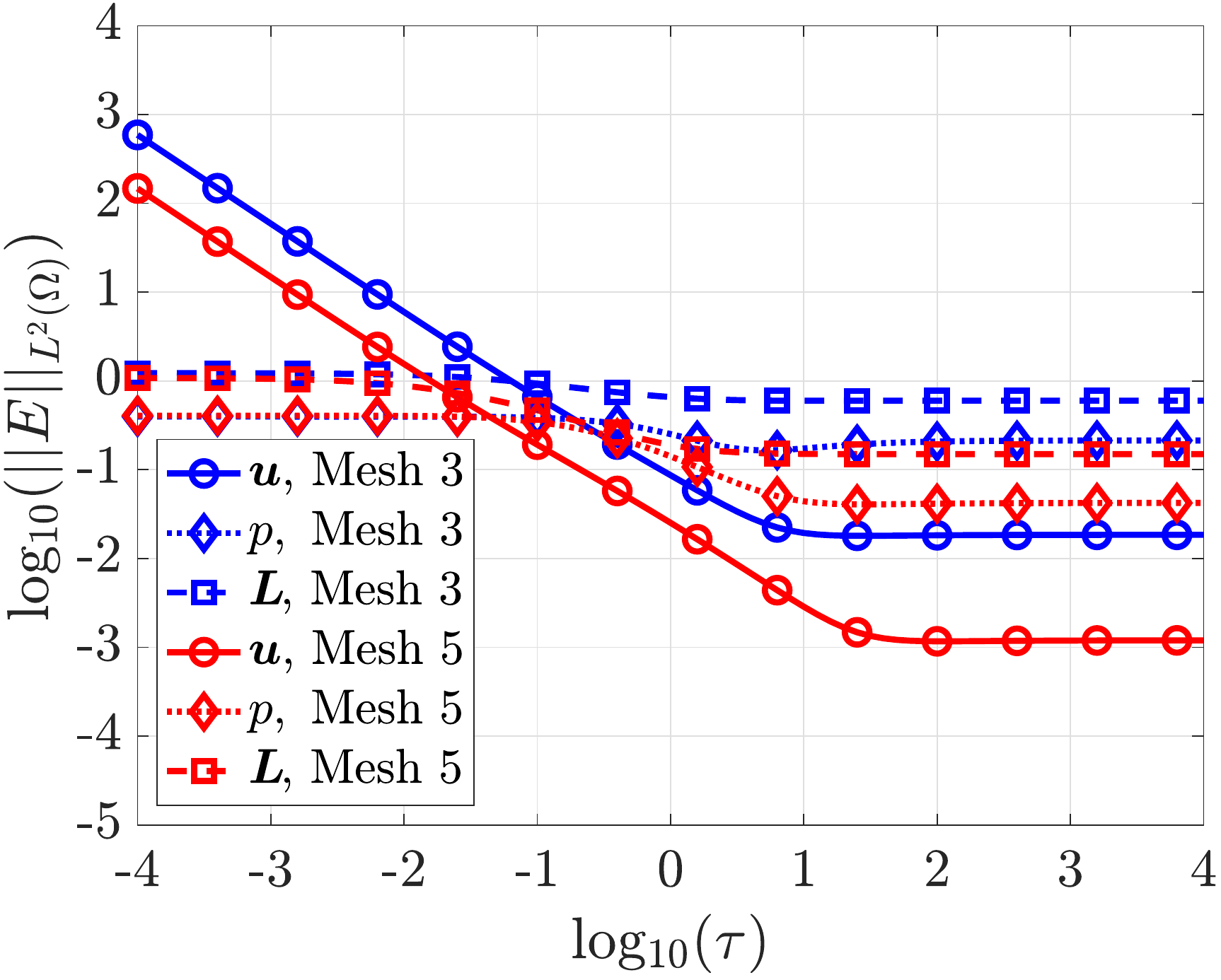}}	
	\caption{Error of velocity $\bu$, pressure $p$ and gradient of velocity $\bL$ in the $\eltwo(\Omega)$ norm as a function of the stabilisation parameter $\tau$ for the Stokes problem using (a) triangular, (b) quadrilateral, (c) tetrahedral, (d) hexahedral, (e) prismatic and (f) pyramidal cells.}
	\label{fig:Stokes_tau}
\end{figure}

Figure~\ref{fig:Stokes_tau} shows the evolution of the relative error of velocity, pressure and gradient of velocity  in the $\eltwo(\Omega)$ norm as a function of the stabilisation parameter $\tau$. 
The results include two different levels of mesh refinement and different cell types in two and three dimensions. It is worth emphasising that the range of values utilised for the experiment in two dimensions is different to the range used in three dimensions. 

\begin{remark}
For the Stokes equations, the usual definition of the stabilisation parameter is $\tau = \kappa \nu / \ell$, where $\nu$ is the viscosity of the fluid, $\ell$ is a characteristic length of the domain and $\kappa$ is a constant scaling factor~\cite{HDGtutorialNS}. For the case under analysis, it holds $\nu = 1$ and $\ell = 1$, the domain being the unit square and the unit cube in two and three dimensions respectively. Hence, figure~\ref{fig:Stokes_tau} is obtained by varying the value of the scaling factor $\kappa$ in the definition of the stabilisation parameter above.
\end{remark}

In all cases, the results show that a low value of the stabilisation parameter leads to a high error for the primal variable, whereas a large value, namely $\tau=10^4$ in two dimensions and $\tau=10^2$ in three dimensions, provides the maximum accuracy for the velocity. When triangular or tetrahedral cells are considered, the error of the gradient of the velocity and the pressure is found to be independent on the value of the stabilisation parameter used. In contrast, for quadrilateral and hexahedral cells, the error of the gradient of the velocity and the pressure decreases as the value of $\tau$ increases. For these cell types, the maximum accuracy is reached for a value of $\tau=10^4$ in two dimensions and $\tau=10^2$ in three dimensions. Finally, for prismatic and pyramidal cells, the accuracy of the gradient of the velocity and the pressure is less sensitive to the choice of $\tau$ and the qualitative behaviour is extremely similar to the one observed for quadrilateral and hexahedral cells.

Henceforth, the stabilisation parameter is selected as $\tau=10^4$ in two dimensions and $\tau=10^2$ in three dimensions, for both Poisson and Stokes problems and for any cell type.

\subsection{Influence of cell distortion and stretching}
\label{sc:InfluenceDistortionStretching}

Previous experiments, employed to test the optimal approximation properties of the proposed method, involved regular meshes. In this section, the effect of cell distortion and stretching on the accuracy of the proposed method is studied. This is of major importance for the method to be applicable to more complicated problems involving complex geometries and for its extension to computational fluid dynamics applications involving boundary layers.

\begin{figure}[!tb]
	\centering
	\subfigure[Distored hexahedrons] {\includegraphics[width=0.24\textwidth]{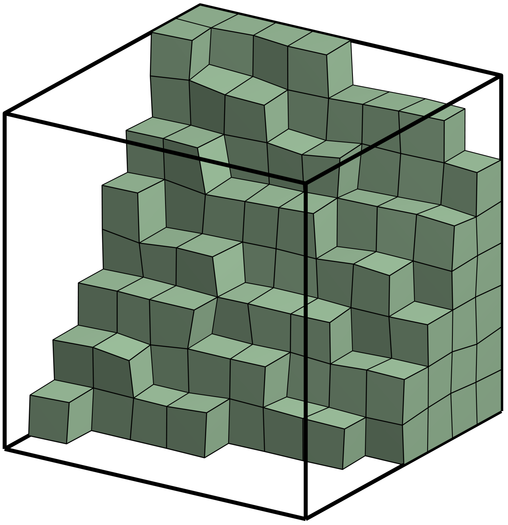}}
	\subfigure[Distored prisms]{\includegraphics[width=0.24\textwidth]{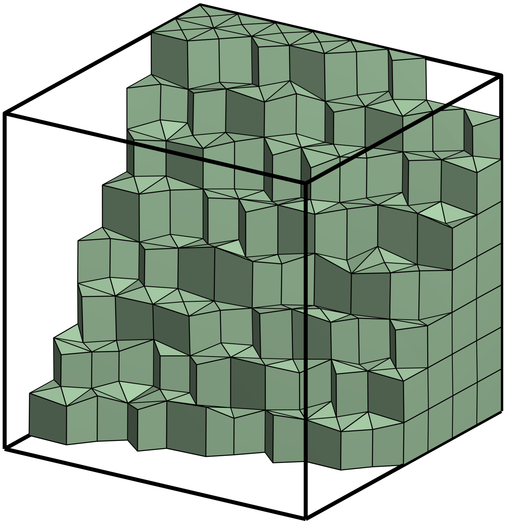}}
	\subfigure[Stretched tetrahedrons]{\includegraphics[width=0.24\textwidth]{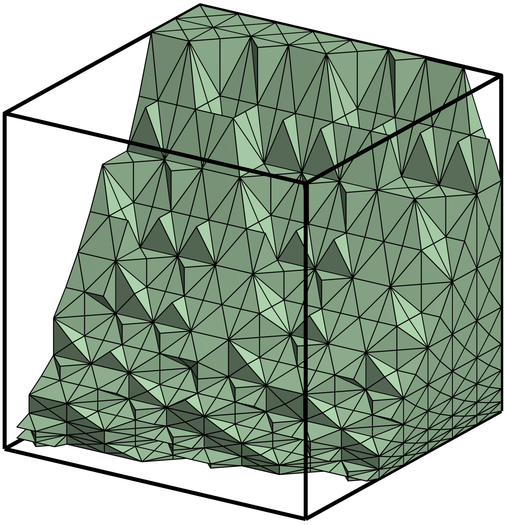}}	
	\subfigure[Stretched pyramids]  {\includegraphics[width=0.24\textwidth]{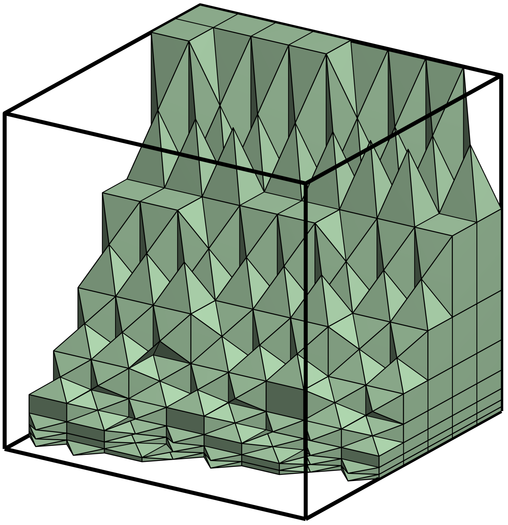}}
	\caption{Internal view of the meshes corresponding to the fourth level of refinement for the domain $\Omega=[0,1]^3$ featuring (a) distorted hexahedral, (b) distorted prismatic, (c) stretched tetrahedral and (d) stretched pyramidal cells.}
	\label{fig:3DmeshesDistortedStretched}
\end{figure}

Cell distortion is introduced by perturbing the internal nodes of the mesh according to a random variation of maximum magnitude $h_{\text{min}}$/4, where $h_{\text{min}}$ is the minimum edge of the regular mesh. For cells with quadrilateral faces, the motion is constrained to ensure that all faces on the distorted mesh are planar~\cite{fovrt2011finite}. Figure~\ref{fig:3DmeshesDistortedStretched} (a) and (b) show two examples of distorted meshes for hexahedral and prismatic cells.
Similarly, cell stretching is introduced by transforming the regular meshes employed in previous experiments. The stretching factor, $s$, is measured as the ratio between the maximum and minimum faces/edges in a cell. Figure~\ref{fig:3DmeshesDistortedStretched} (c) and (d) show two examples of stretched meshes for tetrahedral and pyramidal cells for $s=10$.

\begin{figure}[!tb]
	\centering
	\subfigure[$\bu$]{\includegraphics[width=0.32\textwidth]{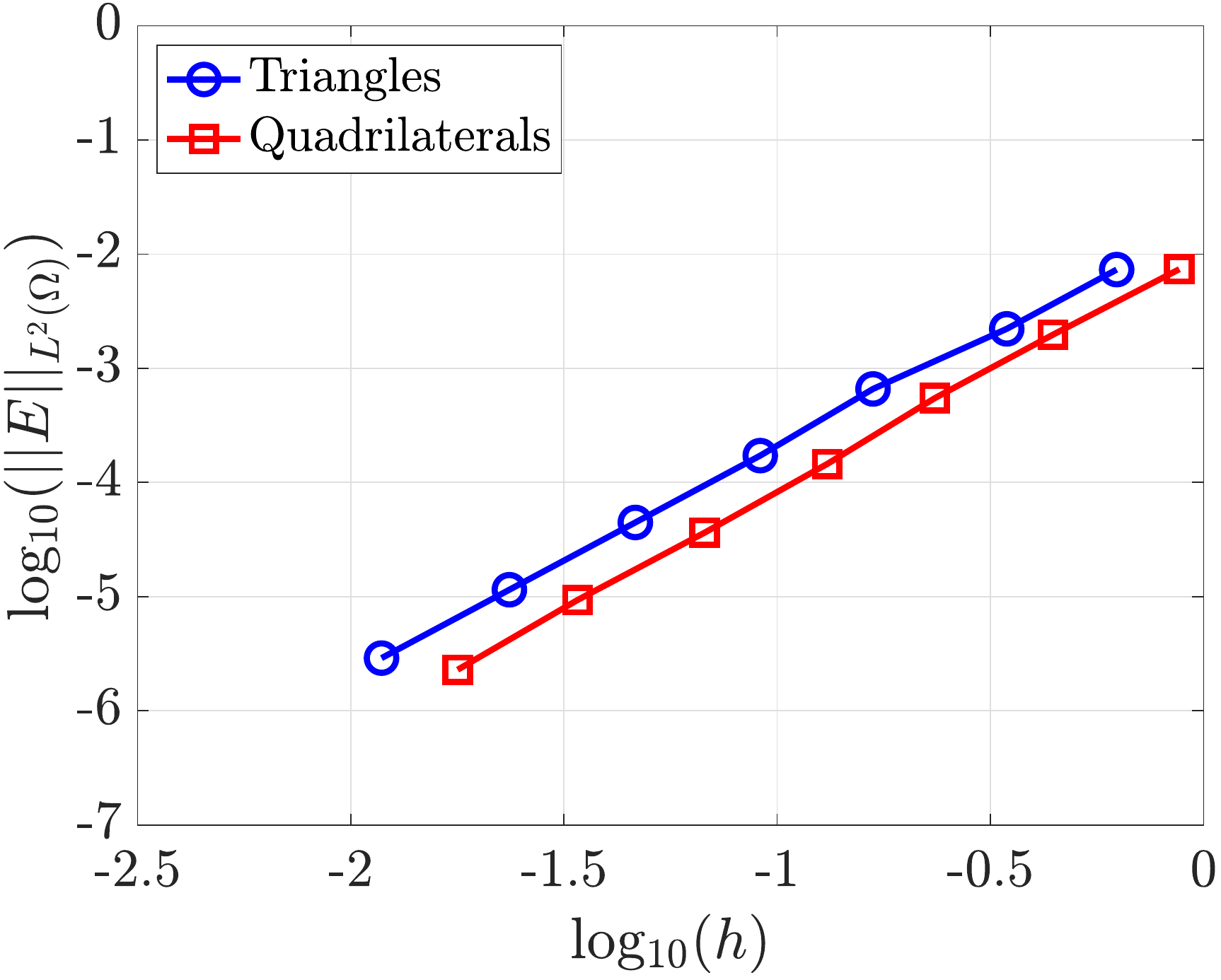}}
	\subfigure[$p$]{\includegraphics[width=0.32\textwidth]{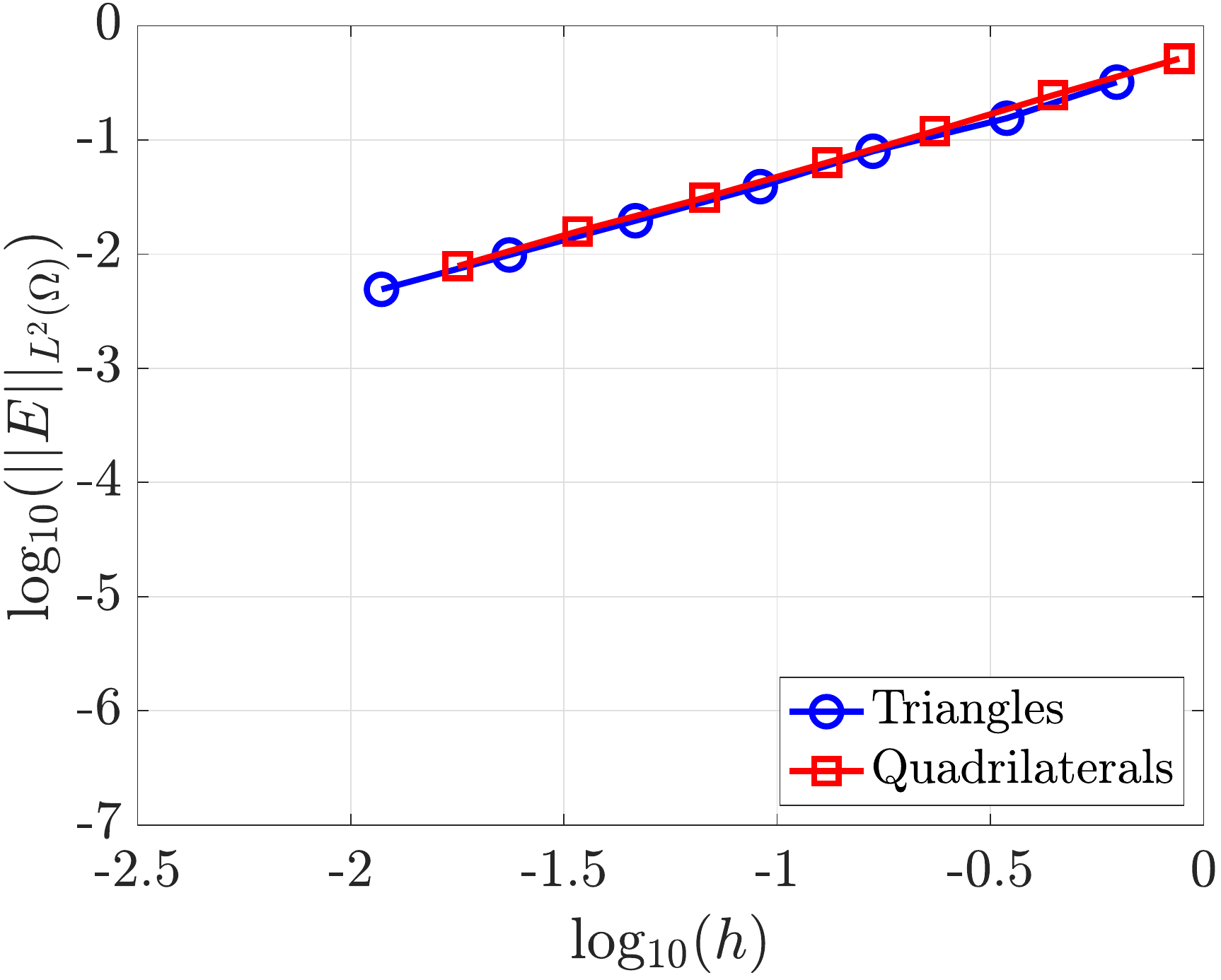}}
	\subfigure[$\bL$]{\includegraphics[width=0.32\textwidth]{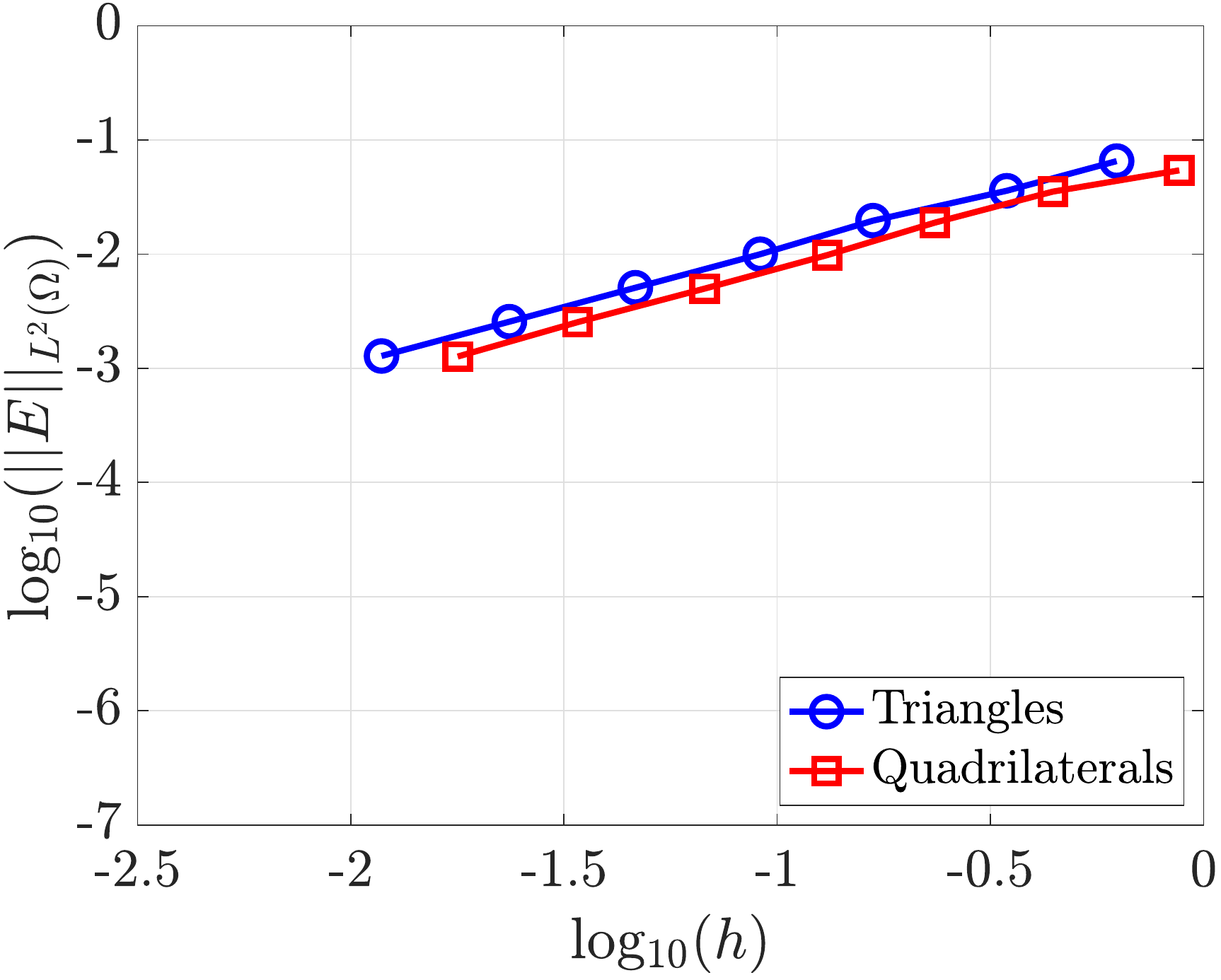}}
	\caption{Mesh convergence of the error of velocity $\bu$, pressure $p$ and gradient of velocity $\bL$ in the $\eltwo(\Omega)$ norm as a function of the mesh size $h$ for two dimensional Stokes problem using meshes of distorted cells.}
	\label{fig:Stokes_Conv_2DPerturb}
\end{figure}

\begin{figure}[!tb]
	\centering
	\subfigure[$\bu$]{\includegraphics[width=0.32\textwidth]{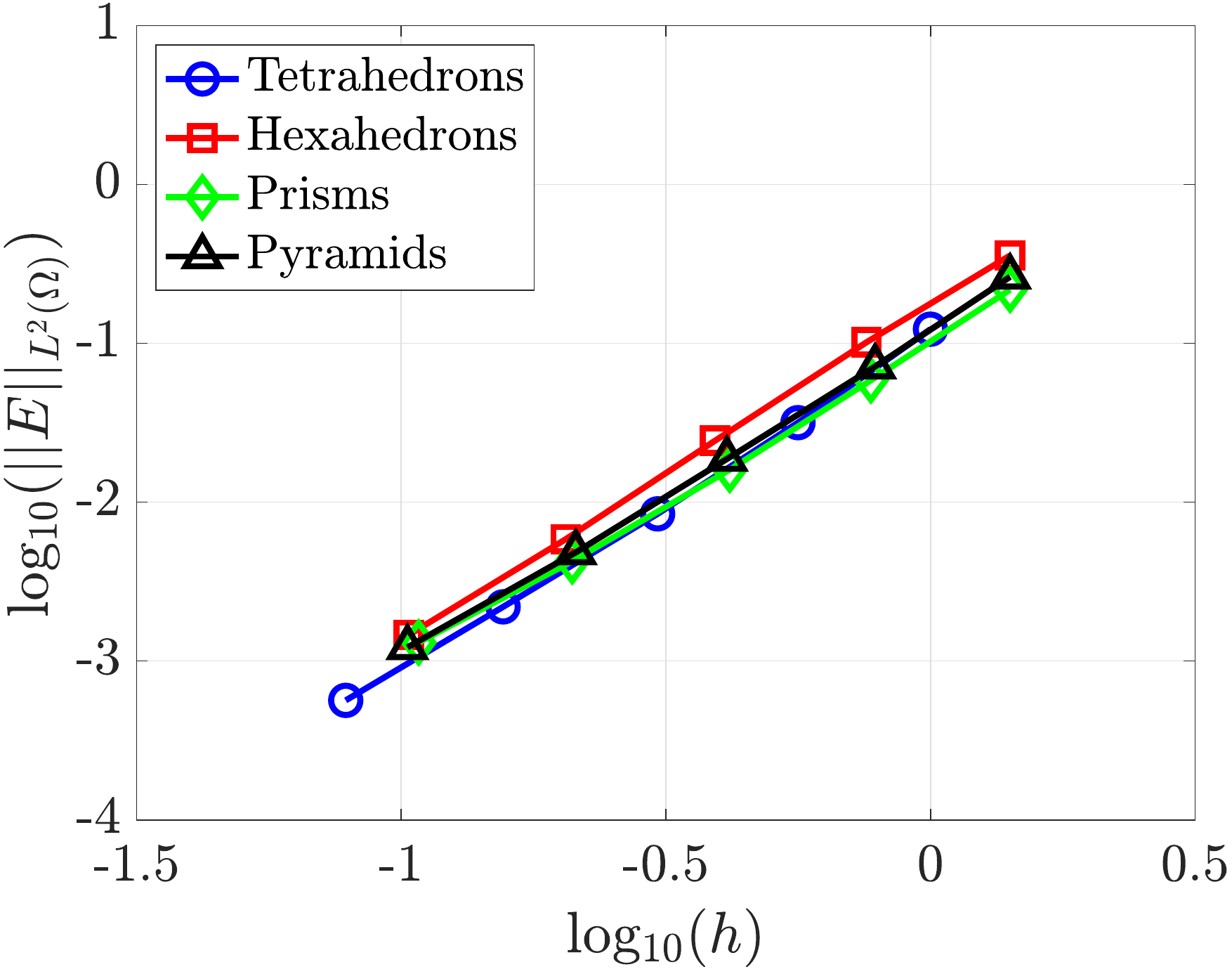}}
	\subfigure[$p$]{\includegraphics[width=0.32\textwidth]{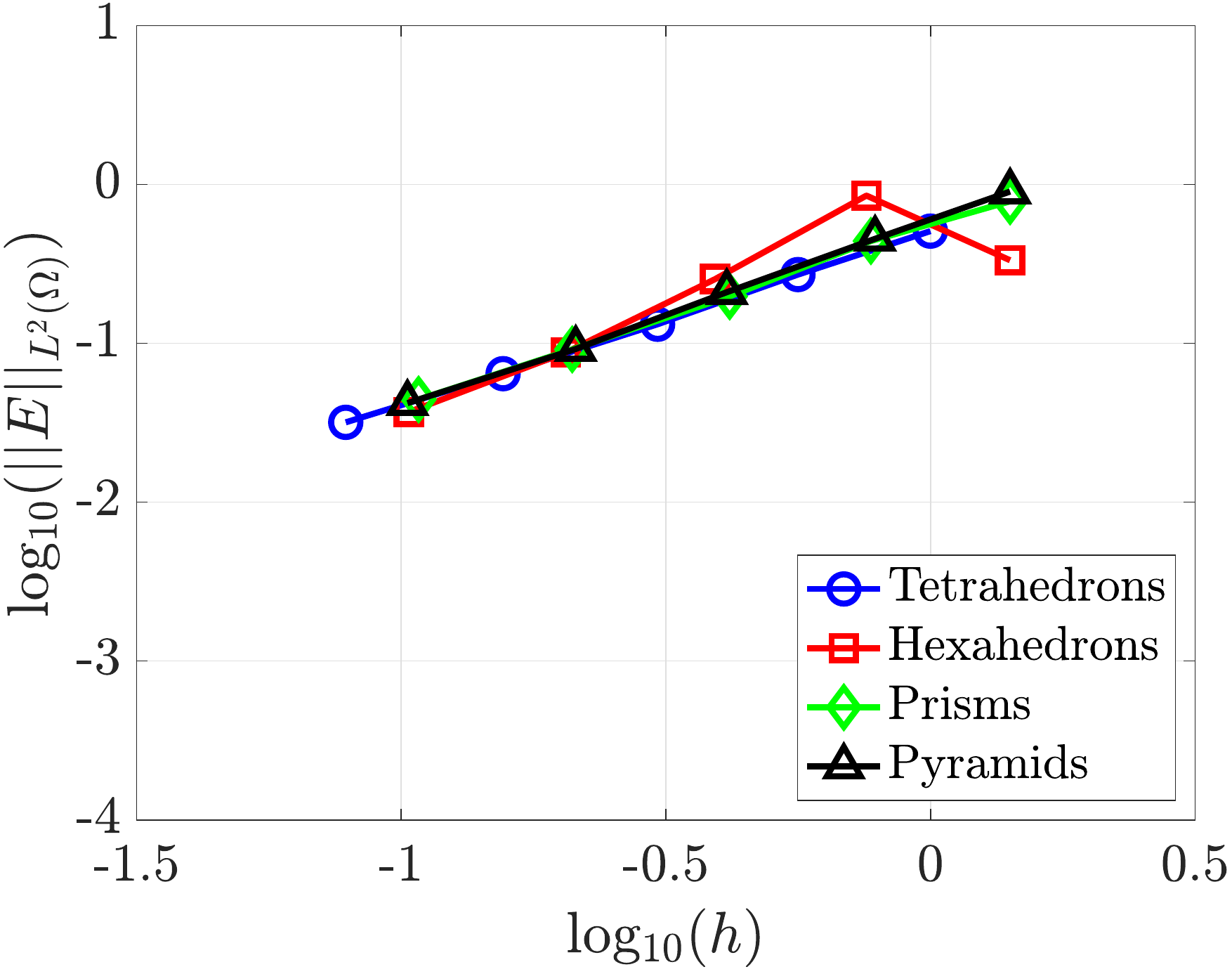}}
	\subfigure[$\bL$]{\includegraphics[width=0.32\textwidth]{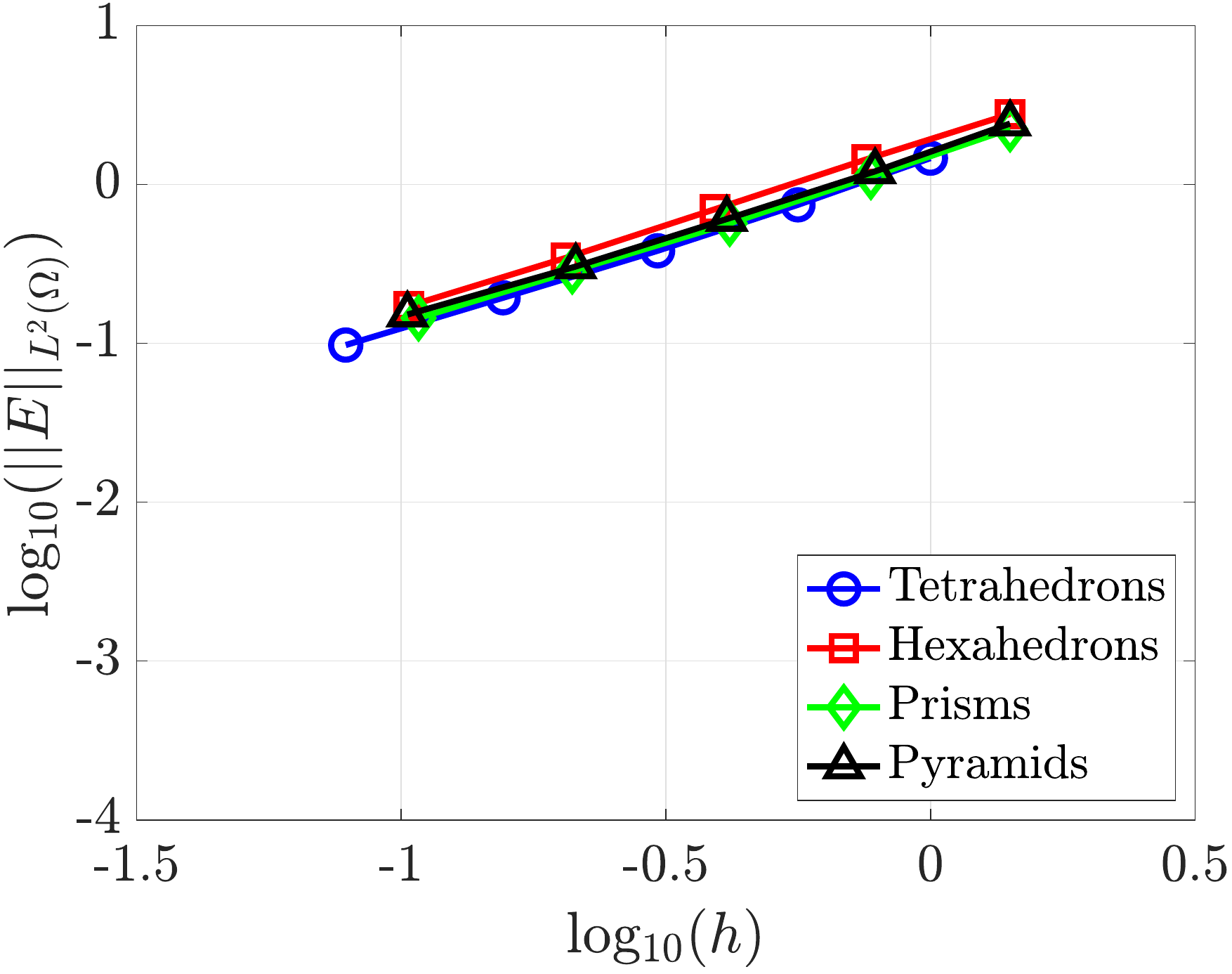}}
	\caption{Mesh convergence of the error of velocity $\bu$, pressure $p$ and gradient of velocity $\bL$ in the $\eltwo(\Omega)$ norm as a function of the mesh size $h$ for three dimensional Stokes problem using meshes of distorted cells.}
	\label{fig:Stokes_Conv_3DPerturbCoplanar}
\end{figure}

The mesh convergence results for the Stokes problem in two and three dimensions using meshes with distorted cells are shown in Figures~\ref{fig:Stokes_Conv_2DPerturb} and \ref{fig:Stokes_Conv_3DPerturbCoplanar}, respectively.
In two dimensions, the optimal convergence properties are observed for velocity, pressure and gradient of velocity both using triangular and quadrilateral meshes. Furthermore, it can be observed that, for the same level of mesh refinement, quadrilateral cells provide more accurate results when compared to meshes with triangular cells. Similar conclusions are obtained in three dimensions, where optimal rate of convergence is achieved in all cases for velocity, pressure and gradient of velocity and comparable accuracy is provided by all cell types. 

\begin{figure}[!tb]
	\centering
	\subfigure[$u$]{\includegraphics[width=0.32\textwidth]{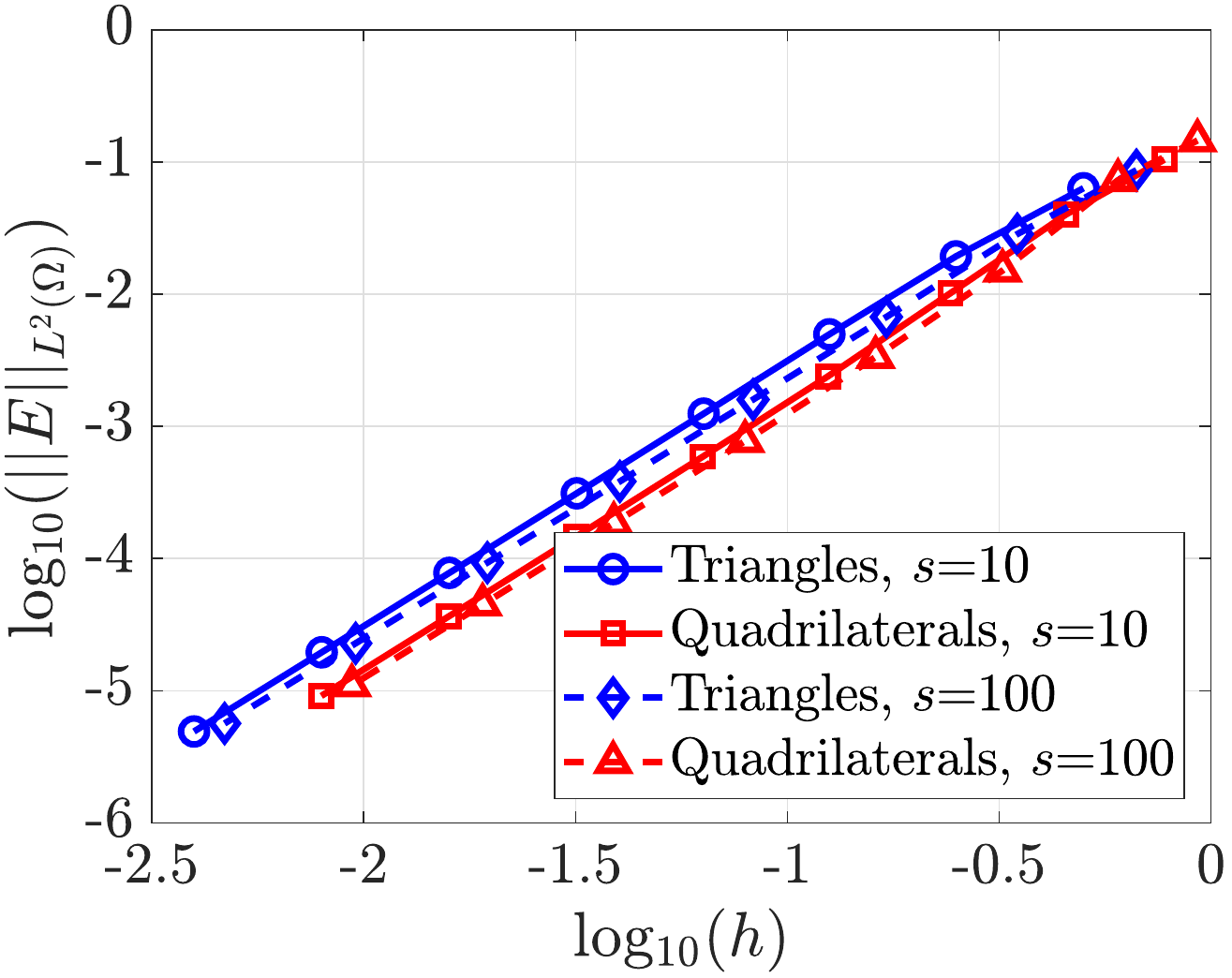}}
	\subfigure[$\bq$]{\includegraphics[width=0.32\textwidth]{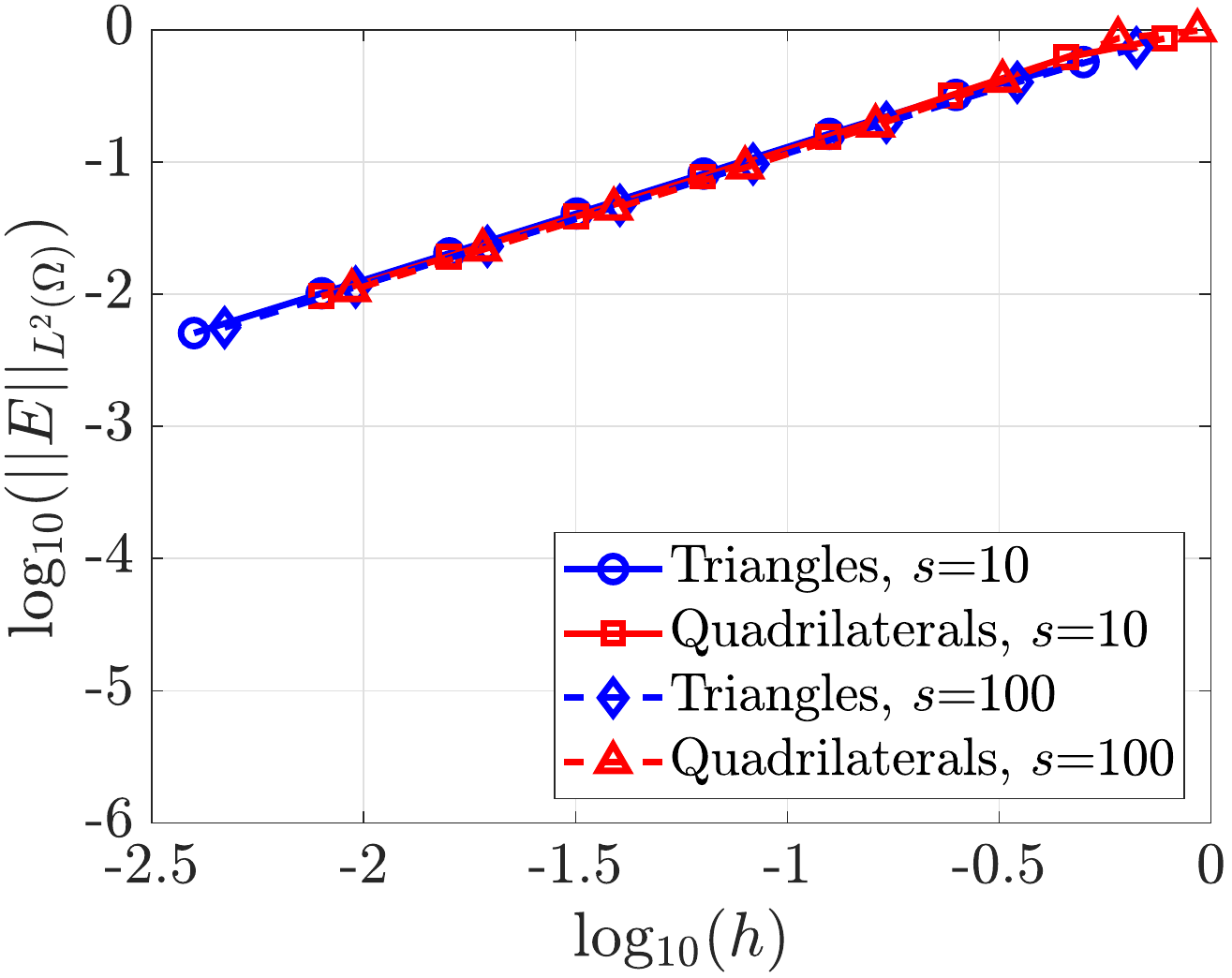}}
	\caption{Mesh convergence of the error of the solution $u$ and its gradient $\bq$ in the $\eltwo(\Omega)$ norm for the two dimensional Poisson problem using meshes of stretched cells, with maximum stretching factor $s{=}10$ and $s{=}100$.}
	\label{fig:Poisson_Conv_2DStretch}
\end{figure}

\begin{figure}[!tb]
	\centering
	\subfigure[$u$, $s=10$]{\includegraphics[width=0.24\textwidth]{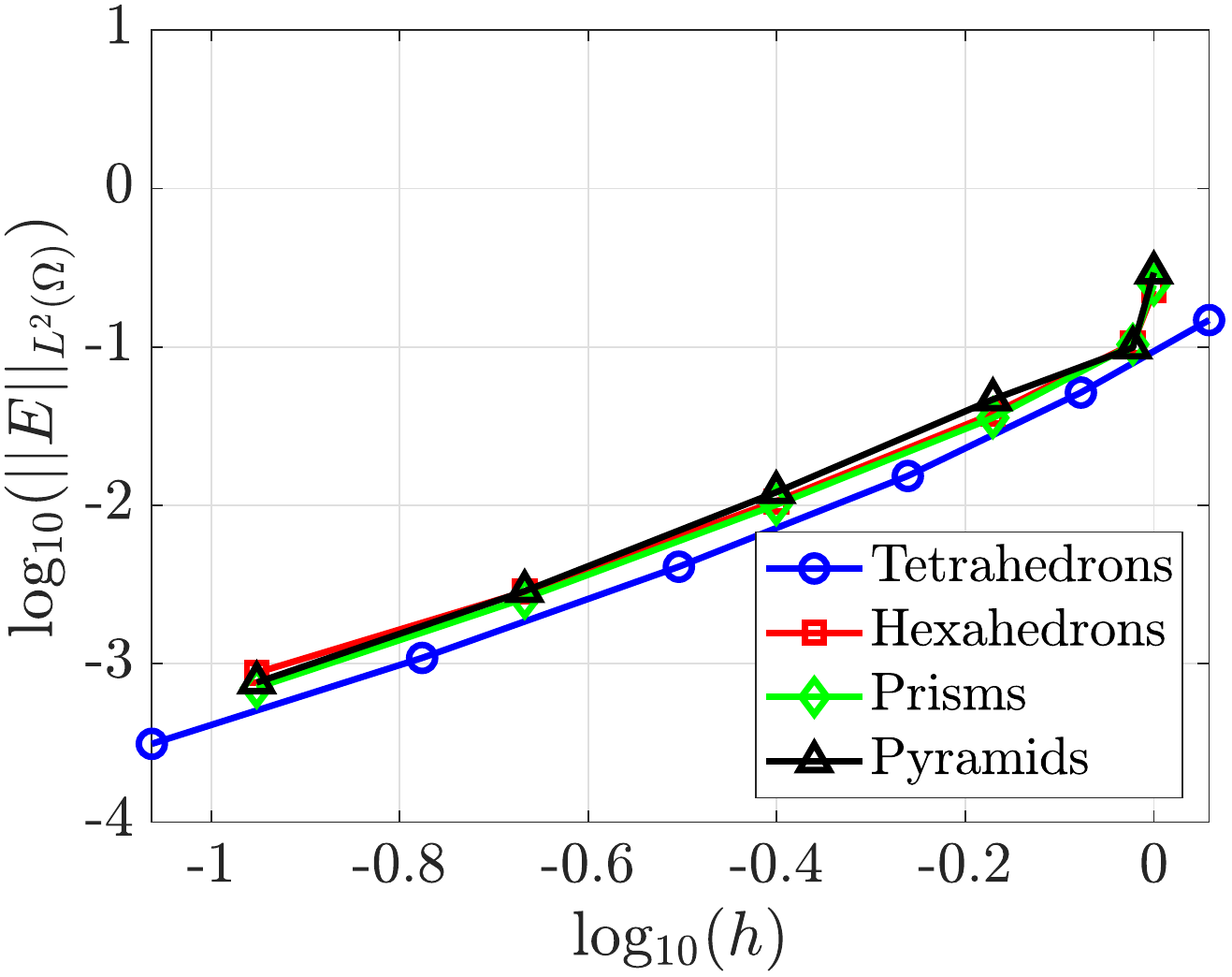}}
	\subfigure[$\bq$, $s=10$]{\includegraphics[width=0.24\textwidth]{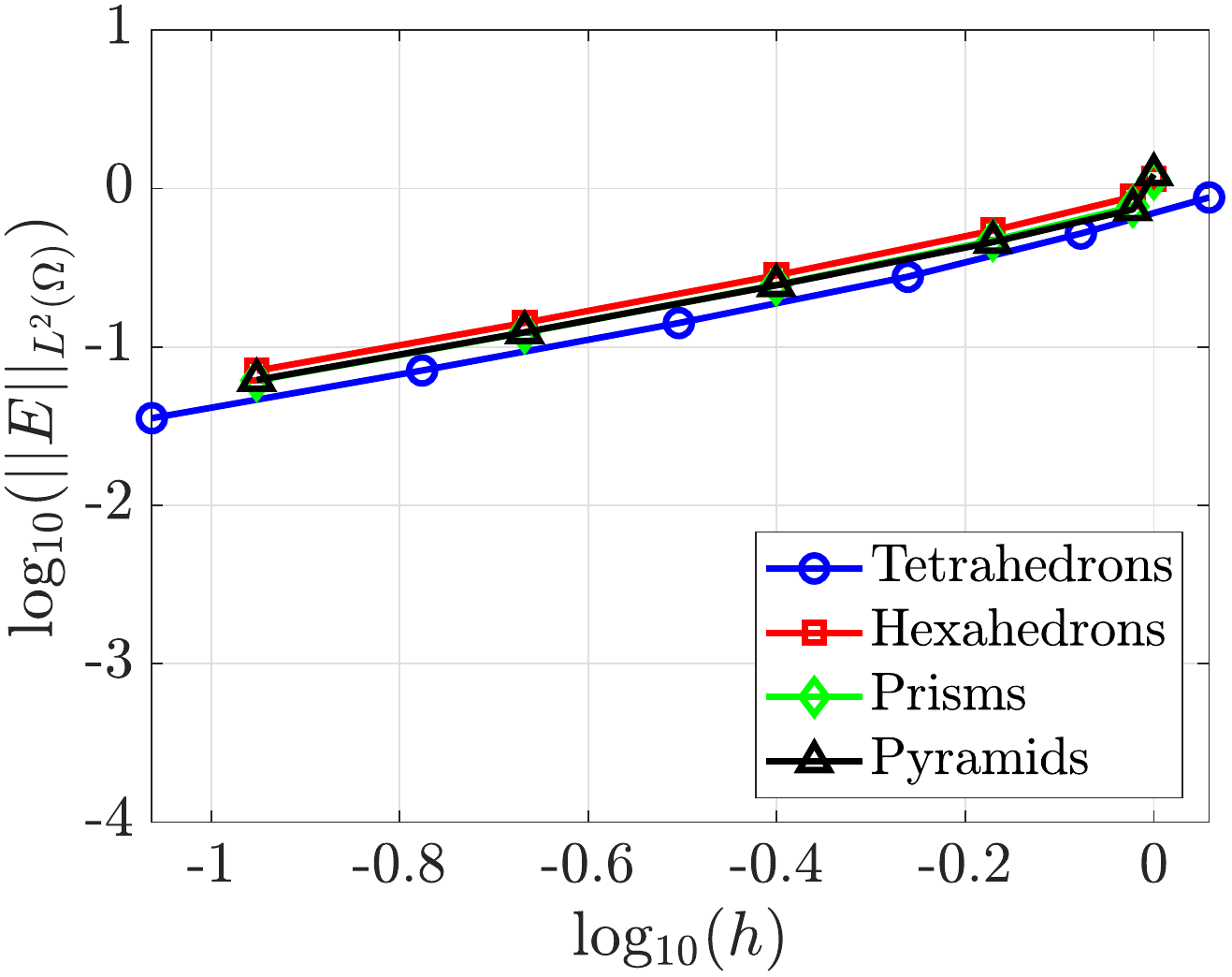}} 
	\subfigure[$u$, $s=100$]{\includegraphics[width=0.24\textwidth]{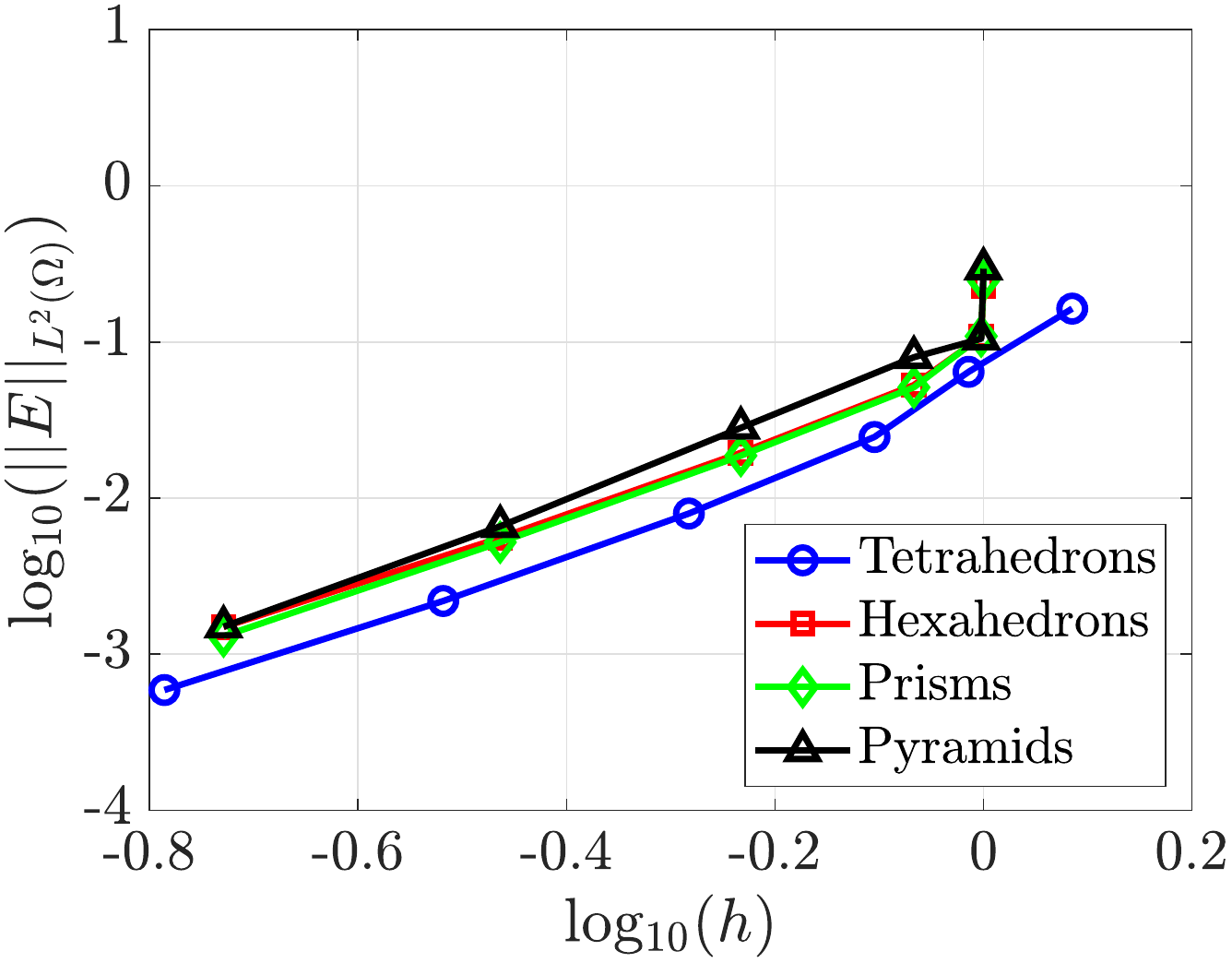}}
	\subfigure[$\bq$, $s=100$]{\includegraphics[width=0.24\textwidth]{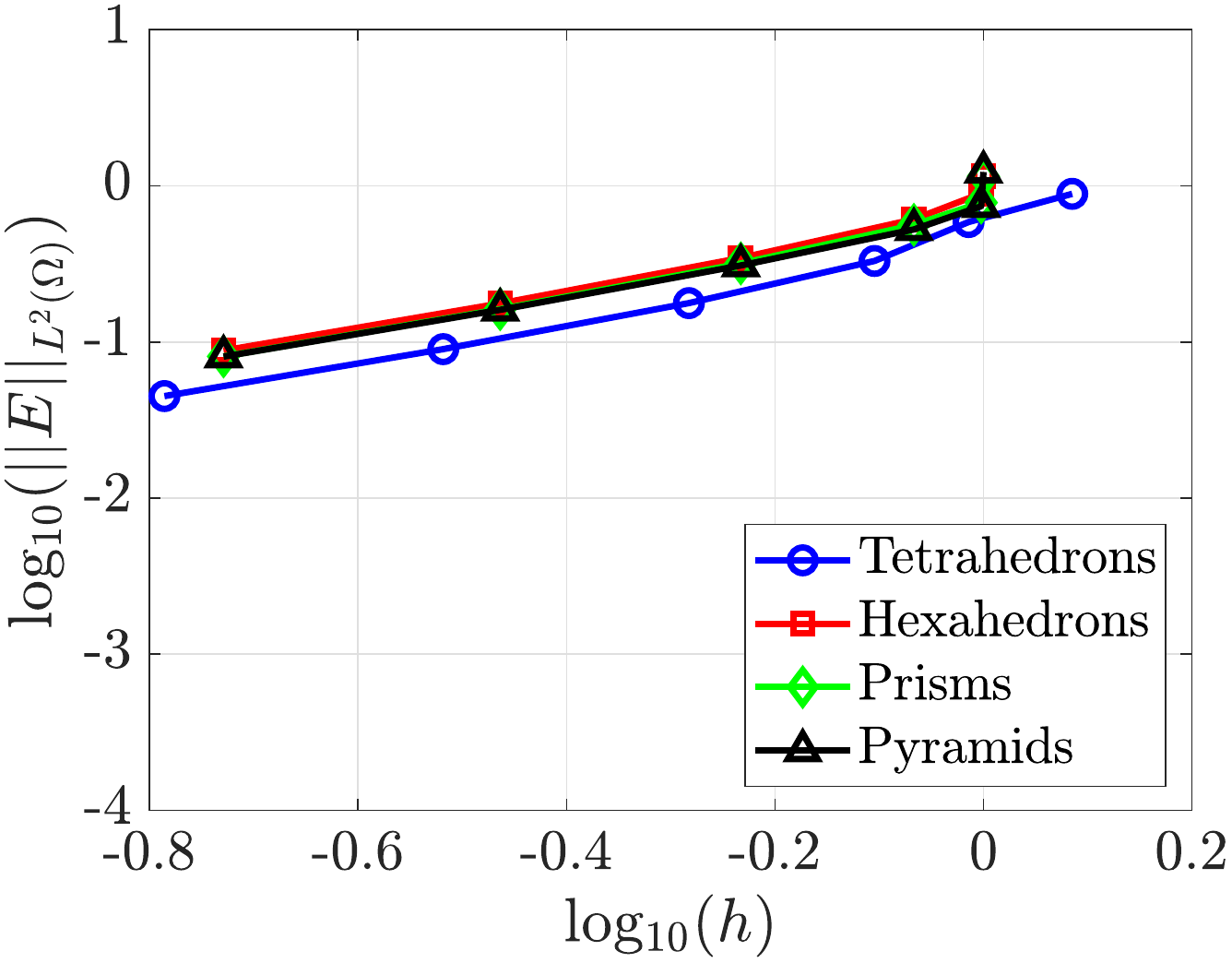}}
	\caption{Mesh convergence of the error of the solution $u$ and its gradient $\bq$ in the $\eltwo(\Omega)$ norm for the three dimensional Poisson problem using meshes of stretched cells, with maximum stretching factor $s=10$ and $s=100$.}
	\label{fig:Poisson_Conv_3DStretch}
\end{figure}

The convergence study on stretched meshes with stretching factor $s=10$ and $s=100$, is performed for the Poisson problem. Figure~\ref{fig:Poisson_Conv_2DStretch} shows the relative error, measured in the $\eltwo(\Omega)$ norm, of the primal and mixed variables as a function of the characteristic cell size. 
The results reveal that the accuracy of the proposed method is not dependent upon the stretching factor. The conclusions also hold for the Poisson problem in three dimensions, as illustrated by the results in figure~\ref{fig:Poisson_Conv_3DStretch}.
Optimal convergence is observed for all the variables and all cell types. The accuracy is again found to be almost insensitive to the stretching factor. 

Further numerical experiments, not reported here for brevity, demonstrated that the same conclusions are obtained when performing the numerical experiments for the Poisson problem on meshes with distorted cells and the Stokes problem on meshes with stretched cells.

\subsection{Computational cost}
\label{sc:ComputationalCost}

The last numerical experiment involves a study of the computational cost of the proposed method for meshes with different cell types. The computational efficiency is compared by directly measuring the CPU time (in seconds) required to assemble and solve the global system of equations, as this is the dominant cost of the proposed methodology. 

\begin{figure}[!tb]
	\centering
	\subfigure[$\bu$]{\includegraphics[width=0.32\textwidth]{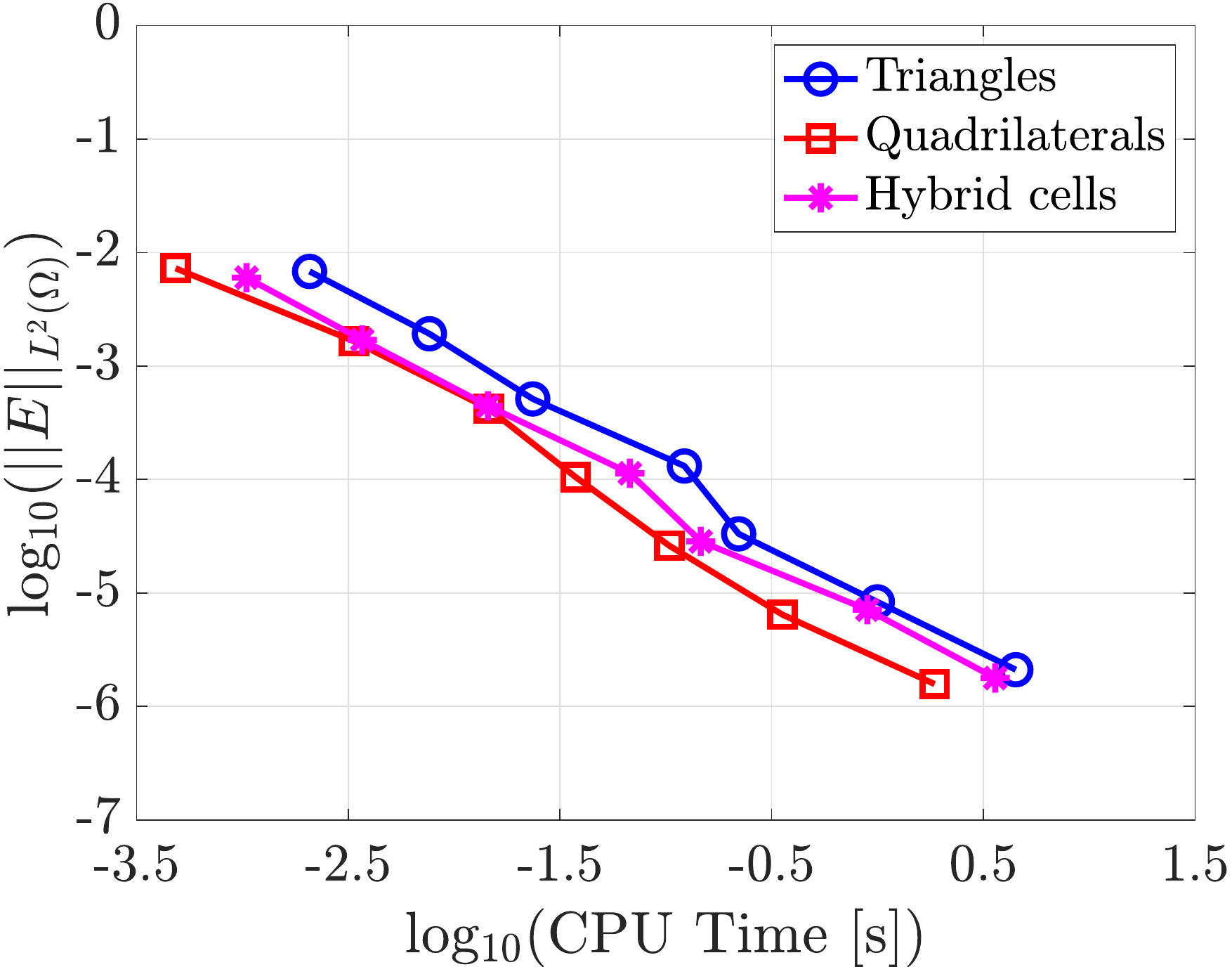}}
	\subfigure[$p$]{\includegraphics[width=0.32\textwidth]{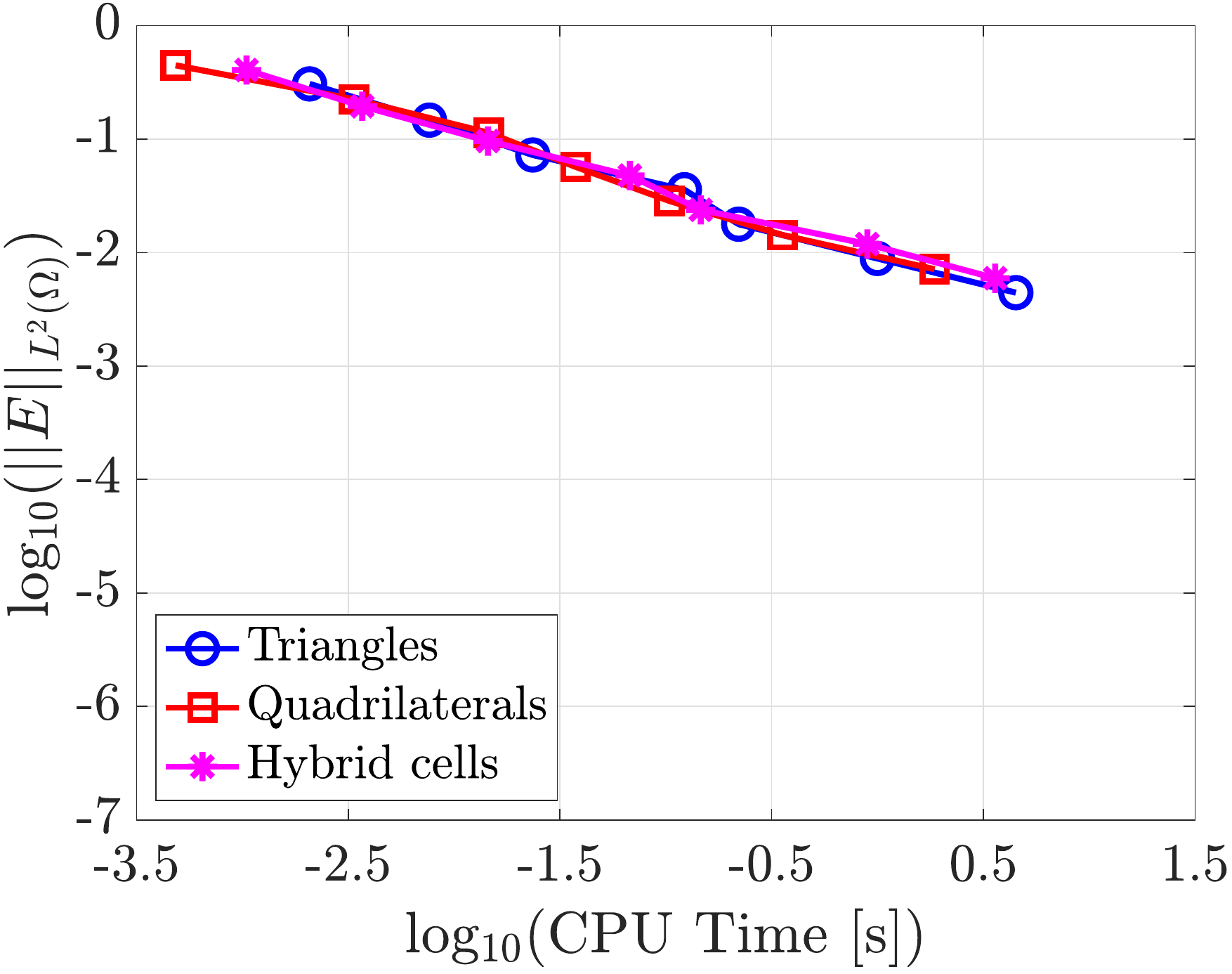}}
	\subfigure[$\bL$]{\includegraphics[width=0.32\textwidth]{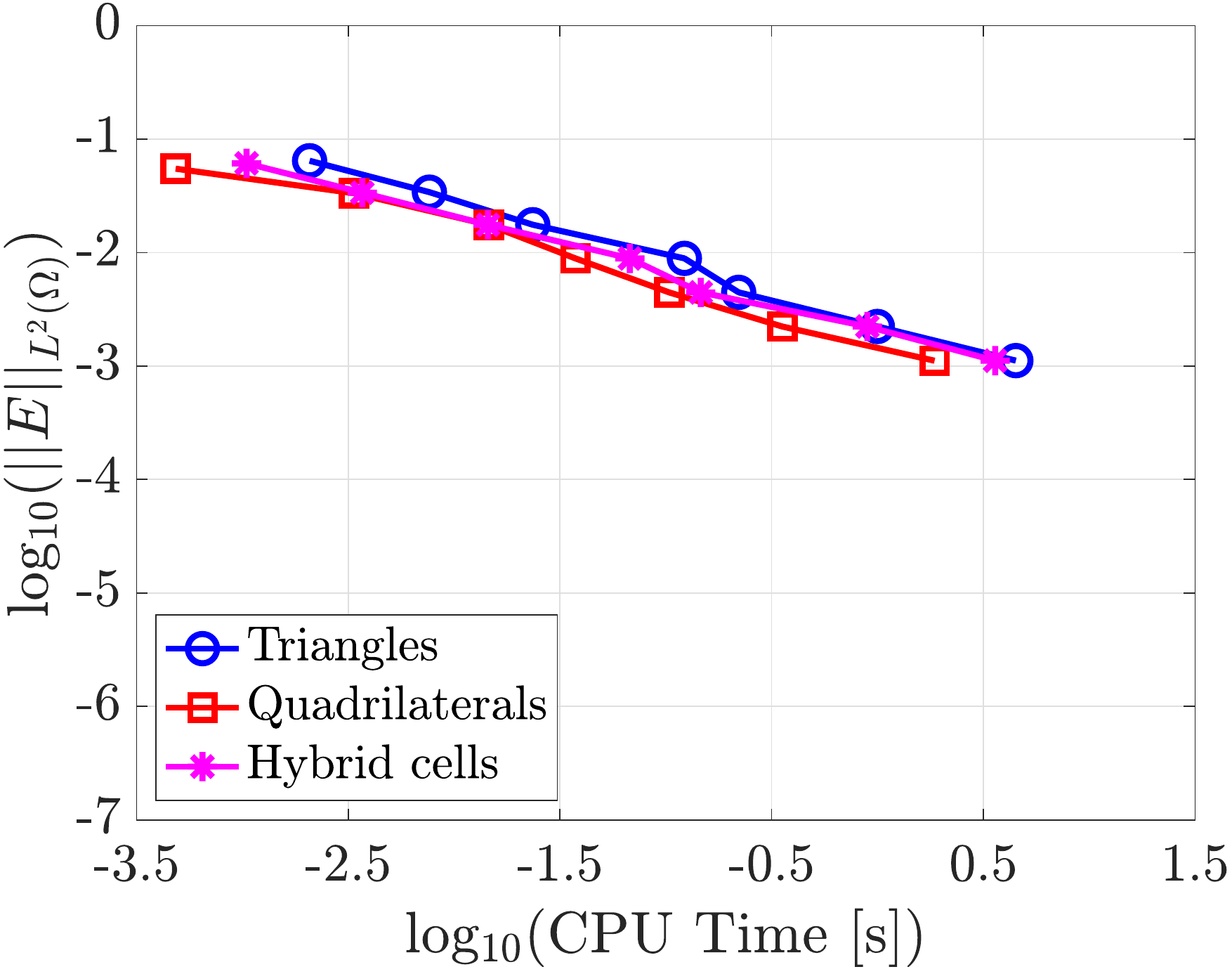}}
	\caption{Error of velocity $\bu$, pressure $p$ and gradient of velocity $\bL$ in the $\eltwo(\Omega)$ norm as a function of the CPU time for two dimensional Stokes problem using meshes of different cell types.}
	\label{fig:Stokes_CPU2D}
\end{figure}

Figure~\ref{fig:Stokes_CPU2D} displays the evolution of the relative error of velocity, pressure and gradient of velocity in the $\eltwo(\Omega)$ norm, as a function of the CPU time. 
The results reveal that quadrilateral and hybrid meshes provide the same accuracy as triangular meshes with slightly less computational effort. The better performance of quadrilateral cells is clearly observed when measuring the error of the velocity, whereas for pressure and gradient of velocity, all types of cells provide the same accuracy with a similar computational effort. It is worth noting that the advantages of using quadrilateral cells are not only observed when high accuracy is required. Even for an accuracy of 1\% quadrilateral cells require nearly one order of magnitude less CPU time than triangles.

\begin{figure}[!tb]
	\centering
	\subfigure[$\bu$]{\includegraphics[width=0.32\textwidth]{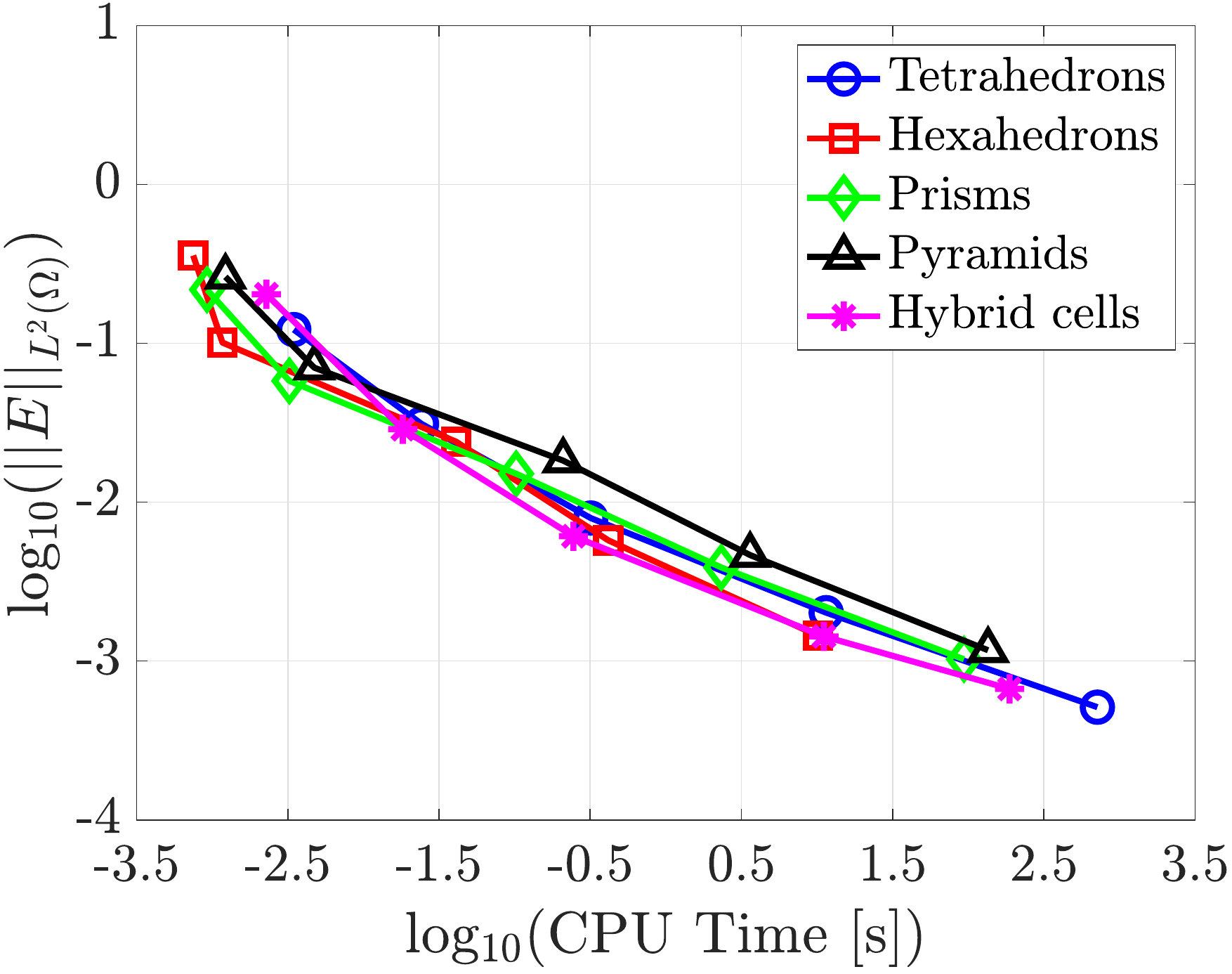}}
	\subfigure[$p$]{\includegraphics[width=0.32\textwidth]{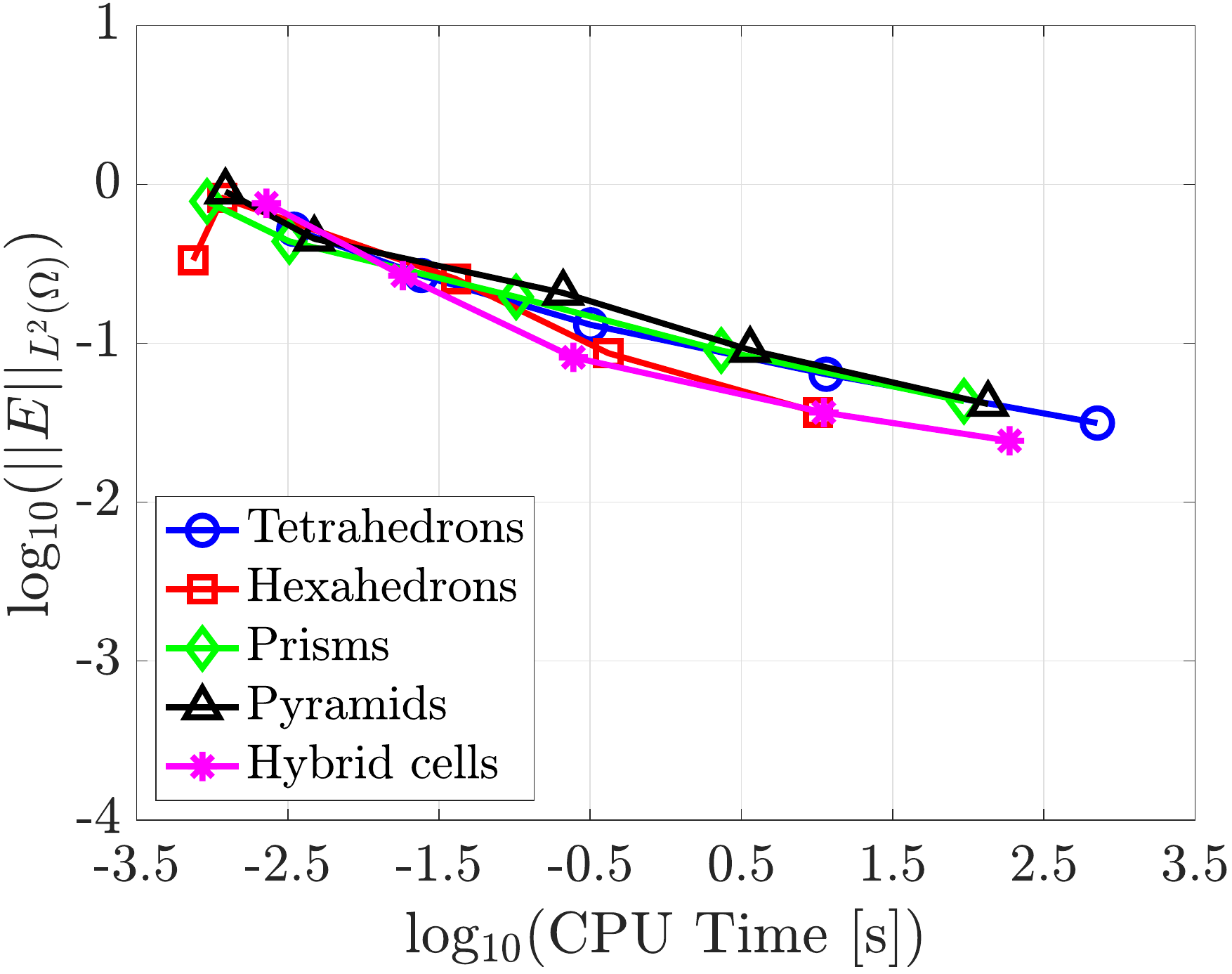}}
	\subfigure[$\bL$]{\includegraphics[width=0.32\textwidth]{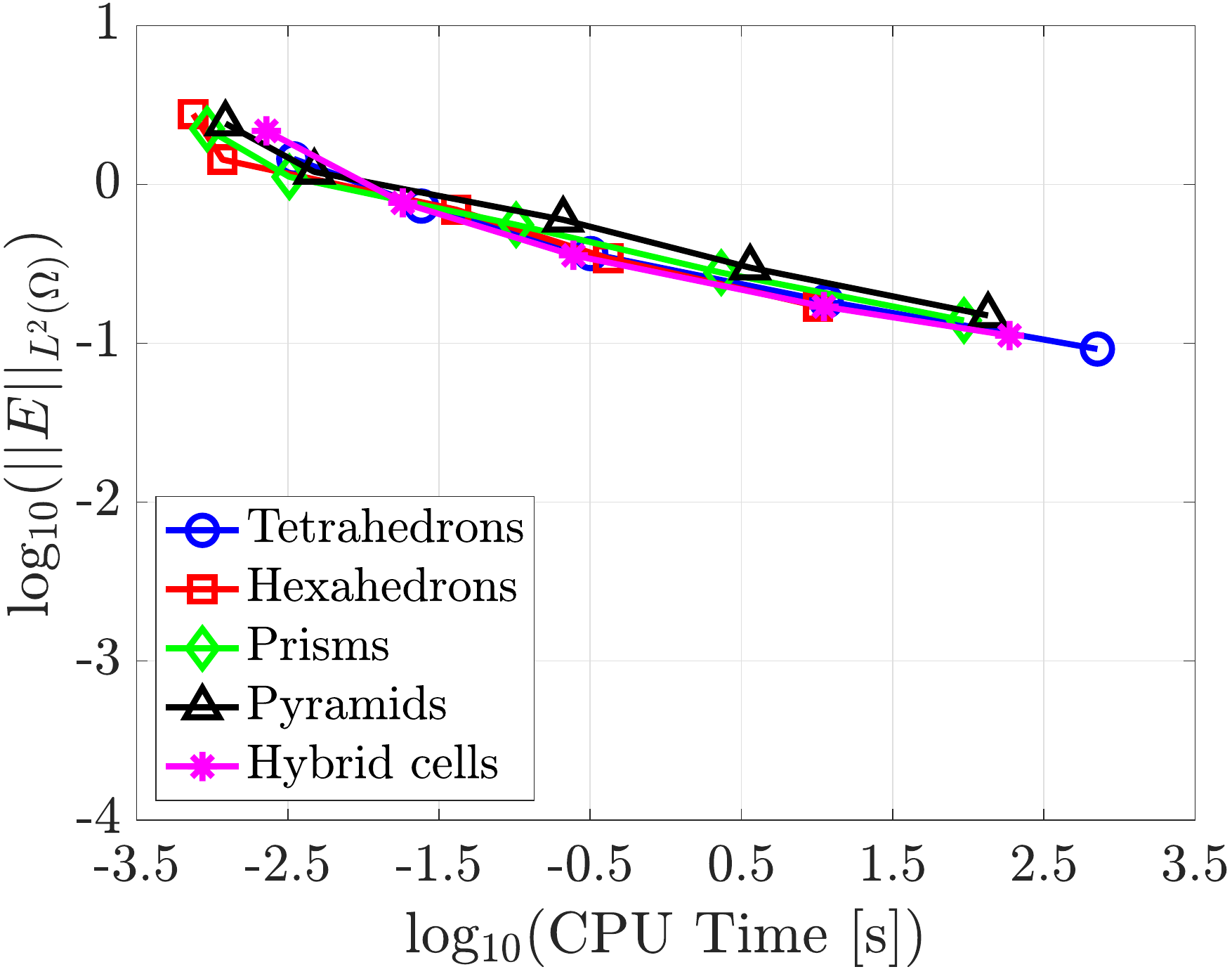}}
	\caption{Error of velocity $\bu$, pressure $p$ and gradient of velocity $\bL$ in the $\eltwo(\Omega)$ norm as a function of the CPU time for three dimensional Stokes problem using meshes of different cell types.}
	\label{fig:Stokes_CPU3D}
\end{figure}

In three dimensions the conclusions are similar, as illustrated in Figure~\ref{fig:Stokes_CPU3D}.
Hexahedral and hybrid meshes are able to provide the solution with a given accuracy with slightly less computational effort when compared to tetrahedral, prismatic and pyramidal meshes. The most important differences are appreciated when the error of the velocity is considered.

It is worth noting that quadrilateral and hexahedral meshes seem to provide the maximum performance, in terms of achieving the desired accuracy with the minimum computational effort. However, it is well known that the mesh generation of complex objects using unstructured hexahedral meshes is still today an open problem~\cite{AIAA-meshChallenges}. In this scenario, the ability of the proposed method to handle hybrid meshes will be of use. As demonstrated by the results of Figures~\ref{fig:Stokes_CPU2D} and~\ref{fig:Stokes_CPU3D}, the use of hybrid meshes is still beneficial when compared to pure triangular or tetrahedral meshes.

\section{Applications of the automatic mesh adaptivity strategy}
\label{sc:examples}

This section presents two numerical examples solved with the proposed second-order FCFV in a mesh adaptivity framework. The first example involves the solution of a Poisson problem with known analytical solution on a simple two dimensional domain. This example is used to evaluate the performance of the error indicator introduced in section~\ref{sc:meshAdaptivity}. The second example involves the solution of the Stokes equations in a complex domain of interest for microfluidics applications~\cite{ZhangNano,keaveny2013optimization} and it is used to demonstrate the potential of the proposed methodology.

\subsection{Two dimensional heat transfer problem with localised source}
\label{sc:PoissonAdapt} 

The model problem~\eqref{eq:PoissonBrokenFirstOrder} in $\Omega=[0,1]^{2}$ is considered, where the source term and boundary data are selected such that the analytical solution is known and given by
\begin{equation}\label{eq:PoissonAdapt}
u^{\text{ex}}(x_1,x_2) = 1 + \exp \left\{ -a \left( (x_1-b)^2 + (x_2-b)^2 \right) \right\},
\end{equation}
with $a = 100$ and $b = 0.7$. The variation of the solution is confined to a small region in the domain, around the point (0.7,0.7), due to the localised source term selected. This example is used to check the performance of the mesh adaptive process described in section~\ref{sc:meshAdaptivity} and to highlight the capability of the error indicator~\eqref{eq:errorMeasure} to identify the region of interest in the domain, where the variation of the solution is localised.

\begin{figure}[!tb]
	\centering
	\subfigure[Initial mesh] {\includegraphics[width=0.2\textwidth]{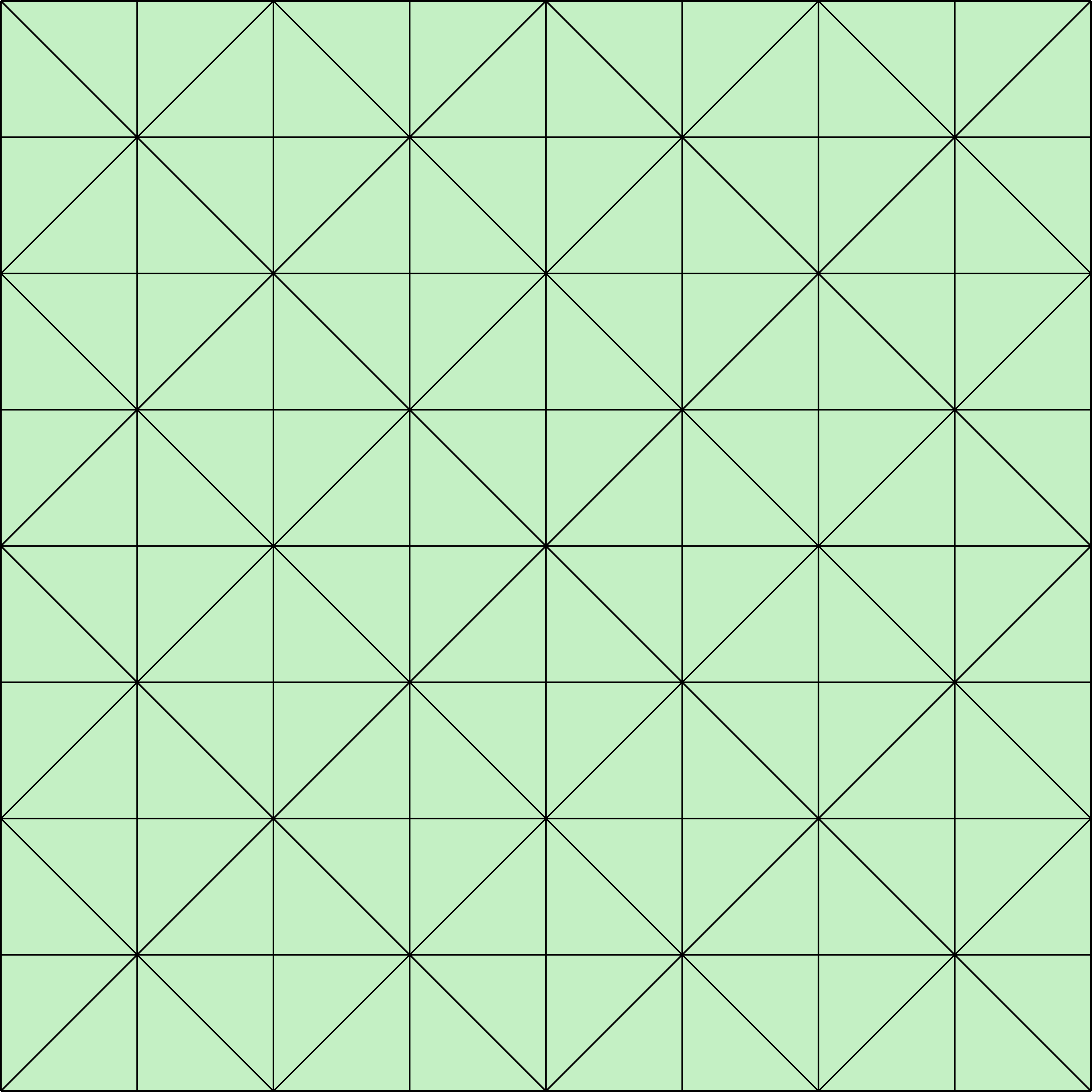}}
	\subfigure[Mesh 2] {\includegraphics[width=0.2\textwidth]{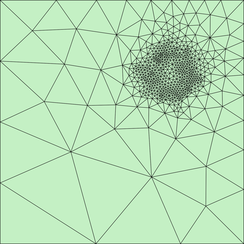}}
	\subfigure[Mesh 7] {\includegraphics[width=0.2\textwidth]{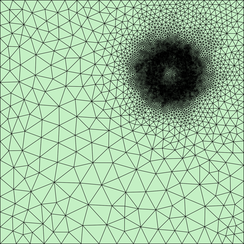}} \\	
	\subfigure[Initial mesh]  {\includegraphics[width=0.2\textwidth]{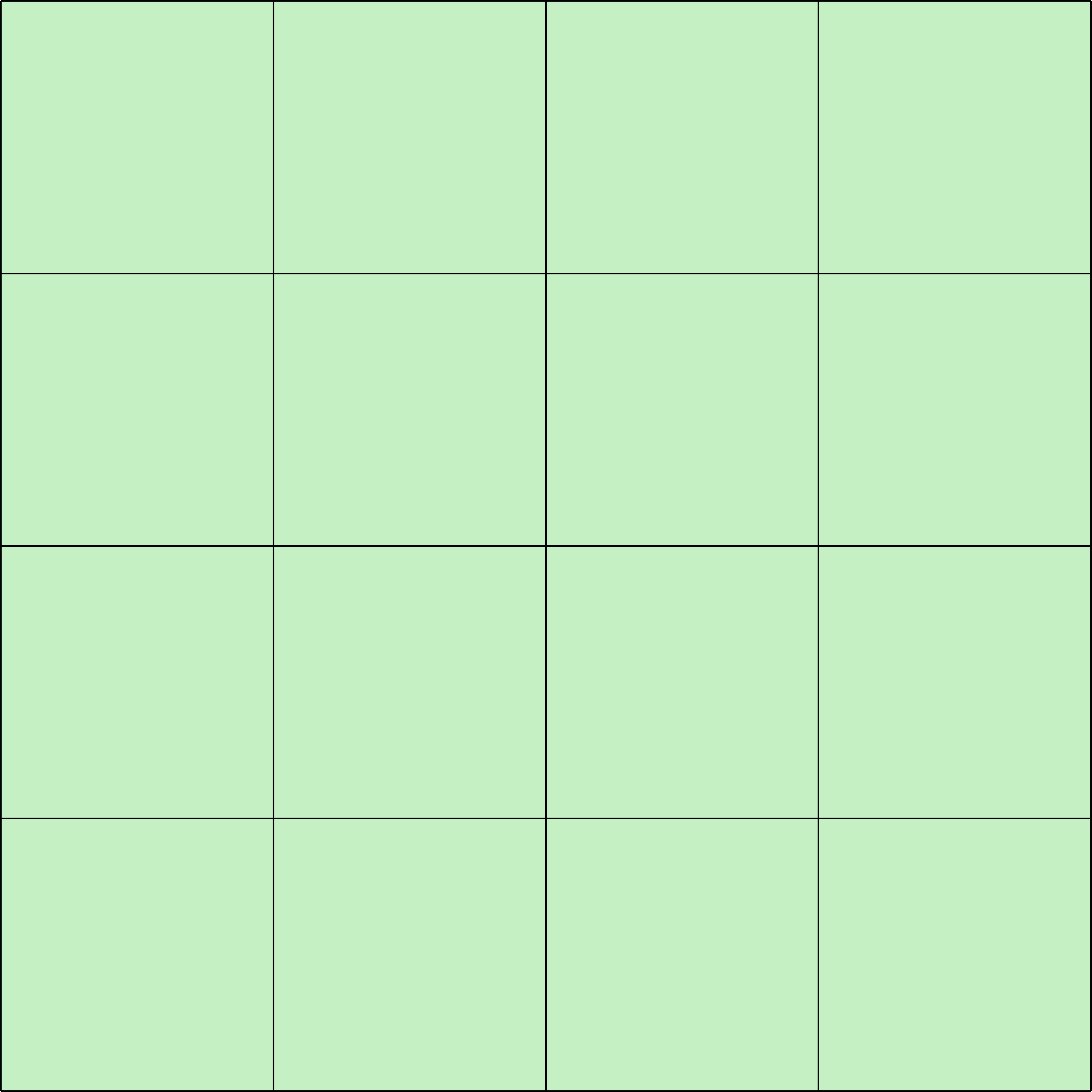}}
	\subfigure[Mesh 2]  {\includegraphics[width=0.2\textwidth]{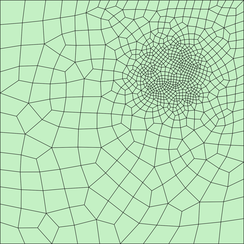}}
	\subfigure[Mesh 6]  {\includegraphics[width=0.2\textwidth]{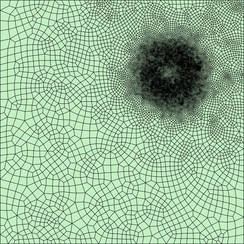}}
	\caption{Initial (left), intermediate (middle) and final (right) meshes generated by the automatic mesh adaptive procedure with a tolerance $\varepsilon = 10^{-2}$ for the Poisson problem using triangular (top) and quadrilateral (bottom) cells.}
	\label{fig:PoissonAdaptivityMesh}
\end{figure}

\begin{figure}[!tb]
	\centering				
	\subfigure[Initial mesh, $u^\star$]{\includegraphics[width=0.16\textwidth]{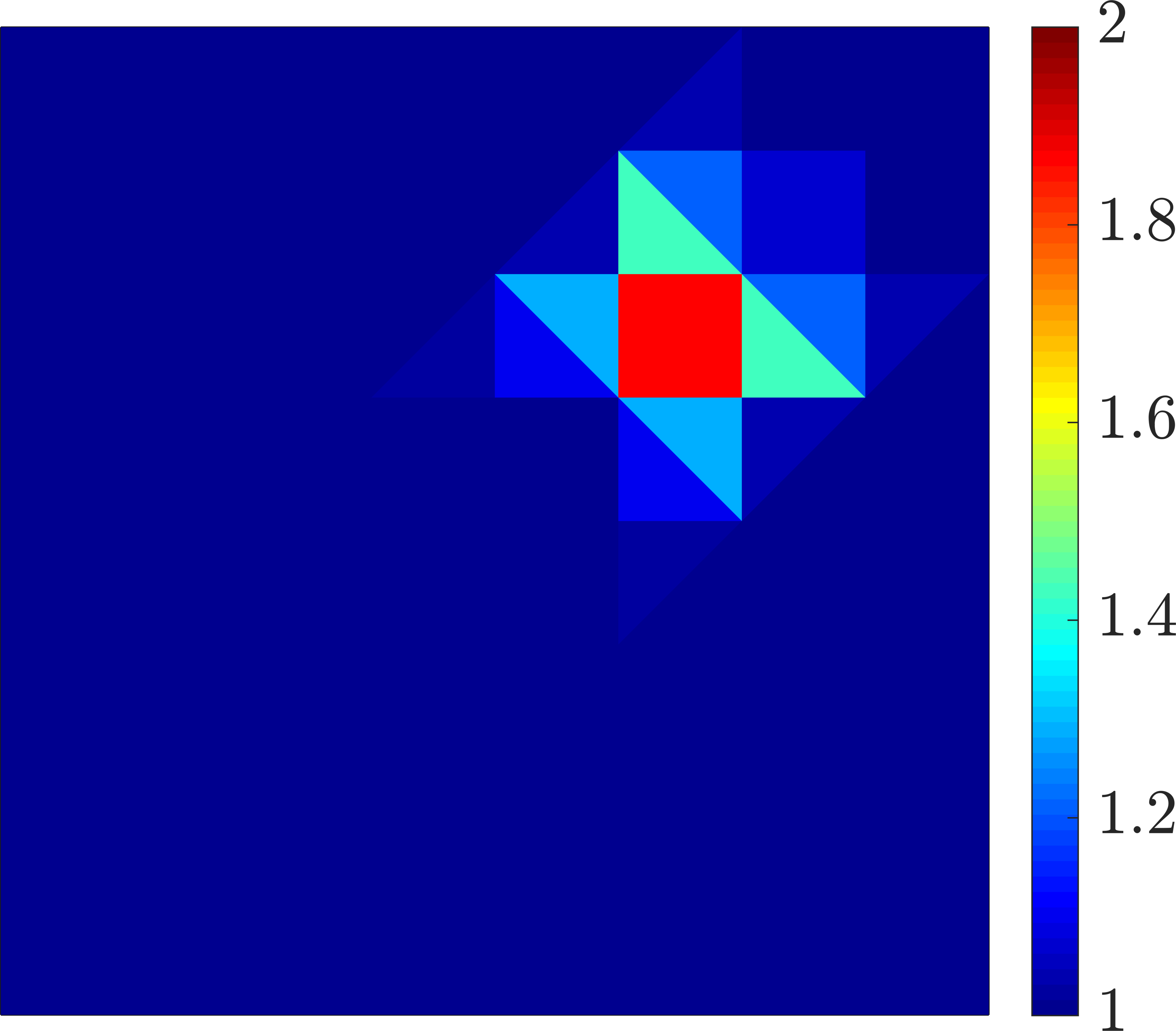}}
	\subfigure[Initial mesh, $u$]{\includegraphics[width=0.16\textwidth]{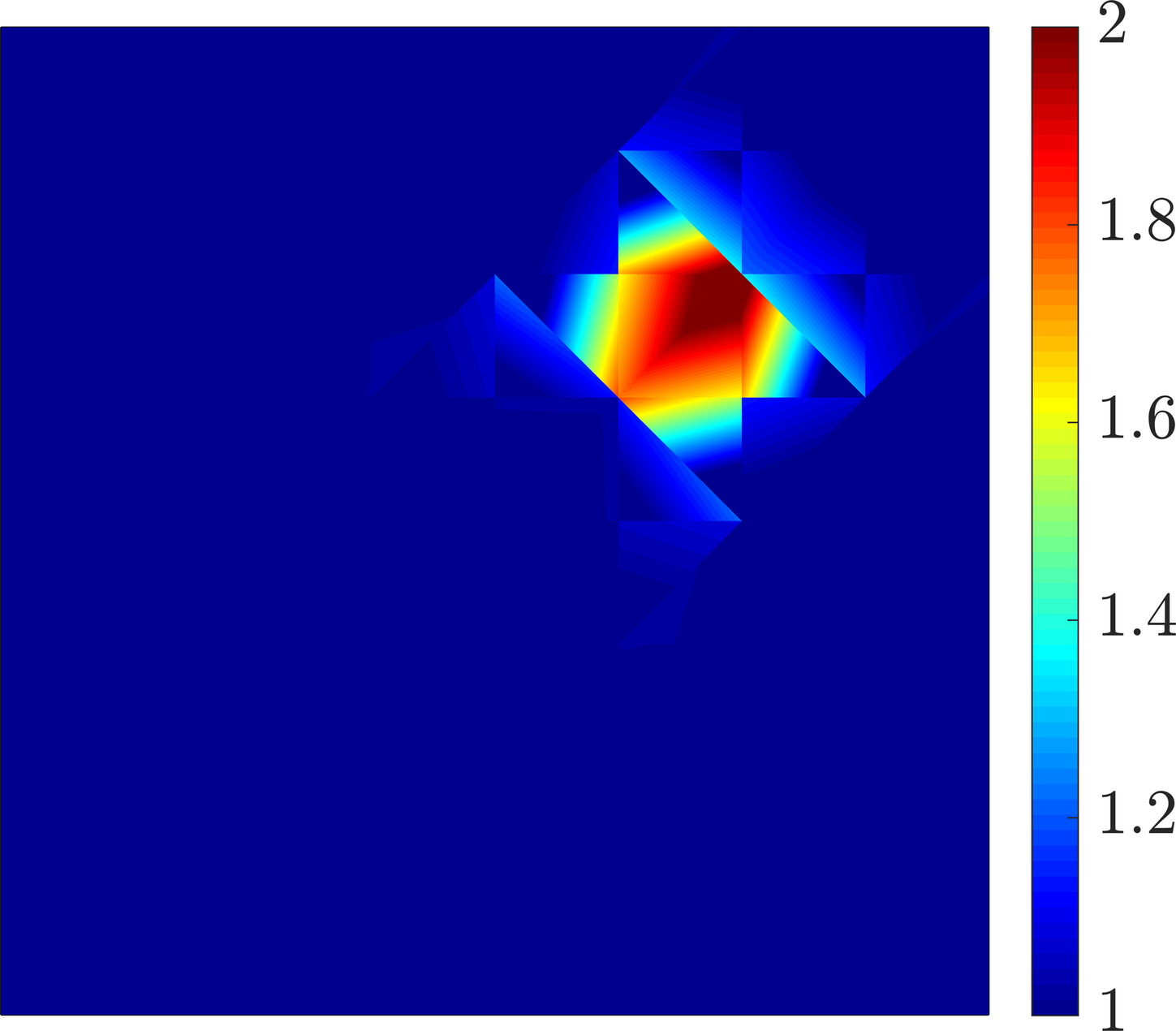}}
	\subfigure[Mesh 2, $u^\star$]{\includegraphics[width=0.16\textwidth]{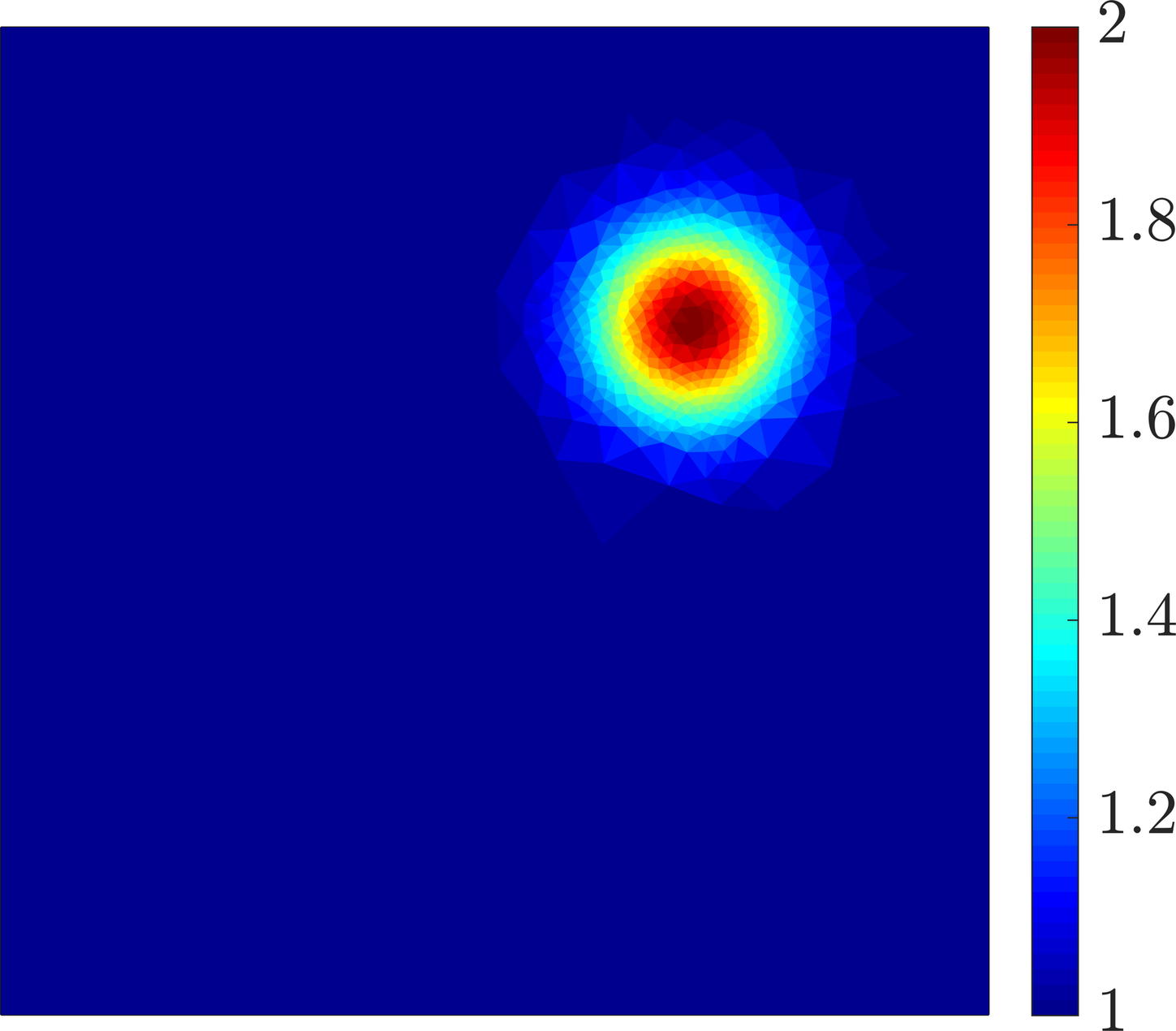}}
	\subfigure[Mesh 2, $u$]{\includegraphics[width=0.16\textwidth]{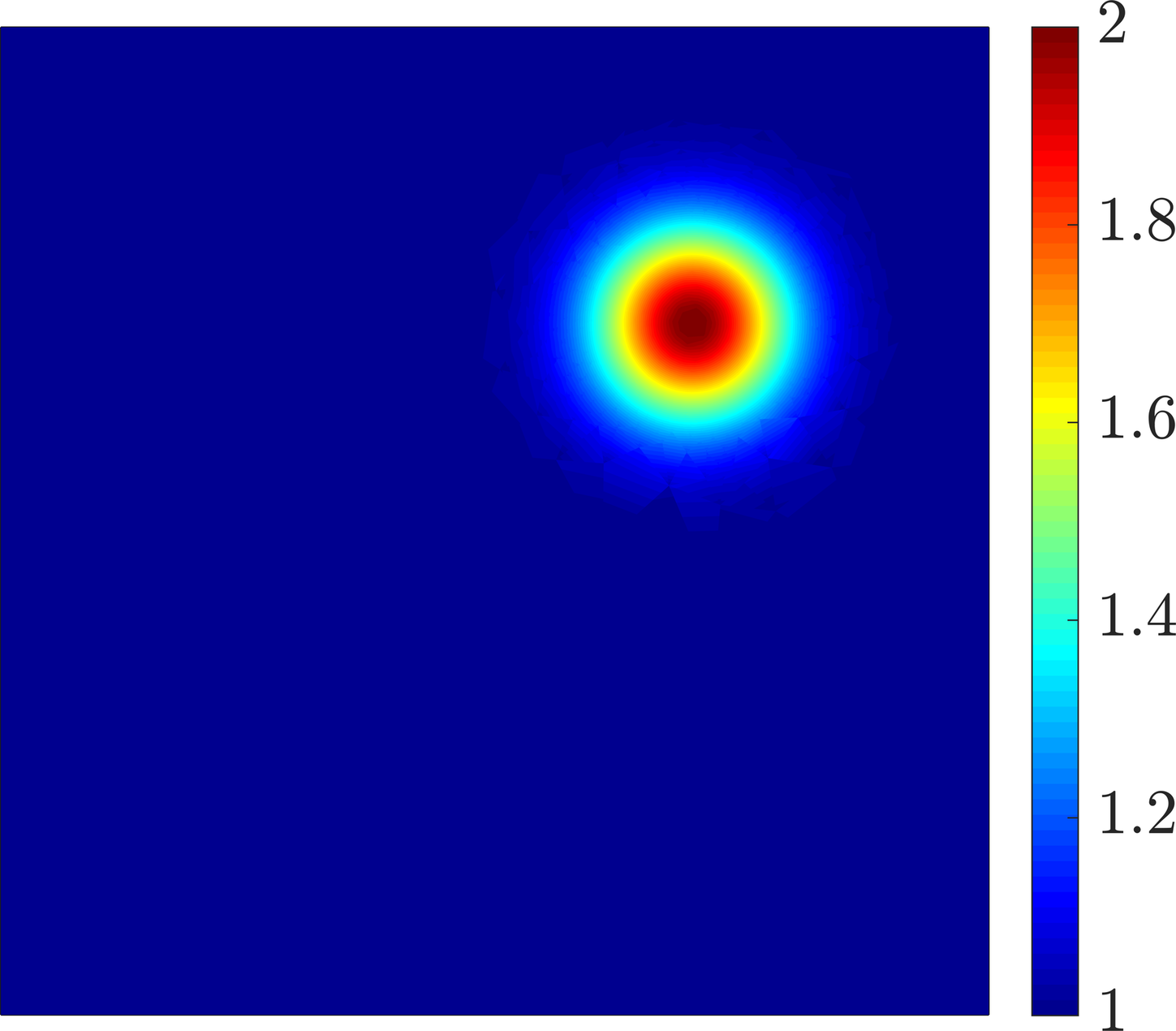}}
	\subfigure[Mesh 7, $u^\star$]{\includegraphics[width=0.16\textwidth]{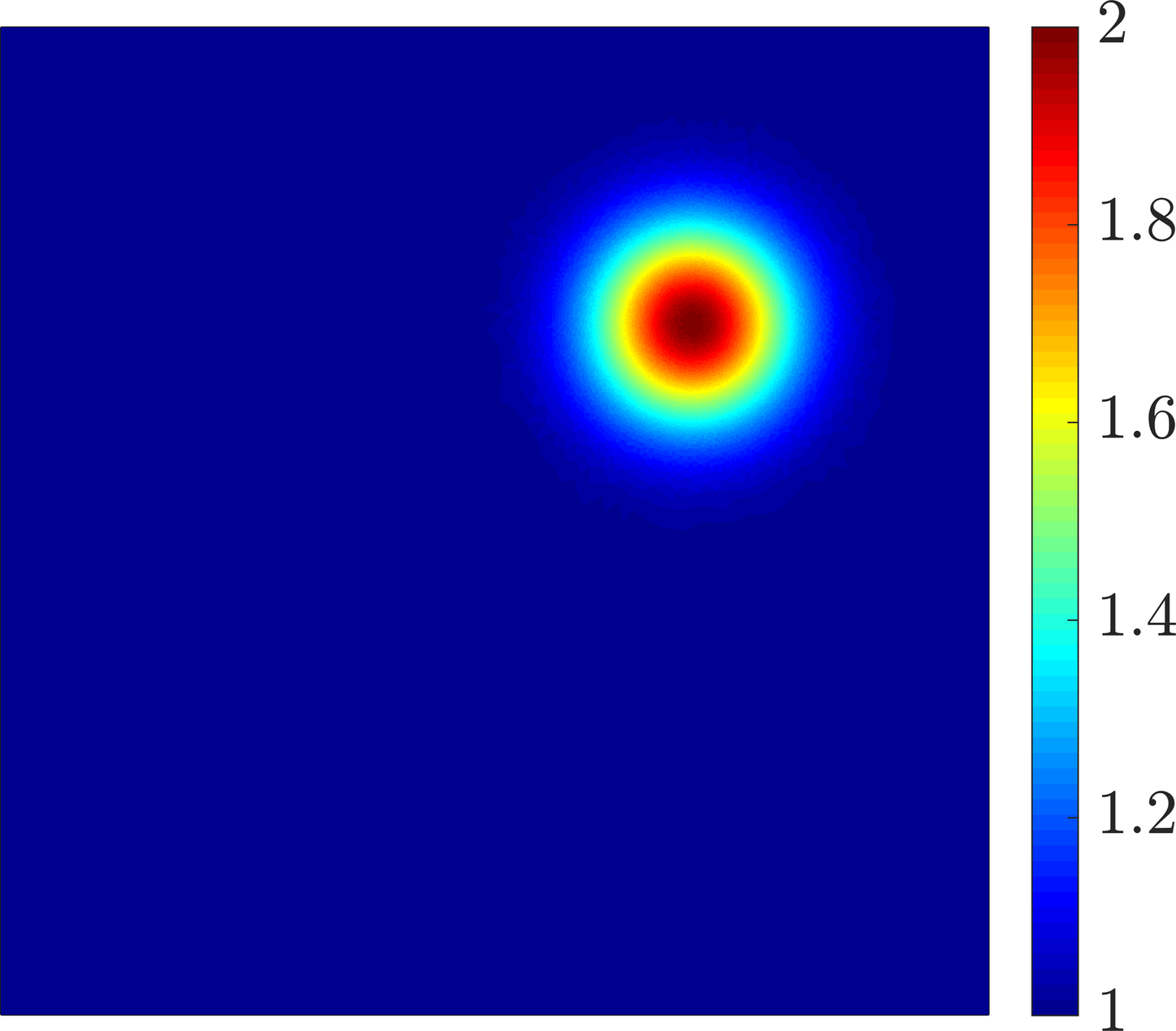}}
	\subfigure[Mesh 7, $u$]{\includegraphics[width=0.16\textwidth]{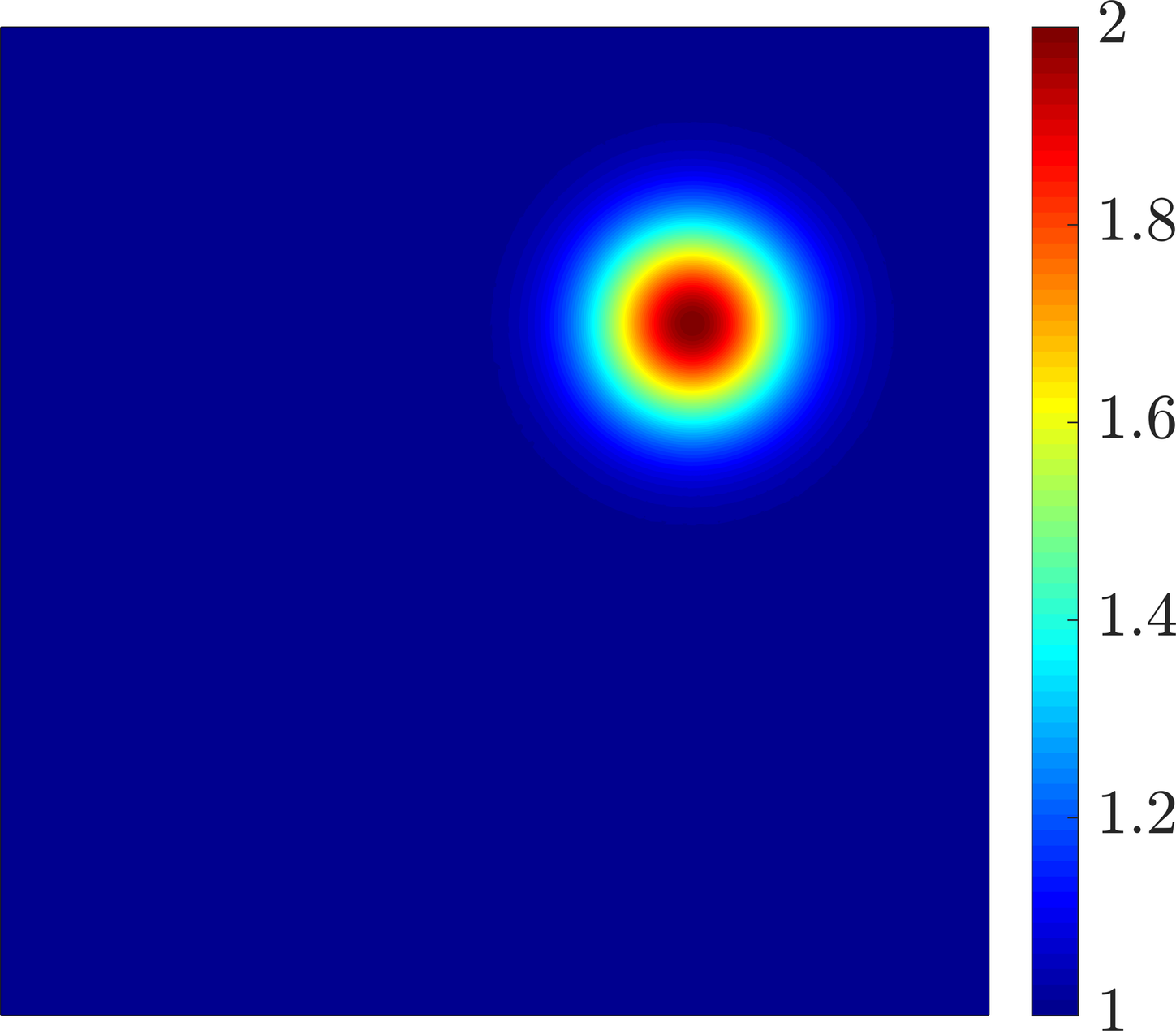}} \\
	\subfigure[Initial mesh, $u^\star$]{\includegraphics[width=0.16\textwidth]{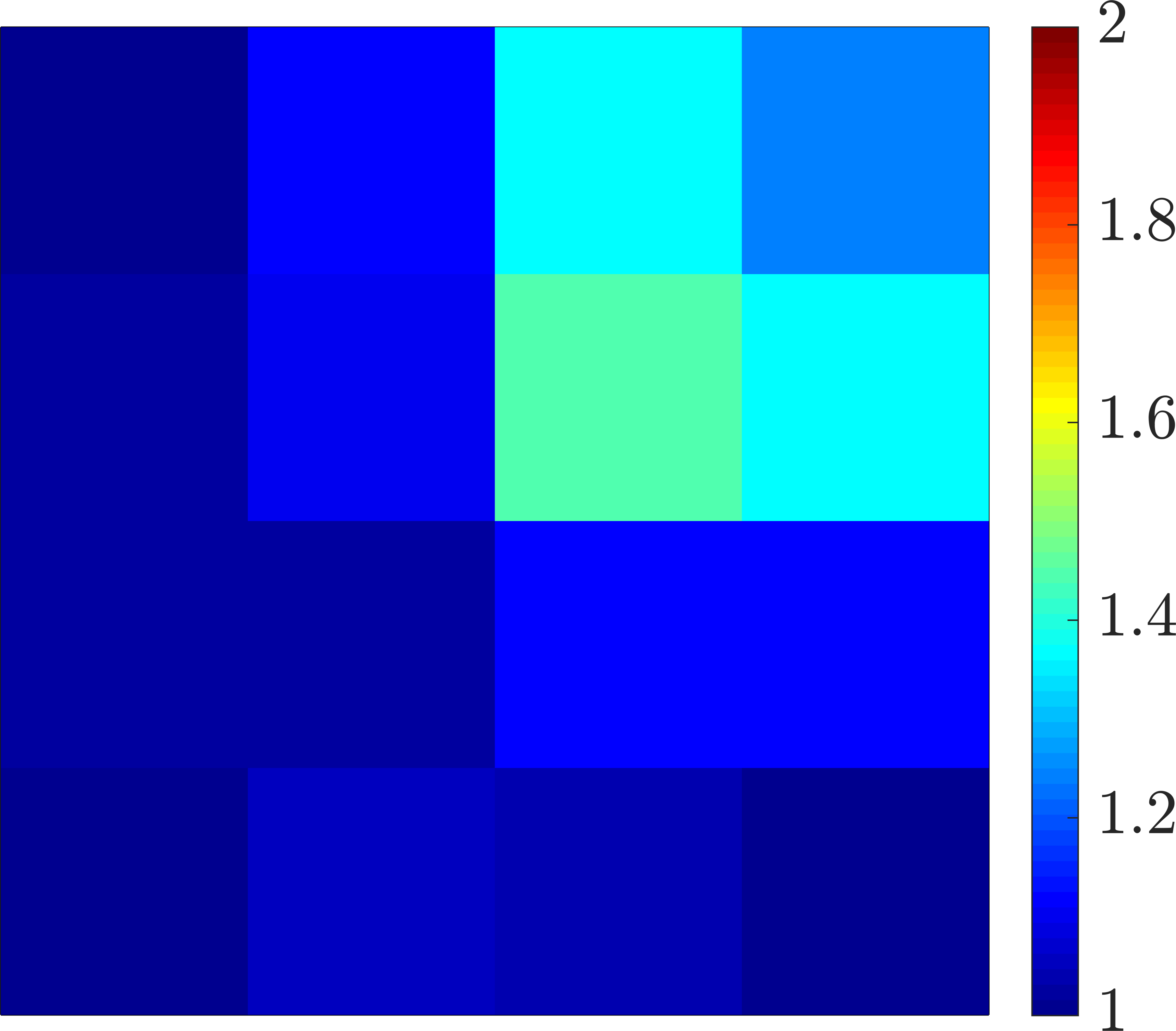}}
	\subfigure[Initial mesh, $u$]{\includegraphics[width=0.16\textwidth]{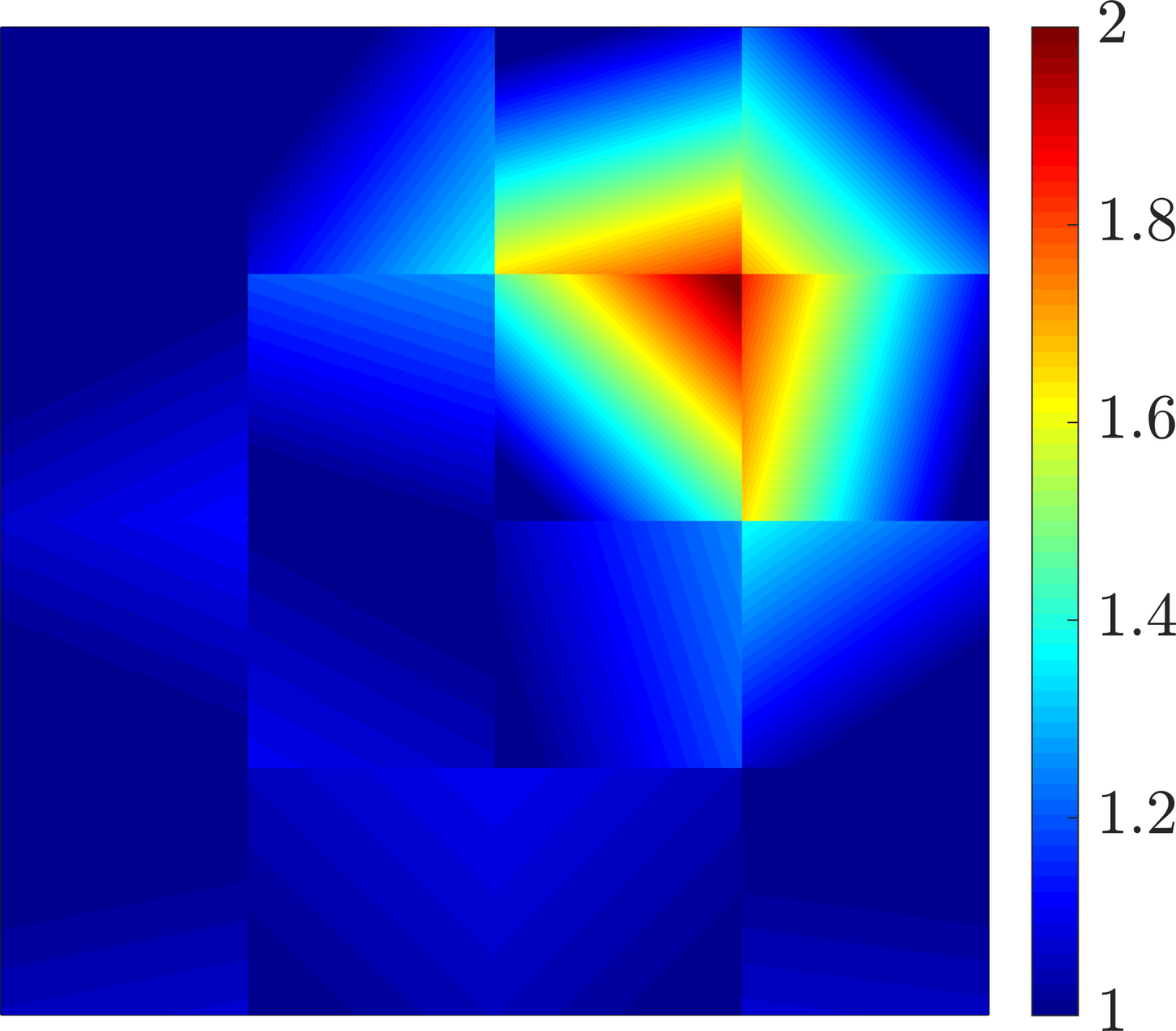}}
	\subfigure[Mesh 2, $u^\star$]{\includegraphics[width=0.16\textwidth]{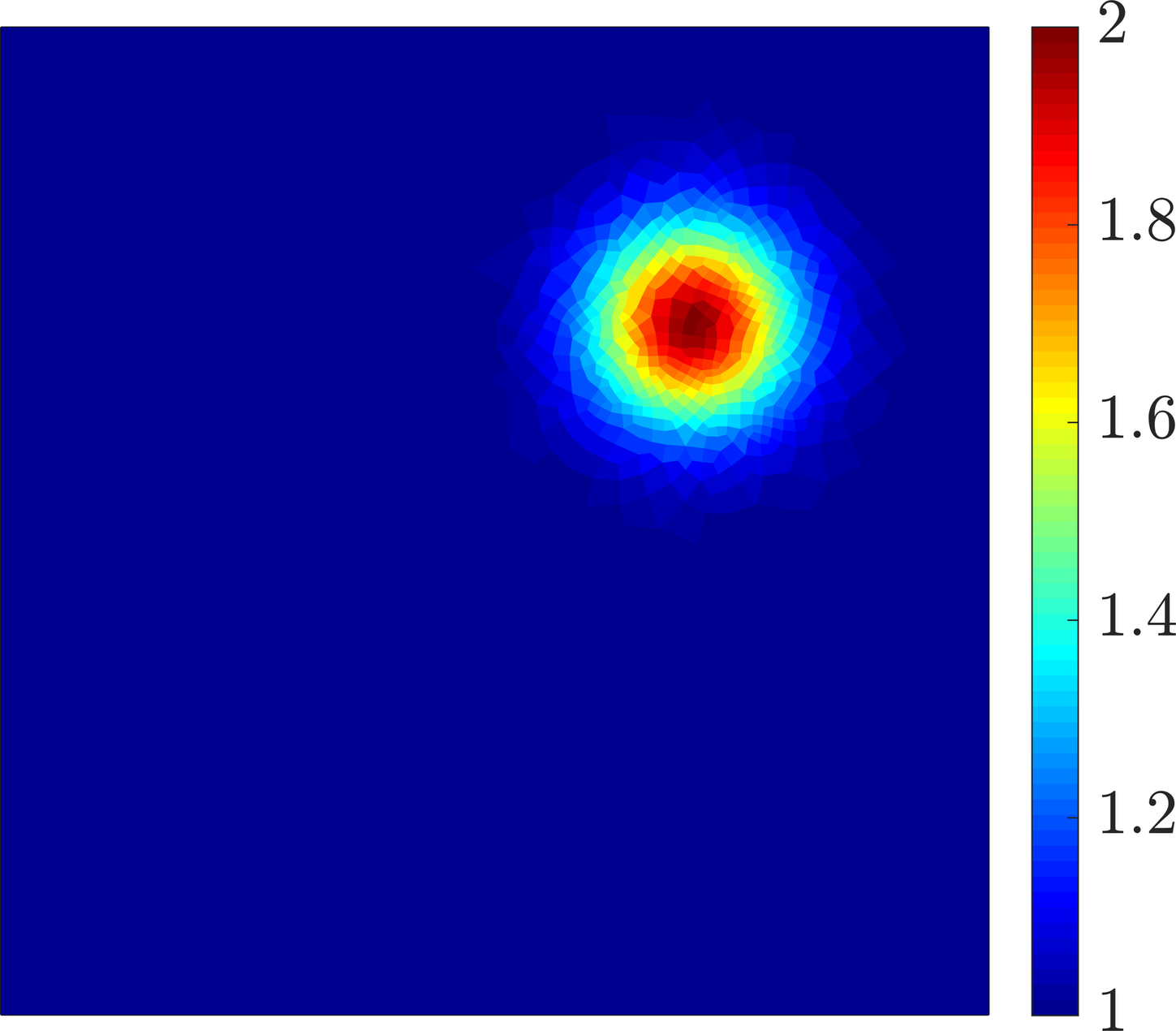}}
	\subfigure[Mesh 2, $u$]{\includegraphics[width=0.16\textwidth]{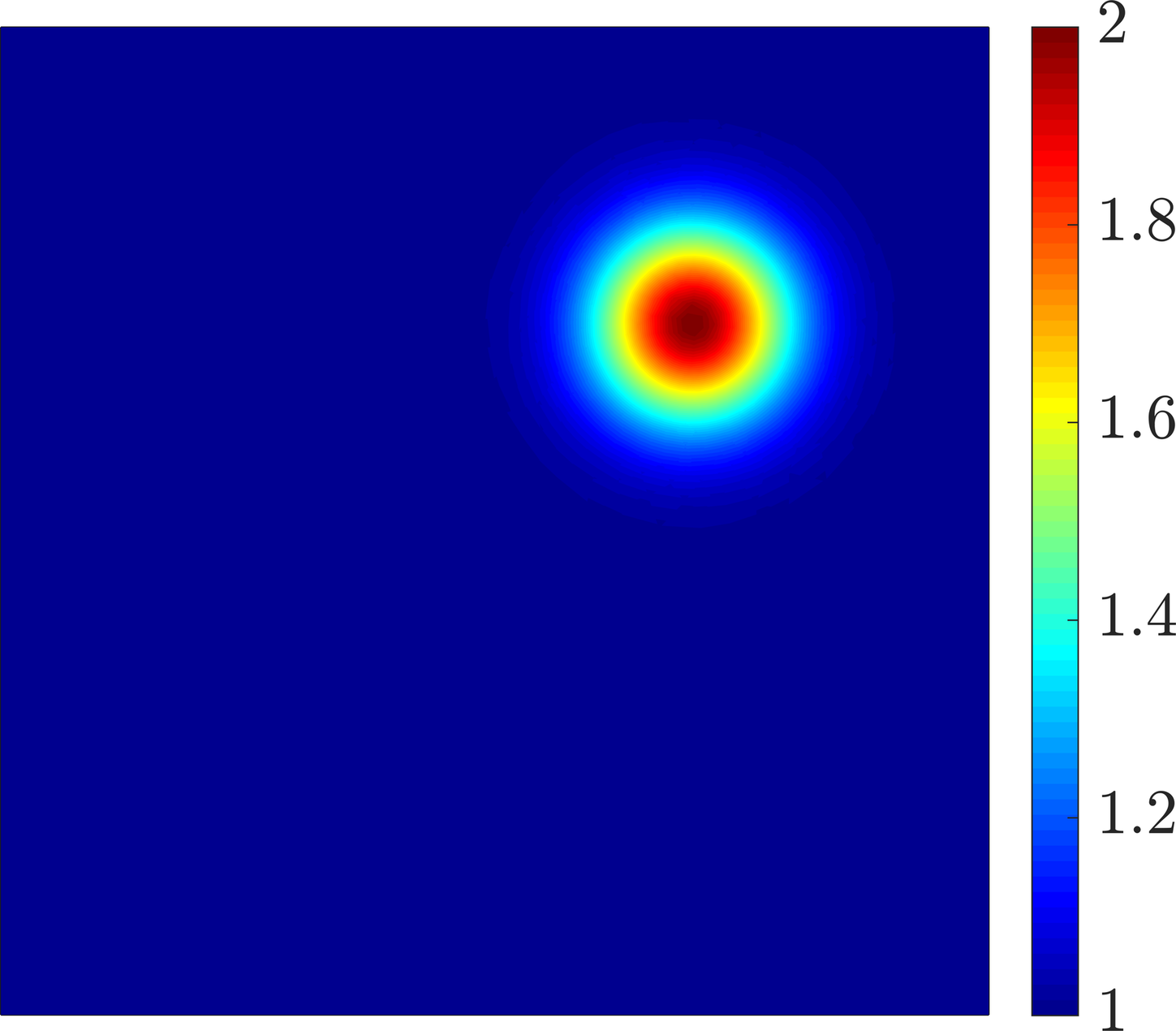}}
	\subfigure[Mesh 6, $u^\star$]{\includegraphics[width=0.16\textwidth]{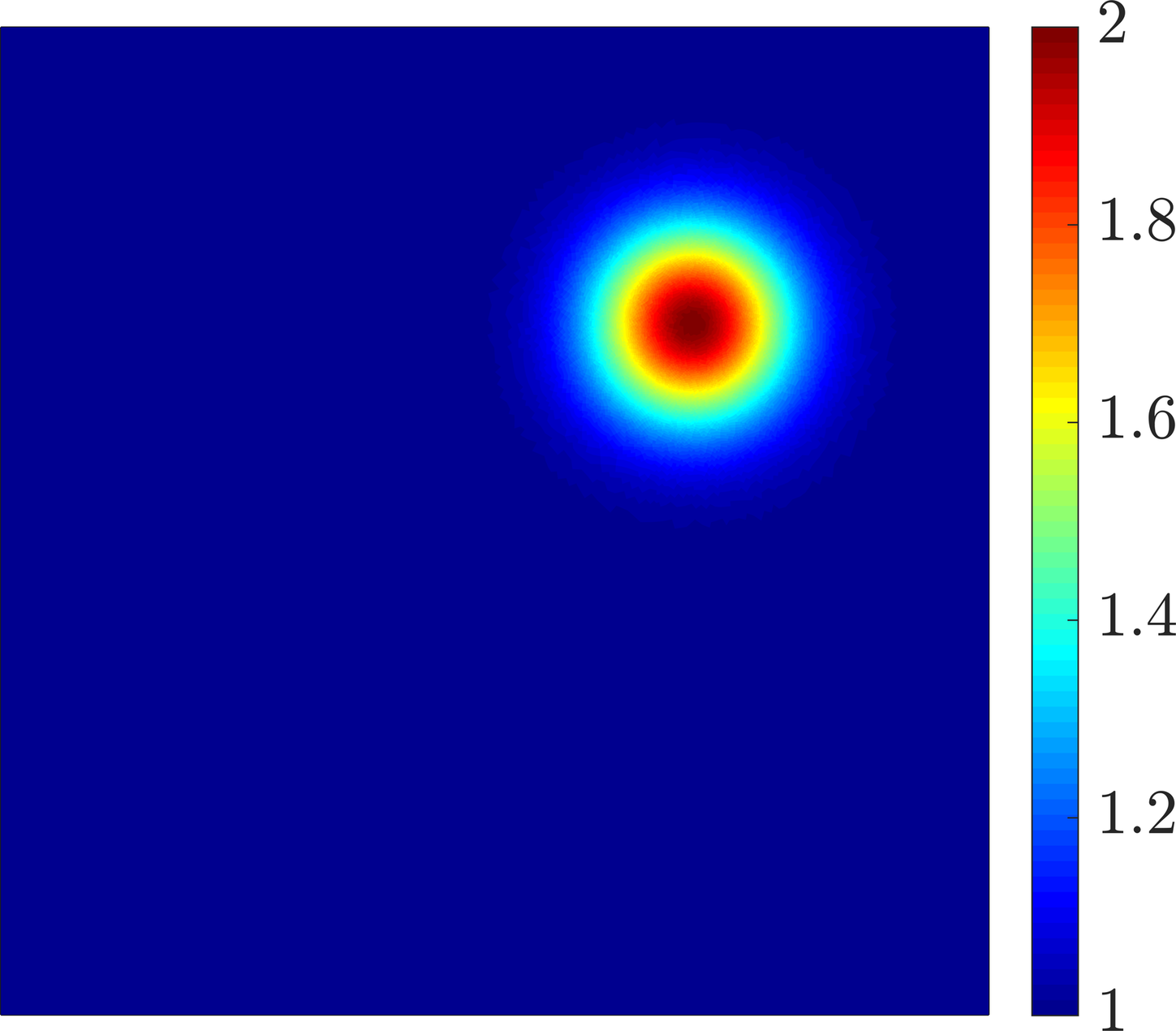}}
	\subfigure[Mesh 6, $u$]{\includegraphics[width=0.16\textwidth]{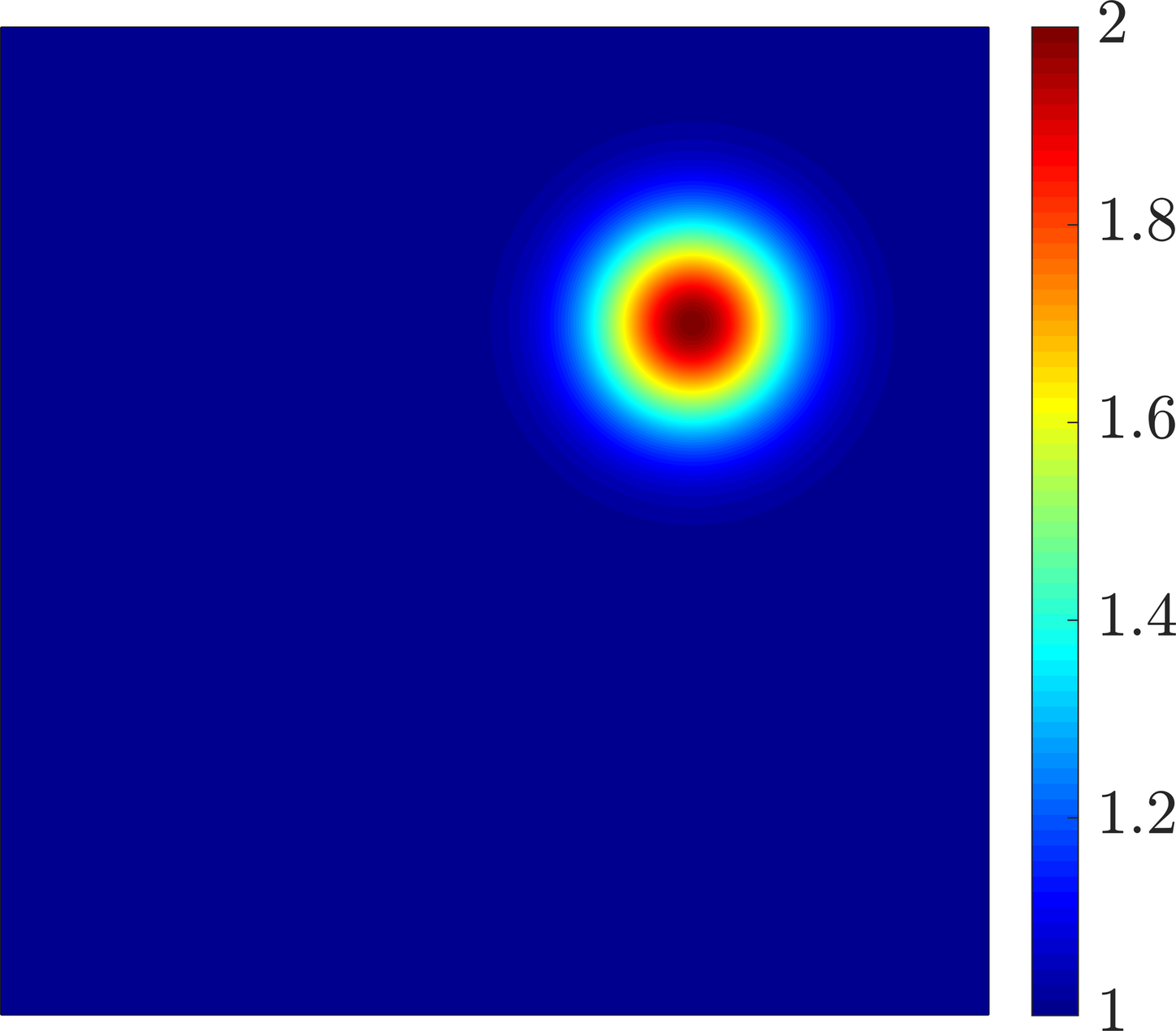}}
		\caption{Initial, intermediate and final FCFV first-order, $u^\star$, and second-order, $u$, approximations for the Poisson problem using triangular (top) and quadrilateral (bottom) meshes depicted in figure~\ref{fig:PoissonAdaptivityMesh} .}
	\label{fig:PoissonAdaptivitySol}
\end{figure}

Two mesh adaptive simulations are performed, with triangular and quadrilateral meshes and imposing a desired error in each cell of $\varepsilon = 10^{-2}$. The initial coarse triangular and quadrilateral meshes, shown in figures~\ref{fig:PoissonAdaptivityMesh} (a) and (d), have 128 and 16 cells respectively. The first-order solutions, $u^\star$, computed with the initial coarse meshes, are shown in figures~\ref{fig:PoissonAdaptivitySol} (a) and (g) and the second-order solutions are displayed in  figures~\ref{fig:PoissonAdaptivitySol} (b) and (h). 

After computing the error indicator, as detailed in section~\ref{sc:meshAdaptivity}, a desired size is determined for each cell of the coarse mesh. With this information, new meshes are generated, the first and second-order solutions are recomputed and the mesh adaptivity procedure is repeated. 
Figures~\ref{fig:PoissonAdaptivityMesh} (b) and (e) display the meshes after two adaptivity iterations. The triangular mesh has 1,271 cells, whereas the quadrilateral mesh has 871 cells. As it can be observed, the approximations computed at the second iteration of the mesh adaptive process already capture the main feature of the solution, which is a Gaussian profile centred at (0.7,0.7). The adaptive process converges in seven iterations for triangular meshes and six iterations for quadrilateral meshes. The triangular and quadrilateral final meshes, shown in figures~\ref{fig:PoissonAdaptivityMesh} (c) and (f), have 14,722 and 17,836 cells, respectively. The corresponding first-order solutions are displayed in figures~\ref{fig:PoissonAdaptivitySol} (e) and (k), whereas the final second-order solutions are reported in figures~\ref{fig:PoissonAdaptivitySol} (f) and (l).

\begin{figure}[!tb]
	\centering
	\subfigure[Triangles] {\includegraphics[width=0.32\textwidth]{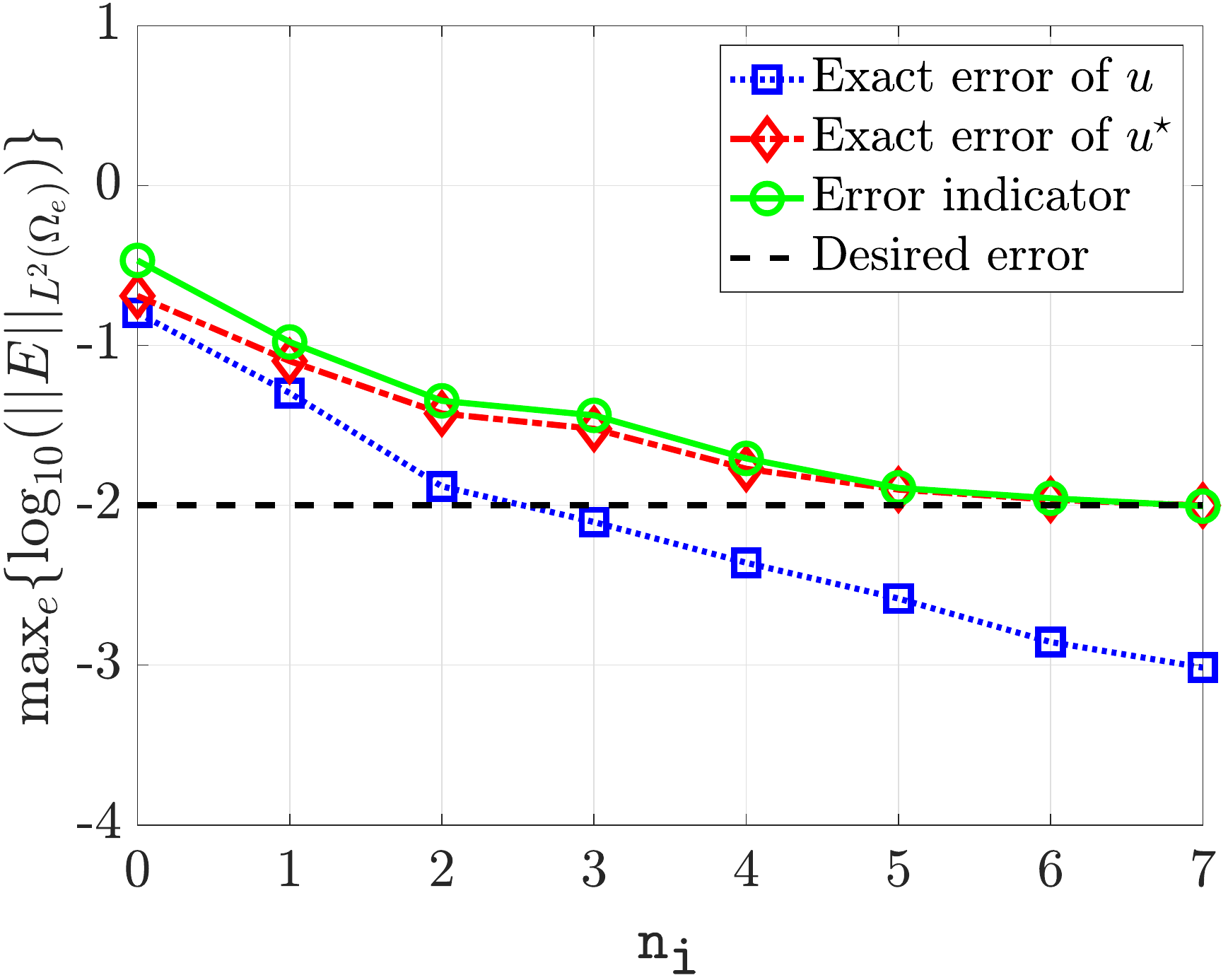}}
	\subfigure[Quadrilaterals] {\includegraphics[width=0.32\textwidth]{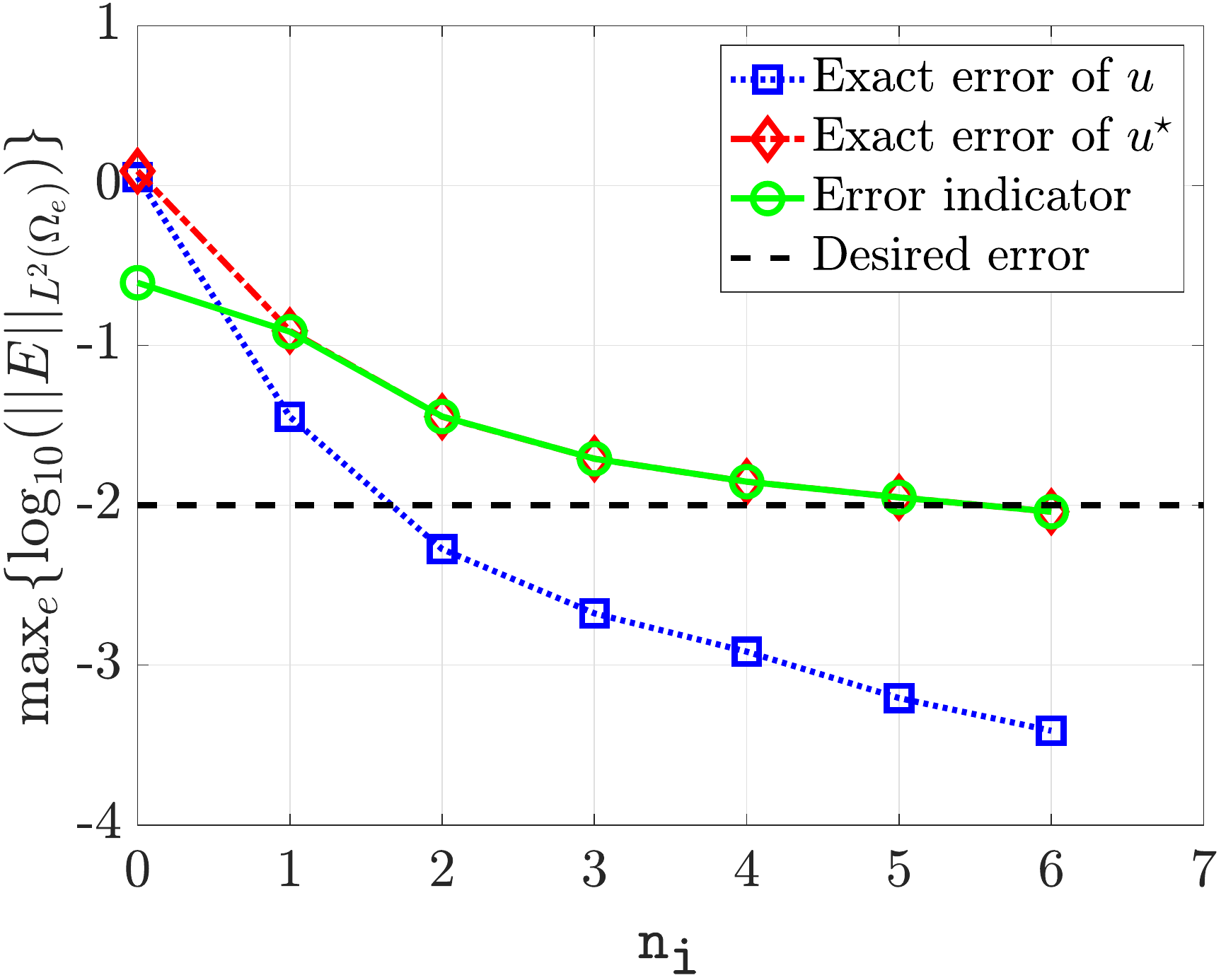}}
	\subfigure[Indicator efficiency] {\includegraphics[width=0.32\textwidth]{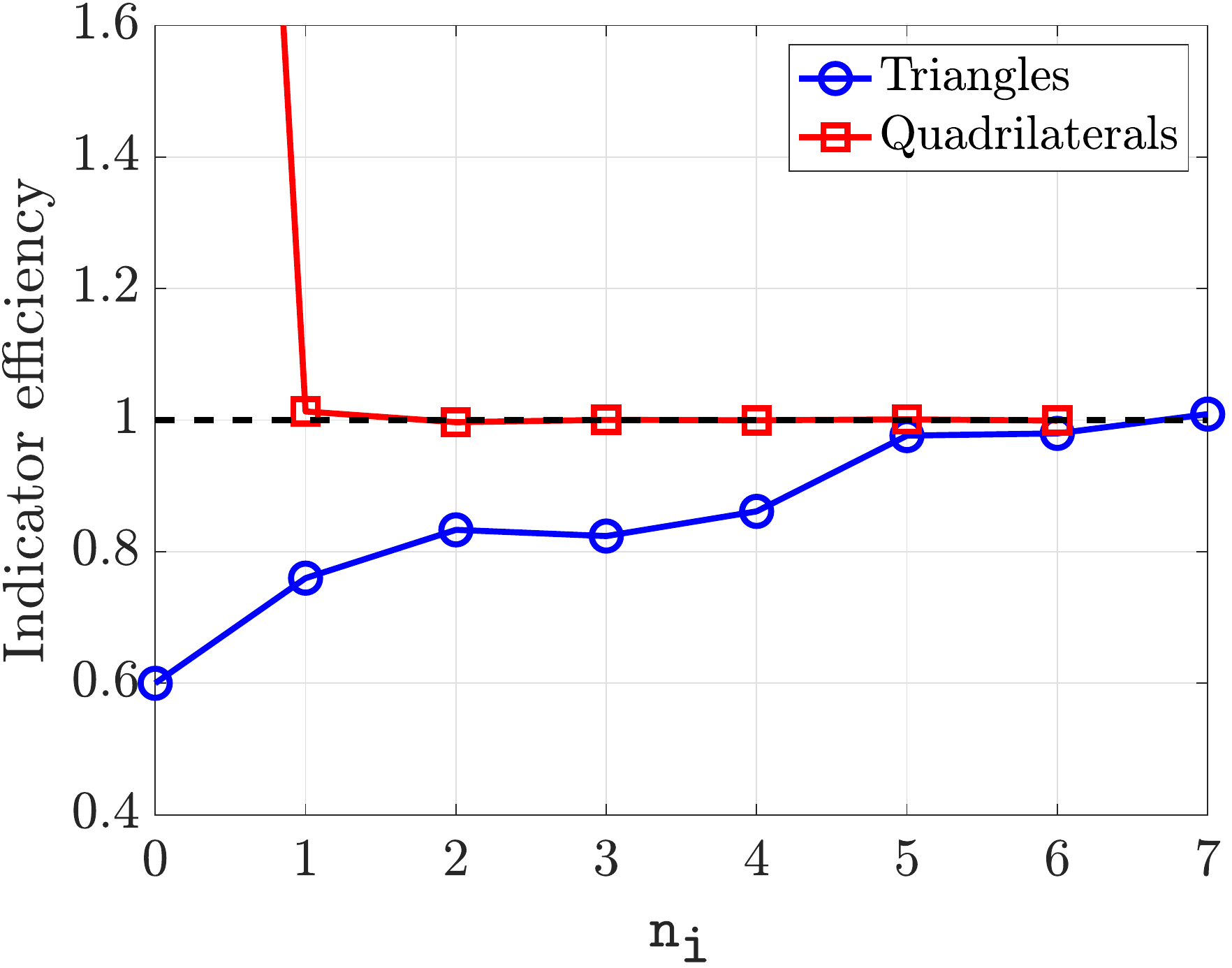}}	
	\caption{Maximum values of the error indicator and the exact error over all the cells as a function of the number of iterations of the mesh adaptive procedure using (a) triangular and (b) quadrilateral cells. (c) Indicator efficiency using triangular and quadrilateral meshes.}
	\label{fig:PoissonAdaptivityIterations}
\end{figure}

To further analyse these results, figures~\ref{fig:PoissonAdaptivityIterations} (a) and (b) show the evolution of the maximum values of the error indicator and the exact error over all the cells as a function of the number of iterations of the mesh adaptive procedure, $\niadapt$, for triangular and quadrilateral meshes respectively.
The results clearly show that the error indicator devised in section~\ref{sc:meshAdaptivity} produces a very accurate estimate of the error of the approximation $u^\star$, computed by solving an inexpensive extra local problem given by equation~\eqref{eq:Poisson-locU}, for both types of meshes. A slight difference is observed between the error indicator and the exact error of $u^\star$ for triangular meshes for the first four iterations of the mesh adaptivity process, whereas for quadrilateral cells a perfect agreement is observed. 

To quantify this difference, the so-called \textit{indicator efficiency}, computed as the ratio between the exact error of $u^\star$ and the error indicator~\eqref{eq:errorMeasure}, is reported in figure~\ref{fig:PoissonAdaptivityIterations} (c). On the first iteration, the efficiency for triangular meshes is 0.76 and slowly improves during the adaptive process, being 0.82 in the third iteration, 0.98 in the fifth iteration and 1.01 in the final iteration. For quadrilateral cells, the indicator efficiency is already 1.01 in the first iteration and takes a value of 1.00 from the second to the sixth iteration. This clearly indicates that in this example, the use of quadrilateral meshes is beneficial compared to triangular meshes. Not only less iterations of the mesh adaptivity process are required but, in addition, the final error of the approximate solution is lower when using quadrilateral cells.

It is worth noting that despite the mesh adaptivity process is driven by computing an error indicator for the approximate solution $u^\star$, the error of the more accurate approximation $u$ also decreases monotonically during the adaptive process, as shown in figures~\ref{fig:PoissonAdaptivityIterations} (a) and (b). This is expected due to the higher accuracy of the approximation computed with the proposed second-order FCFV, when compared to the accuracy of the first-order solution obtained by solving the extra local problem of equation~\eqref{eq:Poisson-locU}. The difference in accuracy can be observed in figure~\ref{fig:CPUtimeStokesK0vsK0K1}, where the error of the first and second-order FCFV is compared for a Stokes problem in two and three dimensions.

\subsection{Three dimensional Stokes flow around complex microswimmers}
\label{sc:Stokes3Dswimmer} 

The last example considers the three dimensional Stokes flow around microswimmers. The geometry of the microswimmers, taken from~\cite{keaveny2013optimization}, is given in parametric form as
\begin{equation} \label{eq:swimeerGeometry}
\bm{S}(\lambda, \theta) = \bm{C}(\lambda) + R_n(\lambda) \sin(\theta) \bn_1 + R_b(\lambda) \cos(\theta) \bn_2, \qquad (\lambda, \theta) \in [-L,L] \times [0, 2\pi) ,
\end{equation}
where the curve $\bm{C}$ is a parametrisation of the centreline of the swimmer, namely
\begin{equation} \label{eq:swimmerCentreline}
\bm{C}(\lambda) = \left( \beta \cos( \kappa\lambda), \beta \sin( \kappa\lambda), \alpha \lambda \right) .
\end{equation}
In~\eqref{eq:swimeerGeometry}, $\bn_1$ and $\bn_2$ denote the unit normal vectors to the centreline tangent and serve as the short and long axis respectively of the propeller cross-section. More precisely, they are defined as
\begin{equation} \label{eq:swimmerN1N2}
\bn_1 = \cos \big( \gamma F_1(\lambda) \big) \bm{N} + \sin \big( \gamma F_1(\lambda) \big) \bm{B}, 
\qquad 
\bn_2 = \cos \big( \gamma F_1(\lambda) \big) \bm{B} - \sin \big( \gamma F_1(\lambda) \big) \bm{N}, 
\end{equation}
in terms of the Serret-Frenet normal $\bm{N}$ and bi-normal $\bm{B}$.
The radii of the long and short axis of the propeller cross-sections of the swimmer, denoted by $R_b$ and $R_n$ respectively, are defined as
\begin{equation} \label{eq:swimmerRadius}
R_b(\lambda) = A_b \left(C_1 + C_2 F_0(\lambda) \right) \left( 1 - \lambda^8 \right)^{1/8}, \qquad R_n(\lambda) = \frac{1}{4} R_b(\lambda),
\end{equation}
where the function $F_s$ is defined, for $s=\{0,1\}$, as
\begin{equation} \label{eq:swimmerFs}
F_s(\lambda) = \frac{1}{2} \left[ 1 - \text{erf} \left( \frac{\lambda - \lambda_s}{\sqrt{2} \sigma} \right) \right].
\end{equation}
All the parameters appearing in the previous expressions are given in table~\ref{tab:swimmerParams}.
\begin{table}[hbt]
	\centering
	\begin{tabular}[hbt]{|c|c|c|c|c|c|c|c|c|c|}
		\hline
		$L$ & $A_b$ & $C_1$ & $C_2$ & $\alpha$ & $\beta$ & $\kappa$ & $\lambda_0$ & $ \lambda_1$ & $\sigma$ \\
		\hline & & & & & & & & &
		\\ [-1em] 
		\hline
		1 & $L/27$ & 1.75 & 2.75 & 0.7 & $(1-\alpha^2)^{1/2}/\kappa$ & $4\pi/L$ & $(1 - 9/54.2)L$ & $(1 - 11/54.2)L$ & $0.02L$ \\
		\hline
	\end{tabular}
	\caption{Parameters used to define the geometry of the microswimmers.}
	\label{tab:swimmerParams}
\end{table}

\begin{figure}[!bt]
	\centering
	\subfigure[$\gamma=0$]     {\includegraphics[width=0.2\textwidth]{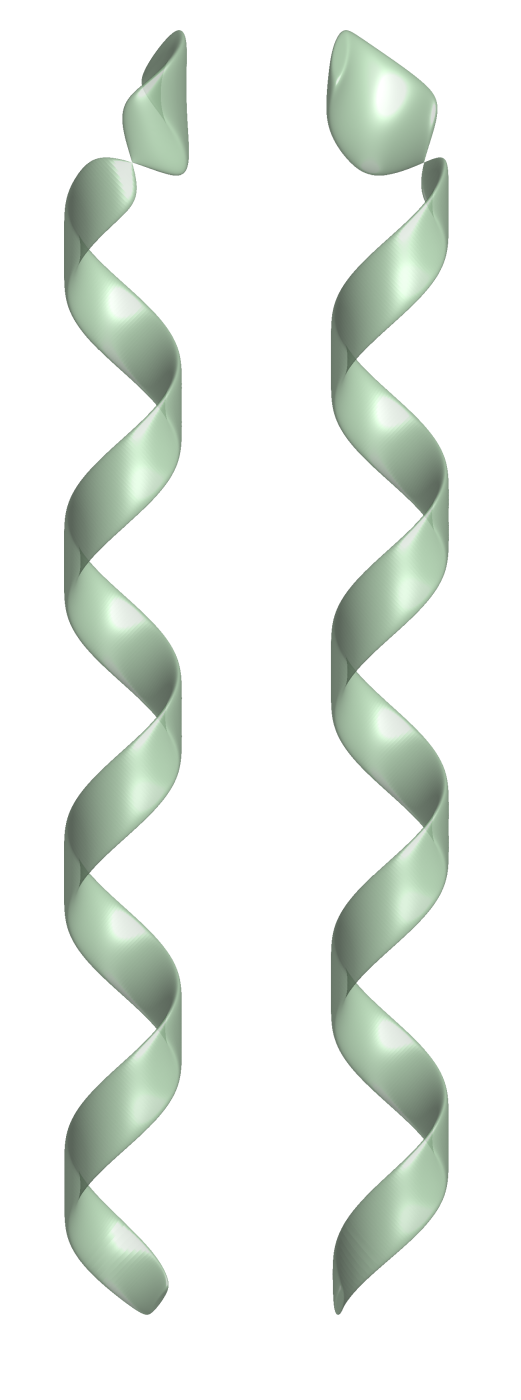}}	
	\hspace{1cm}
	\subfigure[$\gamma=\pi/4$] {\includegraphics[width=0.2\textwidth]{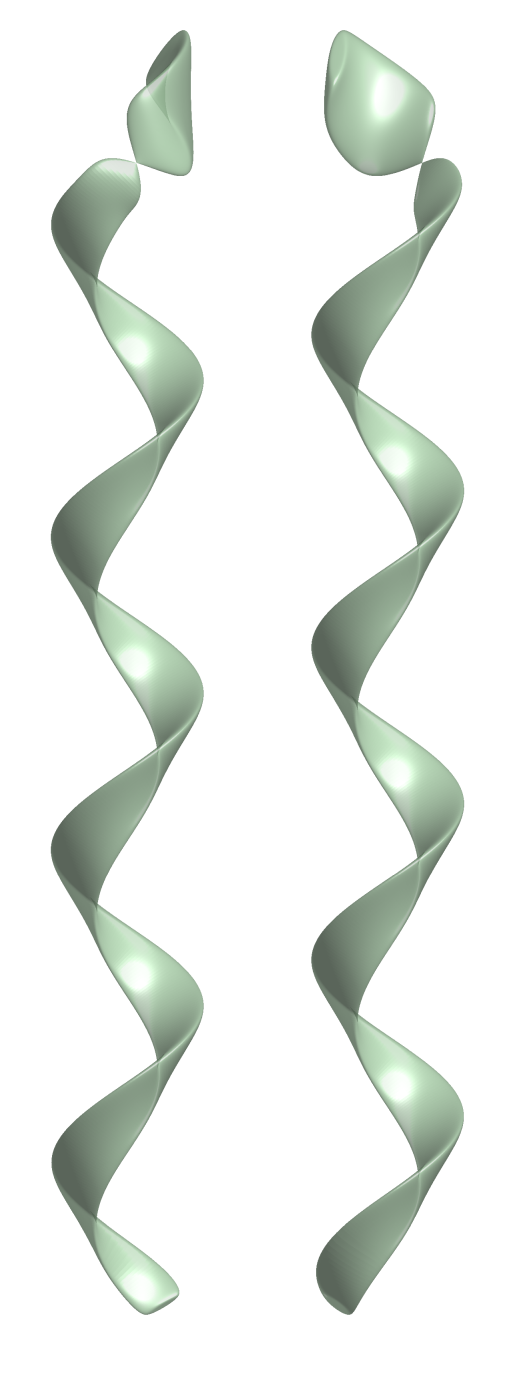}}	
	\hspace{1cm}
	\subfigure[$\gamma=\pi/2$] {\includegraphics[width=0.2\textwidth]{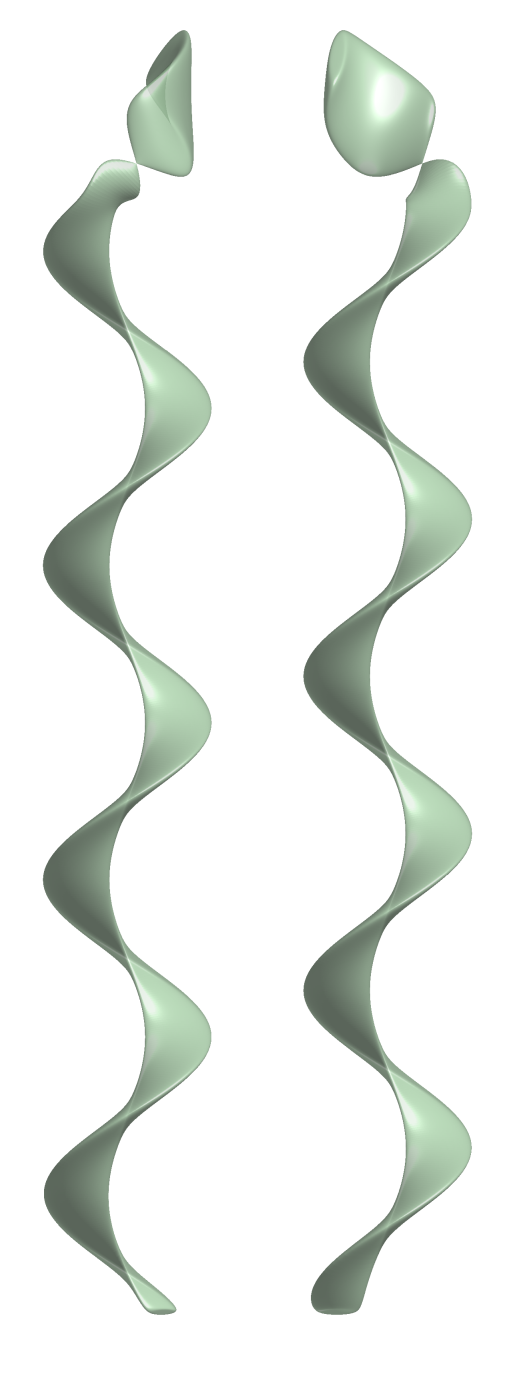}}	
	\caption{Geometry of three microswimmers for the parameters given in table~\ref{tab:swimmerParams}, showing two perspectives of the same geometry.}
	\label{fig:swimmerGeometry}
\end{figure}

Three microswimmers, obtained by varying the parameter $\gamma$ in equation~\eqref{eq:swimmerN1N2}, are considered to demonstrate the potential of the automatic mesh adaptivity framework proposed in section~\ref{sc:meshAdaptivity}. The geometries of the three cases considered are displayed in figure~\ref{fig:swimmerGeometry}.

To perform the Stokes flow simulation, the microswimmers of volume $\mathscr{S}_\gamma$ are placed in the centre of a prismatic channel $\mathscr{B} = [-L_1,L_1] \times [-L_2,L_2] \times [-L_3,L_3]$, with $L_1=L_2=\SI{17.14}{\mu m}$ and $L_3=2L_1$. The resulting computational domain is given by $\Omega_\gamma = \mathscr{B} \setminus \mathscr{S}_\gamma$. A paraboloid velocity profile $\bu_D(x_1,x_2,x_3) = (0,0,-4.1(L_1^2-x_1^2)(L_2^2-x_2^2)/L_1^2 L_2^2)\SI{}{\mu m/s}$ is imposed on the inlet, at $x_3=L_3$, and a free-traction condition is enforced on the outlet, at $x_3=-L_3$. On the remaining lateral walls of $\mathscr{B}$ and on the surface of the microswimmer, a no-slip boundary condition, corresponding to material walls, is enforced. The kinematic viscosity is taken as $\nu = \SI{2.65}{mm^2/s}$, which is the value corresponding to blood at $\SI{37}{\celsius}$. 

\begin{figure}[!bt]
	\centering
	\subfigure[Mesh 1, $\gamma=0$]     {\includegraphics[width=0.2\textwidth]{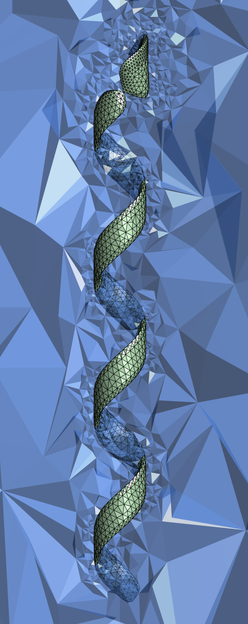}}	
	\hspace{1cm}
	\subfigure[Mesh 1, $\gamma=\pi/4$] {\includegraphics[width=0.2\textwidth]{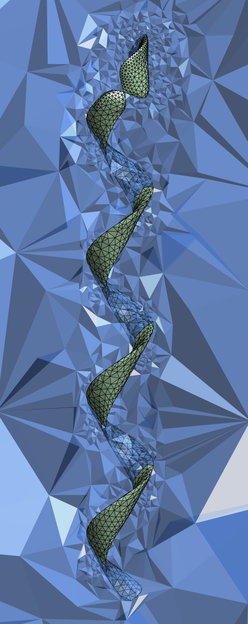}}	
	\hspace{1cm}
	\subfigure[Mesh 1, $\gamma=\pi/2$] {\includegraphics[width=0.2\textwidth]{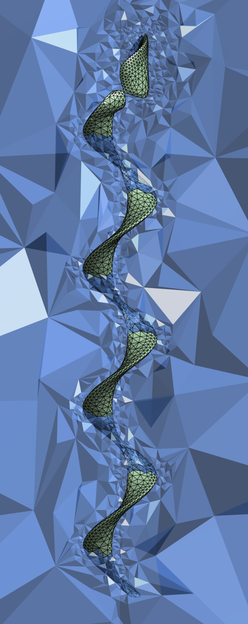}}	
	\caption{Detail of the initial tetrahedral meshes around the microswimmers of figure~\ref{fig:swimmerGeometry}.}
	\label{fig:swimmerInitialMeshes}
\end{figure}

The initial meshes displayed in figure~\ref{fig:swimmerInitialMeshes} are generated, using the technique described in~\cite{OHmeshing}, to start the automatic adaptivity process. 
To ensure a good geometric representation, the initial meshes are generated by imposing a desired element size of $\SI{0.04}{\mu m}$ along the curves $\bm{S}(\lambda, 0)$ and $\bm{S}(\lambda, \pi)$, for $\lambda \in [-L,L]$. The desired element size in the rest of the domain is $\SI{0.5}{\mu m}$. The resulting meshes have 80,024, 80,294 and 77,533 elements for $\gamma=0$, $\gamma=\pi/4$ and $\gamma=\pi/2$ respectively. It is worth noting that, due to the complexity of the geometry, only tetrahedral meshes are considered in this example. 

\begin{figure}[!bt]
	\centering
	\subfigure[Mesh 2, $\gamma=\pi/4$] {\includegraphics[width=0.2\textwidth]{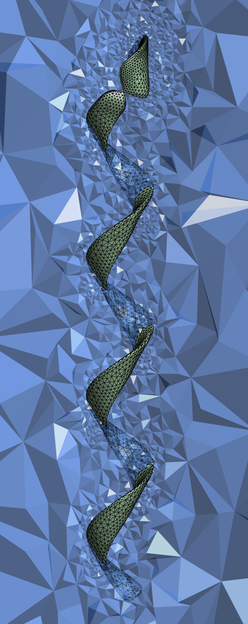}}	
	\hspace{1cm}
	\subfigure[Mesh 3, $\gamma=\pi/4$] {\includegraphics[width=0.2\textwidth]{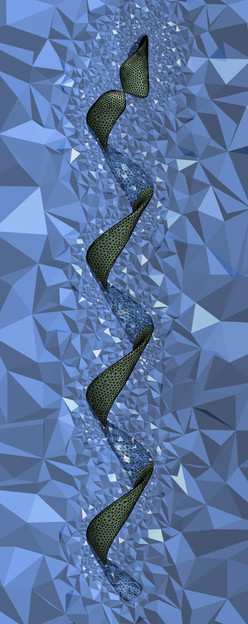}}	
	\hspace{1cm}
	\subfigure[Mesh 4, $\gamma=\pi/4$] {\includegraphics[width=0.2\textwidth]{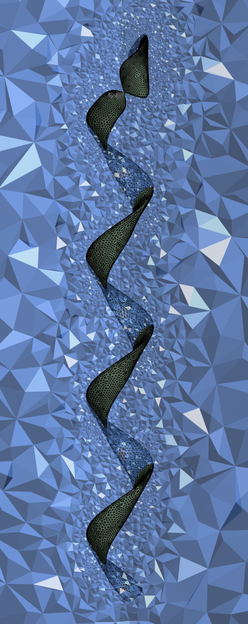}}	
	\caption{Detail of the meshes generated during the automatic mesh adaptive process for the case with initial mesh shown in~figure~\ref{fig:swimmerInitialMeshes}(b).}
	\label{fig:swimmerMeshesPi4}
\end{figure}

The automatic mesh adaptive process is launched with a desired relative error of $\varepsilon = 5 \times 10^{-2}$. Convergence of the adaptivity procedure is achieved in four iterations for the three geometries considered. For the geometry corresponding to $\gamma=\pi/4$, figure~\ref{fig:swimmerMeshesPi4} shows the second, third and fourth mesh obtained during the automatic adaptive process.
The meshes in figure~\ref{fig:swimmerMeshesPi4} have 193,159, 463,342, and 1,101,623 tetrahedrons respectively. As it can be observed, the refinement introduced by the automatic mesh adaptivity process concentrates the cells in the regions where the flow is more complex. The size of the corresponding global system to be solved to compute the velocity on the cell faces and the mean pressure on each cell is 1,338,781, 3,216,565 and 7,662,725, respectively.

\begin{figure}[!bt]
	\centering
	\subfigure[$\gamma=0$]     {\includegraphics[width=0.2\textwidth]{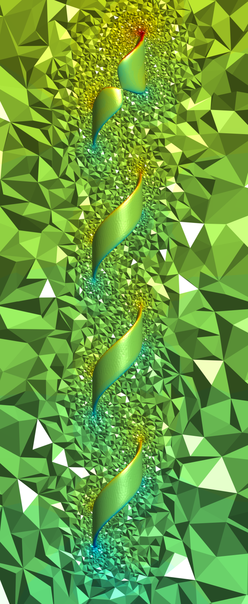}}	
	\hspace{1cm}
	\subfigure[$\gamma=\pi/4$] {\includegraphics[width=0.2\textwidth]{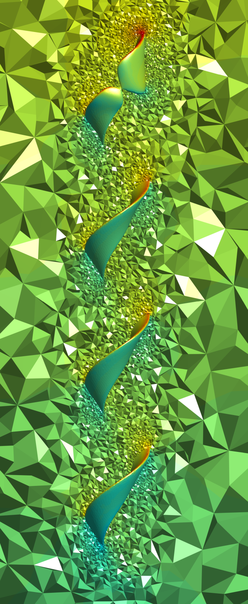}}	
	\hspace{1cm}
	\subfigure[$\gamma=\pi/2$] {\includegraphics[width=0.2\textwidth]{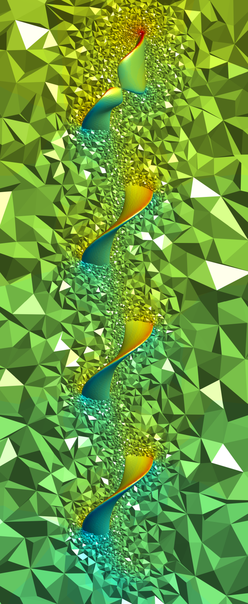}}	
	\caption{Pressure distribution on the final mesh for the microswimmers of figure~\ref{fig:swimmerGeometry}.}
	\label{fig:swimmerPressureFinalMesh}
\end{figure}

The pressure distribution on the last mesh is displayed in figure~\ref{fig:swimmerPressureFinalMesh}. 
To offer a visual comparison between the three cases, the colour scale is adjusted in the three simulations to be in between $\SI{-26.60}{mPa}$ and $\SI{33.24}{mPa}$. The results clearly show the significant variation in the pressure distribution as the geometric configuration, controlled by the parameter $\gamma$, is changed. As the value of $\gamma$ increases the reference area of the body of the swimmer increases, generating higher pressure values in the body for the geometry with $\gamma=\pi/2$ when compared to the geometry with $\gamma=0$. In all cases, the maximum pressure is observed on the head of the swimmer and the magnitude of the maximum pressure is similar in all three configurations.

\begin{figure}[!bt]
	\centering
	\subfigure[$\gamma=0$]     {\includegraphics[width=0.14\textwidth]{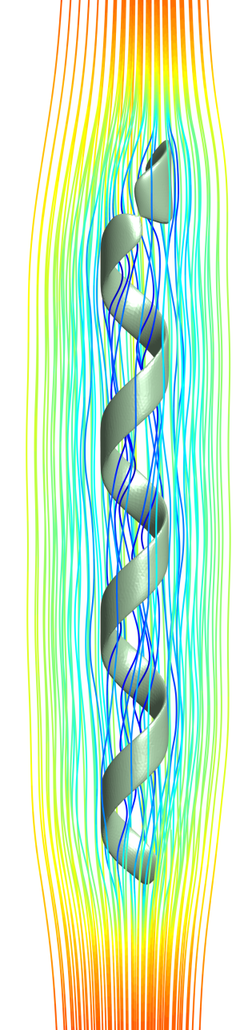}}	
	\hspace{1cm}
	\subfigure[$\gamma=\pi/4$] {\includegraphics[width=0.14\textwidth]{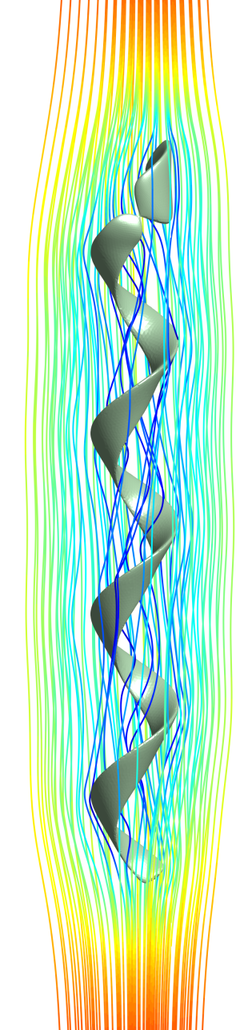}}	
	\hspace{1cm}
	\subfigure[$\gamma=\pi/2$] {\includegraphics[width=0.14\textwidth]{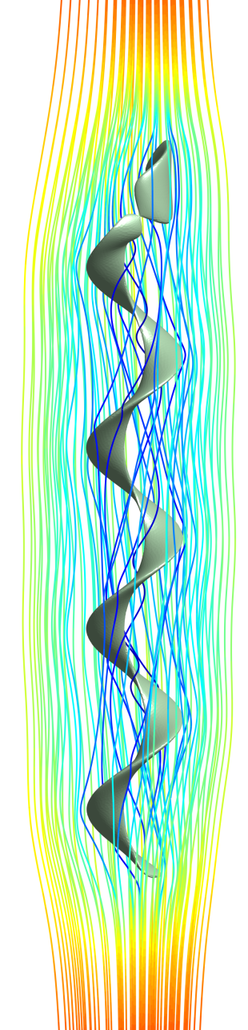}}	
	\caption{Velocity streamlines on the final mesh for the microswimmers of figure~\ref{fig:swimmerGeometry}.}
	\label{fig:swimmerStreamFinalMesh}
\end{figure}

Finally, figure~\ref{fig:swimmerStreamFinalMesh} shows the streamlines coloured with the magnitude of the velocity, with a minimum value of $\SI{0}{\mu m/s}$ and a maximum value of $\SI{4.92}{\mu m/s}$.
The complexity of the flow around the microswimmers can be appreciated as well as the influence of the geometric parameter $\gamma$ on the flow features.

\section{Concluding remarks}
\label{sc:Conclusion}

This paper proposed a formulation of the second-order FCFV method for elliptic problems suitable for application in general meshes of triangular and quadrilateral cells in two dimensions and tetrahedral, hexahedral, prismatic and pyramidal cells in three dimensions.
The computational cost of the resulting problem is comparable to the one of the original first-order FCFV method in terms of number of operations, since in both cases the global unknown is approximated with constant functions on the mesh faces.
As in the original FCFV method, optimal first-order convergence of the gradient of the solution is achieved without the need to perform a reconstruction procedure.
In addition, when CPU time is compared, the proposed method guarantees an improved approximation of the primal variable, which is now second-order accurate, of almost two orders of magnitude in three dimensional problems.

The proposed approach also inherits the robustness of the original first-order FCFV method in the incompressible limit. In addition, the proposed method is insensitive to the choice of the type of cells utilised in the discretisation, to their distortion and stretching. More precisely, successful simulations with hybrid meshes were also presented, paving the path towards the application of the discussed methodology to more complex problems requiring boundary layer meshes.

Finally, an automatic mesh adaptivity strategy is devised by means of a local error indicator obtained via the solution of one extra local problem per cell. The resulting cost of this strategy is thus limited, whereas its advantages are shown in presence of complex geometries and localised phenomenons as, starting from a coarse discretisation, the method is able to automatically construct a set of meshes to achieve a user defined target accuracy.

Extensive numerical simulations in two and three dimensions are presented to validate the proposed FCFV methodology and the mesh adaptivity procedure. Moreover, a three dimensional incompressible Stokes problem featuring geometries of interest in microfluidics applications is presented, showing the potential of the method, enhanced by the automatic mesh adaptivity strategy, to treat complex large scale flow problems.

\section*{Acknowledgements}
This work was partially supported by the European Union's Horizon 2020 research and innovation programme under the Marie Sk\l odowska-Curie Actions (Grant number: 764636) and by the Spanish Ministry of Economy and Competitiveness (Grant agreement No. DPI2017-85139-C2-2-R). The first author is also grateful for the financial support provided by the Generalitat de Catalunya (Grant agreement No. 2017-SGR-1278). The second author also acknowledges the support of the Engineering and Physical Sciences Research Council (EP/P033997/1).

\bibliographystyle{elsarticle-num}
\bibliography{Ref-HDG}

\end{document}